\documentclass[english, hidelinks]{article}
\usepackage{geometry}
\usepackage{setspace}
\usepackage{authblk}
\usepackage[utf8]{inputenc}
\usepackage{color,soul,yhmath}
\usepackage{xcolor}
\usepackage{bm}
\usepackage{graphicx}
\graphicspath{ {.} }
\usepackage{hyperref}
\usepackage{natbib}
\usepackage{bbm}
 \usepackage{xr}
\usepackage{amsmath}
\usepackage{amssymb}
\usepackage{booktabs}
\newcommand{\ra}[1]{\renewcommand{\arraystretch}{#1}}

\usepackage[makeroom]{cancel}

\DeclareMathOperator*{\argmin}{arg\,min}
\hypersetup{
    colorlinks=false,
    linkcolor=black,
    filecolor=magenta,
    urlcolor=cyan,
    }
    
\newtheorem{theorem}{Theorem}

\newtheorem{corollary}[theorem]{Corollary}

\newtheorem{lemma}{Lemma}

\newtheorem{remark}{Remark}

\numberwithin{equation}{section}
\def\c#1{\mathcal{#1}}
\def\b#1{\mathbbm{#1}}
\def\bf#1{\mathbf{#1}}

\def\red#1{\textcolor{red}{#1}}

\def\wt#1{\widetilde{#1}}
\def\wh#1{\widehat{#1}}

\def\vecn#1{\textnormal{\textbf{vec}}\big(#1\big)}

\def\cN{\mathcal{N}}

\def\cQ{\mathcal{Q}}
\def\cI{\mathcal{I}}

\def\A{\mathbf{A}}
\def\B{\mathbf{B}}
\def\C{\mathbf{C}}

\def\Q{\mathbf{Q}}

\def\D{\mathbf{D}}
\def\Y{\mathbf{Y}}
\def\M{\mathbf{M}}
\def\F{\mathbf{F}}
\def\E{\mathbf{E}}
\def\H{\mathbf{H}}
\def\v{\mathbf{v}}
\def\g{\mathbf{g}}
\def\R{\mathbf{R}}
\def\I{\mathbf{I}}
\def\L{\mathbf{L}}

\def\X{\mathbf{X}}
\def\Z{\mathbf{Z}}
\def\0{\mathbf{0}}
\def\1{\mathbf{1}}

\def\diag{\text{diag}}

\def\bSigma{\boldsymbol{\Sigma}}
\def\bgamma{\boldsymbol{\gamma}}
\def\balpha{\boldsymbol{\alpha}}
\def\bbeta{\boldsymbol{\beta}}
\def\boldeta{\boldsymbol{\eta}}
\def\btheta{\boldsymbol{\theta}}
\def\bGamma{\boldsymbol{\Gamma}}
\def\bepsilon{\boldsymbol{\epsilon}}
\def\bmu{\boldsymbol{\mu}}
\def\var{\textnormal{Var}}
\def\tr{\text{tr}}
\DeclareMathOperator*{\plim}{plim}

\def\h{\mathbf{h}}

\def\ddA{\ddot{\A}}
\def\ddF{\ddot{\F}}

\newcommand{\norm}[1]{\big\| #1 \big\|}

\begin{document}
\setlength{\parindent}{18pt}
	\doublespacing
	\begin{titlepage}
		
		\title{Matrix-valued Factor Model with Time-varying Main Effects}
                \author{Clifford Lam\thanks{Clifford Lam is Professor, Department of Statistics, London School of Economics. Email: C.Lam2@lse.ac.uk}}
		        \author{Zetai Cen\thanks{Zetai Cen is PhD student, Department of Statistics, London School of Economics. Email: Z.Cen@lse.ac.uk}}
		
		\affil{Department of Statistics, London School of Economics and Political Science}
		
		\date{}
		
		\maketitle

\begin{abstract}
We introduce the matrix-valued time-varying Main Effects Factor Model (MEFM). MEFM is a generalization to the traditional matrix-valued factor model (FM). We give rigorous definitions of MEFM and its identifications, and propose estimators for the time-varying grand mean, row and column main effects, and the row and column factor loading matrices for the common component. Rates of convergence for different estimators are spelt out, with asymptotic normality shown. The core rank estimator for the common component is also proposed, with consistency of the estimators presented. We propose a test for testing if FM is sufficient against the alternative that MEFM is necessary, and demonstrate the power of such a test in various simulation settings. We also demonstrate numerically the accuracy of our estimators in extended simulation experiments. A set of NYC Taxi traffic data is analysed and our test suggests that MEFM is indeed necessary for analysing the data against a traditional FM.
\end{abstract}
		
		\bigskip
		\bigskip

		\noindent
		{\sl Key words and phrases:} Large-scale dependent data; Time-varying row and column effects; MEFM and FM interchange;  sufficiency of FM over MEFM; Tucker decomposition.

\noindent

	\end{titlepage}
	
	\setcounter{page}{2}

\maketitle


\newpage
\section{Introduction}\label{sec:introduction}
Matrix-valued  time series factor models, a generalization of vector time series factor models \citep{Bai2003,StockWatson2002}, have been utilized a lot for dimension reduction and prediction in recent years in fields such as finance, economics, medical science and meteorology, to name but a few. This is a subject still in its infancy, but important earlier theoretical and methodological developments include \cite{WangLiuChen2019}, \cite{Chenetal2020}, \cite{ChenFan2023} and \cite{Heetal2024}, which are all on matrix-valued factor models using the Tucker decomposition for the common component. \cite{Changetal2023} uses the CP decomposition for the common component, while \cite{Guan2023} considers Tucker decomposition of the common component but taking in covariates in the loadings. Beyond factor modelling, \cite{Chenetal2021}, \cite{WuBi023} and \cite{Zhang2024} propose autoregressive and moving average models for matrix-valued time series data. See \cite{Tsay2023} for a comprehensive review of matrix-valued time series analysis.

With Tucker decomposition, a matrix-valued time series factor model (FM) can be written as
\begin{equation}\label{eqn: FM}
\Y_t = \bmu + \R\F_t\C' + \E_t,
\end{equation}
where $\Y_t \in \mathbb{R}^{p\times q}$ is the observed matrix at time $t$, $\mu \in \mathbb{R}^{p\times q}$ is the mean matrix, $\R \in \mathbb{R}^{p\times k_r}$ and $\C \in \mathbb{R}^{q \times k_c}$ are the row and column factor loading matrices respectively, $\F_t \in \mathbb{R}^{k_r\times k_c}$ is the core factor matrix at time $t$, and finally $\E_t \in \mathbb{R}^{p\times q}$ is the noise matrix at time $t$. If we set $\R := (\balpha_{p\times r}, \wt\R_{p\times (k_r - r - \ell)}, \1_{p\times\ell})$, $\C:= (\1_{q\times r}, \wt\C_{q\times(k_c-r-\ell)}, \bgamma_{q\times\ell})$ and $\F_t := \diag((\g_t)_{r\times r}, (\wt\F_t)_{(k_r-r-\ell)\times(k_c-r-\ell)}, (\h_t)_{\ell\times\ell})$  \citep{Heetal2023},
where $\1_{m\times n}$ is a matrix of ones of size $m\times n$,
then (\ref{eqn: FM}) becomes
\begin{equation}\label{eqn: FM_two_way}
\Y_t = \bmu + \balpha\g_t\1_{r\times q} + \1_{p\times \ell}\h_t\bgamma' + \wt\R\wt\F_t\wt\C' + \E_t.
\end{equation}
If the rows of $\Y_t$ represent different countries and the columns represent different economic indicators, then since the $j$th row of $\balpha\g_t\1_{r\times q}$ is $\balpha_{j\cdot}\g_t\1_{r\times q}$, where $\balpha_{j\cdot}$ is the $j$th row of $\balpha$,
it means that each element in the $j$th row is the same, with value $\balpha_{j\cdot}\g_t\1_r$. Hence we can argue that $\g_t$ represents common global factors affecting all countries, although each country loads differently on $\g_t$. Similarly, $\h_t$ represents latent economic states across different economic indicators. The term $\wt\R\wt\F_t\wt\C'$ can be viewed as an interaction term., while $\balpha\g_t\1_{r\times q}$ the country main effects, and $\1_{p\times\ell}\h_t\bgamma'$ the economic states' main effects.

Three problems arise upon inspecting (\ref{eqn: FM}) and (\ref{eqn: FM_two_way}) however. Firstly, for (\ref{eqn: FM}) to transform to (\ref{eqn: FM_two_way}), $\R$ and $\C$ are both of reduced rank. In the literature for model (\ref{eqn: FM}), we always need $\R$ and $\C$ to be of full rank at least asymptotically (see for example, Assumption (B2) in \cite{Heetal2023} or Equation (8) in \cite{ChenFan2023}) for estimation purpose.

Secondly, model (\ref{eqn: FM_two_way}) is not general enough, unless $r$ and $\ell$ can be large.
For example, if $r$ is small, each country is driven only by few global common factors affecting all countries, on top of the factors in $\wt\F_t$. This will not be a problem, if not for the fact that there can be latent common factors that only drive a small group of countries/economic indicators. For instance, there can be a few small European countries which do not share global common factors with the majority of European countries, but with other middle-Eastern countries. Such ``grouping'' of countries usually comes with their corresponding groups of unique factors. These unique factors become ``weak'' country effects, shared only among ``small'' number of countries, essentially inflating the value of $r$ while inducing a sparse $\balpha$.
Constraint factor modelling in \cite{Chenetal2020} can certainly help, but we do not always know the exact group of countries which share latent common factors.

The final problem is related to the second one. The inability of (\ref{eqn: FM_two_way}) to accommodate ``weak'' country/economic states effects originates from the fact that the common component in (\ref{eqn: FM}), $\wt\R\wt\F_t\wt\C'$, contains only pervasive factors, which is essentially assumed across all past works in factor models for matrix-valued time series. In a tensor setting however, \cite{CenLam2024} and \cite{ChenLam2024} have both allowed weak factors in the common component in the factor model.

One way to generalize (\ref{eqn: FM_two_way}) to address all aforementioned problems is to note that
\[\balpha\g_t\1_{r\times q} = (\balpha\g_t\1_r)\1_q' =: \balpha_t\1_q',\;\;\; \1_{p\times \ell}\h_t\bgamma' = \1_p(\bgamma\h_t'\1_\ell)' =: \1_p\bbeta_t',\]
where $\balpha_t$ and $\bbeta_t$ are the time-varying row and column main effects respectively. If we are able to estimate the two vectors $\balpha_t$ and $\bbeta_t$ without any low-rank constraints as in the equation above, then the second problem is naturally solved. Formally allowing for weak factors in the row and column loading matrices, like those in \cite{LamYao2012} for a vector factor model, solves the third problem. Finally, with these problems solved, we can go back to assuming full rank row and column factor loading matrices (see Assumption (L1) in Section \ref{sec: assumption} in this paper) to solve the first problem.

In this paper, we contribute to the literature in several important ways. Firstly, we generalize model (\ref{eqn: FM_two_way}) to (\ref{eqn: model}) which is the time-varying main effects factor model (MEFM), incorporating all relaxations described in the previous paragraph. Secondly, we provide estimation and inference methods and the corresponding theoretical guarantees, on top of a separate ratio-based method for identifying the core rank of $\F_t$, with consistency proved . Third and perhaps the most important of all, we provide a statistical test to test if FM in (\ref{eqn: FM}), with $\R$ and $\C$ both of full rank, is sufficient against the more general MEFM in (\ref{eqn: model}). A rejected null hypothesis of FM being sufficient for the data then means that there are row and/or column main effects that is not of a low rank structure like those in (\ref{eqn: FM_two_way}), essentially pointing to the existence of ``weak'' main effects.

The rest of the paper is organised as follows. Section 2 introduces the notations used in this paper. Section \ref{sec: model_estimation} introduces MEFM formally, laying down important identification conditions and estimation methodologies for all the components in the model. Section \ref{sec: assumption_theories} presents the assumptions for MEFM and the consistency and asymptotic normality results for its estimators. The test for FM versus MEFM is detailed in Section \ref{subsec: testFMvsMEFM}, while the core rank estimator for $\F_t$ is presented in Section \ref{sec: estimation_rank}. Finally, Section \ref{sec: numerical_results} presents our extensive simulation results and details the NYC Taxi traffic data analysis, pinpointing the presence of weak hourly main effects in the data. Our method is available in the R package \texttt{MEFM}, with instruction in its reference manual on CRAN. All proofs of the theorems are relegated to the supplementary materials of this paper.

\section{Notations}\label{sec:notations}
Throughout this paper, we use the lower-case letter, bold lower-case letter and bold capital letter, i.e. $a,\bf{a},\bf{A}$, to denote a scalar, a vector and a matrix respectively.
We also use $a_i, A_{ij}, \bf{A}_{i\cdot}, \bf{A}_{\cdot i}$ to denote, respectively, the $i$-th element of $\bf{a}$, the $(i,j)$-th element of $\bf{A}$, the $i$-th row vector (as a column vector) of $\bf{A}$, and the $i$-th column vector of $\bf{A}$. We use $\circ$ to denote the Hadamard product. We use $a\asymp b$ to denote $a=O(b)$ and $b=O(a)$, while $a\asymp_P b$ to denote $a=O_P(b)$ and $b = O_P(a)$.

Given a positive integer $m$, define $[m]:=\{1,2,\dots,m\}$.  The vector $\1_m$ denotes a vector of ones of length $m$. The $i$-th largest eigenvalue of a matrix $\A$ is denoted by $\lambda_i(\bf{A})$. The notation $\bf{A}\succcurlyeq 0$ (resp. $\bf{A} \succ 0$) means that $\bf{A}$ is positive semi-definite (resp. positive definite). We use $\bf{A}'$ to denote the transpose of $\bf{A}$, and $\diag(\A)$ to denote a diagonal matrix with the diagonal elements of $\A$, while $\diag(\{a_1,\dots,a_n\})$ represents the diagonal matrix with $\{a_1,\dots,a_n\}$ on the diagonal.

Some norm notations.
For a given set, we denote by $|\cdot|$ its cardinality. We use $\|\bf{\cdot}\|$ to denote the spectral norm of a matrix or the $L_2$ norm of a vector, and $\|\bf{\cdot}\|_F$ to denote the Frobenius norm of a matrix. We use $\|\cdot\|_{\max}$ to denote the maximum absolute value of the elements in a vector or a matrix. The notations $\|\cdot\|_1$ and $\|\cdot\|_{\infty}$ denote the $L_1$ and $L_{\infty}$-norm of a matrix respectively, defined by $\|\A\|_{1} := \max_{j}\sum_{i}|(\A)_{ij}|$ and $\|\A\|_{\infty} := \max_{i}\sum_{j}|(\A)_{ij}|$. WLOG, we always assume the eigenvalues of a matrix are arranged by descending orders, and so are their corresponding eigenvectors.

\section{Model and Estimation}\label{sec: model_estimation}

\subsection{Main effect matrix factor model}\label{sec: model}

We propose the time-varying \textbf{M}ain \textbf{E}ffect matrix \textbf{F}actor \textbf{M}odel (MEFM) such that for $t\in[T]$,
\begin{equation}
\label{eqn: model}
    \Y_t = \mu_t \1_p\1_q' + \balpha_t\1_q' + \1_p\bbeta_t'  +\C_t +\E_t ,
\end{equation}
where $\Y_t$ is a $p\times q$ observed matrix at time $t$, $\mu_t$ is a scalar representing the grand mean of $\Y_t$, $\balpha_t\in\b{R}^p$ and $\bbeta_t\in\b{R}^q$ are the row and column main effects at time $t$, respectively. The common component $\C_t:= \A_r\F_t\A_c'$ is latent, where $\F_t\in\b{R}^{k_r\times k_c}$ is the core factor series with unknown number of factors $k_r$ and $k_c$, and $\A_r$ and $\A_c$ are the row and column factor loading matrices, with size $p\times k_r$ and $q\times k_c$, respectively. Lastly, $\E_t$ is the idiosyncratic noise series with the same dimension as $\Y_t$.

Unlike FM in (\ref{eqn: FM_two_way}), the main effects $\b\alpha_t$ and $\bbeta_t$ in MEFM are not restricted to be of low rank, which significantly improves the flexibility of FM, and allows for a test of FM in (\ref{eqn: FM}) in the end.
In fact, setting concatenated matrices $\ddA_r=(\1_p, \I_p, \A_r, \1_p)$ and $\ddA_c=(\1_q, \1_q, \A_c, \I_q)$, block diagonal matrix $\ddF_t=\diag\{\mu_t, \balpha_t, \F_t, \bbeta_t'\}$, then we can read (\ref{eqn: model}) as
\[
\Y_t = \ddA_r \ddF_t \ddA_c' + \E_t .
\]
However, we observe that the dimension of the factor series is now $(2+p+k_r) \times (2+q+k_c)$, and hence there is not much dimension reduction for $\Y_t$, and both $\ddA_r$ and $\ddA_c$ have no full column ranks. This observation suggests again that MEFM is more general than FM, and numerical results in Section \ref{sec: numerical_results} actually show that even an approximate estimation by FM in general comes at a cost of using very large number of factors.

Given the above motivation of MEFM, we point out that the form of MEFM can be obtained by FM in general, see Remark \ref{remark: FM_as_MEFM} for details. For generality purpose, $\Y_t$ can have non-zero mean but we can always demean the data as the sample mean is not our main parameter of interest. The right hand side of (\ref{eqn: model}) is entirely latent and hence we propose Assumption (IC1) below to identify the grand mean and the row and column effects.

\begin{itemize}
  \item[(IC1)] (Identification)
\textit{For any $t\in[T]$, we assume $\1_p'\balpha_t =\1_q'\bbeta_t =0$, $\1_p' \A_r = \0$ and $\1_q' \A_c =\0$.
}
\end{itemize}

However, we require further identification between the factors and the factor loading matrices. To do this, we normalise the loading matrices to $\Q_r= \A_r\Z_r^{-1/2}$ and $\Q_c= \A_c\Z_c^{-1/2}$, where $\Z_r = \diag(\A_r' \A_r)$ and $\Z_c = \diag(\A_c' \A_c)$, measuring the sparsity of each column of loading matrices and hence the factor strength. For example, $\F_t$ pervasive in the $j$-th row will have the $j$-th column of $\A_r$ dense and hence the $j$-th diagonal entry of $\Z_r$ will be of order $p$. For technical details, see Assumption (L1). We leave the identification to Section \ref{sec: assumption}. Assumption (IC1) also facilitates the estimation of $\mu_t$, $\balpha_t$ and $\bbeta_t$, and we discuss in the next section how to estimate the grand mean, the row and column effects, and the row and column factor loading matrices in (\ref{eqn: model}).

\subsection{Estimation of the main effects and factor components}\label{sec: estimation_fac}
The factor structure is hidden in $\Y_t$ and we need to estimate the time-varying grand mean and main effects first. For the grand mean, right-multiplying by $\1_q$ and left-multiplying by $\1_p'$ on both sides of (\ref{eqn: model}) results in $\1_p'\Y_t\1_q=pq \mu_t + \1_p'\E_t\1_q$ by Assumption (IC1). Hence for each $t\in[T]$, we obtain the moment estimator for the time-varying grand mean as
\begin{equation}
\label{eqn: mu_estimator}
    \wh\mu_t := \1_p'\Y_t\1_q /pq .
\end{equation}
Also, right-multiplying by $\1_q$ and left-multiplying by $\1_p'$ lead respectively to $\Y_t\1_q = q \mu_t \1_p + q\balpha_t + \E_t\1_q$ and $\1_p'\Y_t = p \mu_t \1_q' + p\bbeta_t'  + \1_p'\E_t$. Therefore, we obtain the time-varying row and column effect estimators as
\begin{equation}
\label{eqn: alpha_beta_estimator}
    \wh\balpha_t := q^{-1}\Y_t\1_q - \wh\mu_t \1_p
    , \;\;\;
    \wh\bbeta_t':=p^{-1}\1_p'\Y_t - \wh\mu_t \1_q' .
\end{equation}
Finally, we introduce the following to estimate the factor structure,
\begin{align}
    \wh\L_t := \Y_t - \wh\mu_t \1_p\1_q' - \wh\balpha_t\1_q' - \1_p\wh\bbeta_t'
    &= \Y_t + (pq)^{-1} \1_p'\Y_t\1_q \1_p\1_q' - q^{-1}\Y_t\1_q \1_q' -  p^{-1}\1_p\1_p'\Y_t \notag\\
    &= \M_p\Y_t\M_q, \label{eqn: Ct_estimator}
\end{align}
where $\M_m := \I_m - m^{-1}\1_m\1_m'$  for any positive integer $m$.
From the above, $\wh\L_t\wh\L_t'$ admits $\1_p$ in its null space, and $\wh\L_t'\wh\L_t$ admits $\1_q$ instead.
The factor structure can hence be estimated, with $\wh\Q_r$ constructed as the eigenvectors of $T^{-1}\sum_{t=1}^T\wh\L_t\wh\L_t'$ corresponding to the first $k_r$ largest eigenvalues, and $\wh\Q_c$ the eigenvectors of $T^{-1} \sum_{t=1}^T\wh\L_t'\wh\L_t$ corresponding to the first $k_c$ largest eigenvalues.

We can then estimate the factor time series $\F_{Z,t} = \Z_r^{1/2}\F_t\Z_c^{1/2}$, and the common component $\C_t$, respectively as
\begin{align}
  \wh{\F}_{Z,t} := \wh{\Q}_r'\wh\L_t\wh\Q_c = \wh\Q_r'\Y_t\wh\Q_c, \;\;\; \wh{\C}_t := \wh\Q_r\wh\F_{Z,t}\wh\Q_c' = \wh\Q_r\wh\Q_r'\Y_t\wh\Q_c\wh\Q_c'. \label{eqn: F_Zt_C_t_estimators}
\end{align}
Finally, the residuals $\E_t$ is estimated by
\begin{equation}\label{eqn: Ehat_t}
\wh\E_t := \wh\L_t - \wh\C_t.
\end{equation}

\begin{remark}\label{remark: FM_as_MEFM}
Suppose we have a traditional matrix-valued factor model such that $\Acute{\Y}_t = \Acute{\C}_t +\Acute{\E}_t$ where $\Acute{\Y}_t$, $\Acute{\C}_t$, and $\Acute{\E}_t$ are $p\times q$ matrices representing the observation, common component, and noise, respectively. Suppose also $\Acute{\C}_t = \A_r\F_t\A_c'$. Then we can construct
\begin{equation*}
    \Acute{\mu}_t := (pq)^{-1} \1_p'\Acute{\C}_t\1_q,
    \;\;\;
    \Acute{\balpha}_t := q^{-1}\Acute{\C}_t\1_q
    - \Acute{\mu}_t\1_p = q^{-1}\M_p\Acute{\C}_t\1_q,
    \;\;\;
    \Acute{\bbeta}_t := p^{-1}\Acute{\C}_t'\1_p
    - \Acute{\mu}_t\1_q = p^{-1}\M_q\Acute{\C}_t'\1_p.
\end{equation*}
Hence we can express FM in the following MEFM form satisfying (IC1):
\[
\Acute{\Y}_t=\Acute{\mu}_t \1_p\1_q' + \Acute{\balpha}_t \1_q' + \1_p\Acute{\bbeta}_t'
+(\Acute{\C}_t - \Acute{\mu}_t \1_p\1_q' - \Acute{\balpha}_t\1_q' - \1_p\Acute{\bbeta}_t' )
+\Acute{\E}_t,
\]
where
\[\Acute{\C}_t - \Acute{\mu}_t \1_p\1_q' - \Acute{\balpha}_t\1_q' - \1_p\Acute{\bbeta}_t'  = (\M_p\A_r)\F_t(\M_q\A_c)',\]
is the common component. Since $\M_m\1_m = \0$, it is easy to see that
\[\1_p'(\M_p\A_r) = \0, \;\;\; \1_q'(\M_q\A_c) = \0.\]
It is also easy to verify that $\1_p'\Acute{\balpha}_t=\1_q'\Acute{\bbeta_t} = 0$.
Hence a traditional matrix-valued factor model can be expressed as MEFM in (\ref{eqn: model}) to satisfy (IC1).
\end{remark}

\section{Assumptions and Theoretical Results}\label{sec: assumption_theories}

\subsection{Assumptions}\label{sec: assumption}
A set of assumptions on the factor structure is imposed below, and in particular, we allow factors to have different strengths, as in \cite{LamYao2012} and \cite{ChenLam2024}.

\begin{itemize}
    \item[(M1)] (Alpha mixing)
\textit{The elements in $\vecn{\F_t}$ and $\vecn{\E_t}$ are $\alpha$-mixing. A vector process $\{\bf{x}_t: t=0,\pm 1, \pm2,\dots\}$ is $\alpha$-mixing if, for some $\gamma>2$, the mixing coefficients satisfy the condition that
\begin{equation*}
    \sum_{h=1}^\infty
    \alpha(h)^{1-2/\gamma}<\infty,
\end{equation*}
where $\alpha(h)=\sup_\tau\sup_{A\in\c{H}_{-\infty}^\tau,
B\in\c{H}_{\tau+h}^\infty}|\b{P}(A\cap B)-\b{P}(A)\b{P}(B)|$ and $\c{H}_\tau^s$ is the $\sigma$-field generated by $\{\bf{x}_t: \tau\leq t\leq s\}$.
}
\end{itemize}
\begin{itemize}
    \item[(F1)] (Time Series in $\F_t$)
\textit{
There is $\X_{f,t}$ the same dimension as $\F_t$, such that $\F_t = \sum_{w\geq 0}a_{f,w}\X_{f,t-w}$. The time series $\{\X_{f,t}\}$ has i.i.d. elements with mean $0$ and variance $1$, with uniformly bounded fourth order moments. The coefficients $a_{f,w}$ are such that $\sum_{w\geq 0}a_{f,w}^2=1$ and $\sum_{w\geq 0}|a_{f,w}|\leq c$ for some constant $c$.
}
\item[(L1)] (Factor strength)
\textit{We assume that $\A_r$ and $\A_c$ are of full rank and independent of factors and errors series. Furthermore, as $p,q\to\infty$,
\begin{equation}
\label{eqn: L1}
    \Z_r^{-1/2}\A_r'\A_r\Z_r^{-1/2}
    \to \bSigma_{A,r},
    \;\;\;
    \Z_c^{-1/2}\A_c'\A_c\Z_c^{-1/2}
    \to \bSigma_{A,c},
\end{equation}
where $\Z_r=\textnormal{diag}(\A_r'\A_r)$, $\Z_c=\textnormal{diag}(\A_c'\A_c)$, and both $\bSigma_{A,r}$ and $\bSigma_{A,c}$ are positive definite with all eigenvalues bounded away from 0 and infinity. We assume $(\Z_r)_{jj}\asymp p^{\delta_{r,j}}$ for $j\in[k_r]$ and $1/2<\delta_{r,k_r}\leq \dots\leq \delta_{r,2}\leq \delta_{r,1}\leq 1$. Similarly, we assume $(\Z_c)_{jj}\asymp p^{\delta_{c,j}}$ for $j\in[k_c]$, with $1/2< \delta_{c,k_c}\leq \dots\leq \delta_{c,2}\leq \delta_{c,1}\leq 1$.
}
\end{itemize}
With Assumption (L1), we can denote $\Q_r:=\A_r\Z_r^{-1/2}$ and $\Q_c:=\A_c\Z_c^{-1/2}$. Hence  $\Q_r'\Q_r\to\bSigma_{A,r}$ and $\Q_c'\Q_c\to\bSigma_{A,c}$.

\begin{itemize}
  \item[(E1)] (Decomposition of $\E_t$)
{\em We assume that
\begin{equation}
    \label{eqn: E1}
    \E_t = \A_{e,r}\F_{e,t}\A_{e,c}' + \bSigma_{\epsilon}\circ \bepsilon_t ,
\end{equation}
where $\F_{e,t}$ is a matrix of size $k_{e,r}\times k_{e,c}$, containing independent elements with mean $0$ and variance $1$. The matrix $\bepsilon_t\in\b{R}^{p\times q}$ contains independent elements with mean 0 and variance 1, with $\{\bepsilon_t\}$ independent of $\{\F_{e,t}\}$. The matrix $\bSigma_{\epsilon}$ contains the standard deviations of the corresponding elements in $\bepsilon_t$, and has elements uniformly bounded away from 0 and infinity.

Moreover, $\A_{e,r}$ and $\A_{e,c}$ are (approximately) sparse matrices with sizes $p\times k_{e,r}$ and $q\times k_{e,c}$ respectively, such that $\|\A_{e,r}\|_1, \|\A_{e,c}\|_1=O(1)$, with $k_{e,r}, k_{e,c} = O(1)$.
}
\end{itemize}
\begin{itemize}
  \item[(E2)] (Time Series in $\E_t$)
{\em
There is $\X_{e,t}$ the same dimension as $\F_{e,t}$, and $\X_{\epsilon,t}$ the same dimension as $\bepsilon_t$, such that $\F_{e,t} = \sum_{w\geq 0}a_{e,w}\X_{e,t-w}$ and $\bepsilon_t = \sum_{w\geq 0}a_{\epsilon,w}\X_{\epsilon,t-w}$, with $\{\X_{e,t}\}$ and $\{\X_{\epsilon,t}\}$ independent of each other, and each time series has independent elements with mean $0$ and variance $1$ with uniformly bounded fourth order moments. Both $\{\X_{e,t}\}$ and $\{\X_{\epsilon,t}\}$ are independent of $\{\X_{f,t}\}$ from (F1).

The coefficients $a_{e,w}$ and $a_{\epsilon,w}$ are such that $\sum_{w\geq 0}a_{e,w}^2 = \sum_{w\geq 0}a_{\epsilon,w}^2 = 1$ and $\sum_{w\geq 0}|a_{e,w}|, \sum_{w\geq 0}|a_{\epsilon,w}|\leq c$ for some constant $c$.
}
\end{itemize}
\begin{itemize}
  \item[(R1)] (Rate assumptions)
{\em
We assume that,
\begin{align*}
    & T^{-1} p^{2(1-\delta_{r,k_r})}q^{1-2\delta_{c,1}} = o(1) , \;\;\;
    p^{1-2\delta_{r,k_r}}q^{2(1-\delta_{c,1})} = o(1), \\
    & T^{-1} q^{2(1-\delta_{c,k_c})}p^{1-2\delta_{r,1}} = o(1) , \;\;\;
    q^{1-2\delta_{c,k_c}}p^{2(1-\delta_{r,1})} = o(1) .
\end{align*}
}
\end{itemize}

Assumption (F1) introduces serial dependence into the factors, and (E1) and (E2) introduce both cross-sectional and temporal dependence in the noise. The factor structure depicted by (F1), (E1) and (E2) is the same as the one in \cite{CenLam2024}. Note that although Assumption (M1) also features in serial dependence, it is mainly used to construct asymptotic normality of estimators.

By Assumption (L1), we have $\A_r\F_t\A_c'= \Q_r\Z_r^{1/2}\F_t\Z_c^{1/2}\Q_c'$, so we aim to estimate $(\Q_r, \Q_c,\F_{Z,t})$ where $\F_{Z,t}:=\Z_r^{1/2}\F_t\Z_c^{1/2}$. Unlike the traditional approximate factor models which assumes all factors are pervasive, we allow factors to have different strength similar to \cite{LamYao2012} and \cite{ChenLam2024}. To be precise, a column of $\A_r$ (resp. $\A_c$) is dense (i.e., a pervasive factor) if the corresponding $\delta_{r,j}=1$ (resp. $\delta_{c,j}=1$), otherwise the column represents a weak factor as it is sparse to certain extent.

Due to the presence of potentially weak factors, we require rate conditions in Assumption (R1) for consistency to hold. If all factors are pervasive (R1) holds trivially. We point out that the first and second conditions in (R1) are exactly the same as the first and third conditions of Assumption (R1) in \cite{CenLam2024} for matrix time series.

\subsection{Identification of the model}
With Assumptions (IC1) and (L1), the model (\ref{eqn: model}) is identified according to Theorem \ref{thm: identification} below.

\begin{theorem}\label{thm: identification}
(Identification) With Assumption (IC1), each $\mu_t$, $\balpha_t$, and $\bbeta_t$ can be identified. The common component is hence identified, and if (L1) is also satisfied, the factor structure is identified up to some rotation such that $(\Q_r,\Q_c, \F_{Z,t})=(\Q_r\M_r,\Q_c\M_c, \M_r^{-1}\F_{Z,t}\M_c^{-1})$ for some invertible matrices $\M_r\in\b{R}^{k_r\times k_r}$ and $\M_c\in\b{R}^{k_c\times k_c}$.
\end{theorem}

\subsection{Rate of convergence for various estimators}\label{subsec: rate_of_convergence}
To present the consistency of the loading estimators, define
\begin{align}
    \H_r &:= T^{-1}\wh\D_r^{-1} \wh\Q_r' \Q_r \sum_{t=1}^T (\F_{Z,t} \Q_c' \Q_c \F_{Z,t}') ,
    \label{def: H_r} \\
    \H_c &:= T^{-1}\wh\D_c^{-1} \wh\Q_c' \Q_c \sum_{t=1}^T (\F_{Z,t}' \Q_r' \Q_r \F_{Z,t}) ,
    \label{def: H_c}
\end{align}
where $\wh\D_r := \wh\Q_r'(T^{-1} \sum_{t=1}^T \wh\L_t\wh\L_t')\wh\Q_r$ is the $k_r \times k_r$ diagonal matrix of eigenvalues of $T^{-1} \sum_{t=1}^T \wh\L_t\wh\L_t'$, and similarly $\wh\D_c := \wh\Q_c' (T^{-1} \sum_{t=1}^T \wh\L_t'\wh\L_t)\wh\Q_c$ is the $k_c \times k_c$ diagonal matrix of eigenvalues of $T^{-1} \sum_{t=1}^T \wh\L_t'\wh\L_t$.

\begin{theorem}\label{thm: consistency}
Under Assumptions (IC1), (M1), (F1), (L1), (E1), (E2) and (R1), we have
\begin{align}
    (\wh\mu_t - \mu_t)^2 &= O_P(p^{-1} q^{-1}) ,
    \label{eqn: consistency_mu} \\
    p^{-1} \|\wh\balpha_t - \balpha_t\|^2 &= O_P(q^{-1}) ,
    \label{eqn: consistency_alpha} \\
    q^{-1} \|\wh\bbeta_t - \bbeta_t\|^2 &= O_P(p^{-1}) ,
    \label{eqn: consistency_beta} \\
    p^{-1} \|\wh\Q_r - \Q_r\H_r' \|_F^2 &=
    O_P\Big(T^{-1}p^{1-2\delta_{r,k_r}}q^{1-2\delta_{c,1}} + p^{-2\delta_{r,k_r}}q^{2(1-\delta_{c,1})} \Big) ,
    \label{eqn: consistency_row_loading} \\
    q^{-1} \|\wh\Q_c - \Q_c\H_c' \|_F^2 &=
    O_P\Big(T^{-1}q^{1-2\delta_{c,k_c}}p^{1-2\delta_{r,1}} + q^{-2\delta_{c,k_c}}p^{2(1-\delta_{r,1})} \Big) .
    \label{eqn: consistency_column_loading}
\end{align}
\end{theorem}

From Theorem \ref{thm: consistency}, the consistency for the loading matrix estimators requires Assumption (R1). If all factors are pervasive, the (squared) convergence rates for the row (resp. column) loading matrix will be \linebreak $\max(1/(Tpq), 1/p^2)$ (resp. $\max(1/(Tpq), 1/q^2)$), which are consistent with those in \cite{ChenFan2023} after the same normalization of the loading matrices.

\begin{theorem}\label{thm: consistency_F_C}
Under the assumptions in Theorem \ref{thm: consistency}, we have the following:
\begin{itemize}
    \item [1.] The error of the estimated factor series has rate
    \begin{align*}
    \|\wh\F_{Z,t} - (\H_r^{-1})' \F_{Z,t} \H_c^{-1}\|_F^2 &=
    O_P\big( p^{1-\delta_{r,k_r}} q^{1-\delta_{c,k_c}}
    + T^{-1} p^{1+ 2\delta_{r,1} - 2\delta_{r,k_r}} q^{1-\delta_{c,1}} +  p^{1 +\delta_{r,1} -3\delta_{r,k_r}} q^{2-\delta_{c,1}} \\
    & \hspace{36pt}
    + T^{-1} q^{1+ 2\delta_{c,1} - 2\delta_{c,k_c}} p^{1-\delta_{r,1}} +  q^{1 +\delta_{c,1} -3\delta_{c,k_c}} p^{2-\delta_{r,1}}\big) .
    \end{align*}
    \item [2.] For any $t\in[T], i\in[p], j\in[q]$, the squared error of the estimated individual common component is
    \begin{align*}
    (\wh{C}_{t,ij} - C_{t,ij})^2
    &= O_P\big(p^{1-2\delta_{r,k_r}} q^{1-2\delta_{c,k_c}}
    + T^{-1} p^{1 + 2\delta_{r,1} -3\delta_{r,k_r}} q^{1 -\delta_{c,1} -\delta_{c,k_c}} + p^{1 +\delta_{r,1} -4\delta_{r,k_r}} q^{2 -\delta_{c,1} -\delta_{c,k_c}} \\
    & \hspace{36pt}
    + T^{-1} q^{1 + 2\delta_{c,1} -3\delta_{c,k_c}} p^{1 -\delta_{r,1} -\delta_{r,k_r}} + q^{1 +\delta_{c,1} -4\delta_{c,k_c}} p^{2 -\delta_{r,1} -\delta_{r,k_r}}
    \big) .
    \end{align*}
\end{itemize}
\end{theorem}
We state the above results separating from Theorem \ref{thm: consistency} since they have used some arguments from the proof of Theorem \ref{thm: asymp_loading}. If all factors are pervasive, it is clear that individual common components are consistent with rate $(pq)^{-1/2} + T^{-1/2}(q^{-1/2} + p^{-1/2}) + p^{-1} + q^{-1} = \max(1/(Tq)^{1/2}, 1/(Tp)^{1/2}, 1/p, 1/q)$.
This rate coincides with Theorem 4 of \cite{ChenFan2023} for instance.

\subsection{Asymptotic normality of estimators}\label{subsec:normality}
We present the asymptotic normality of various estimators in this section, together with the estimation of the corresponding covariance matrices for practical inferences. But before that, we need three more assumptions.

\begin{itemize}
  \item[(L2)] (Eigenvalues)
{\em
The eigenvalues of the $k_r\times k_r$ matrix $\bSigma_{A,r}\Z_r$ from Assumption (L1) are distinct, and so are those of the $k_c\times k_c$ matrix $\bSigma_{A,c}\Z_c$.
}

\item[(AD1)] \textit{There is a constant $C$ such that for any $k\in[K],j\in[d_k]$, as $p, q,T \to\infty$,
\begin{align*}
    & \sqrt{\frac{1}{Tq p^{\delta_{r,1}}}} \cdot
    \b{E}\Big\{\Big\| \H_r^\ast \sum_{i=1}^p \Q_{r,i\cdot} \sum_{t=1}^T (\C_t\E_t')_{ij}\Big\|\Big\} \geq C > 0,\\
    & \sqrt{\frac{1}{Tp q^{\delta_{c,1}}}} \cdot
    \b{E}\Big\{\Big\| \H_c^\ast \sum_{i=1}^q \Q_{c,i\cdot} \sum_{t=1}^T (\C_t'\E_t)_{ij}\Big\|\Big\} \geq C > 0,
\end{align*}
where
$\H_r^{\ast} :=  \textnormal{tr}(\A_c'\A_c))^{1/2} \cdot \D_r^{-1/2} (\bGamma_r^\ast)'\Z_r^{1/2}$ with $\D_r := \textnormal{tr}(\A_c'\A_c) \textnormal{diag}\{ \lambda_1(\A_r'\A_r), \ldots, \lambda_{k_r}(\A_r'\A_r)\}$, and $\bGamma_r^\ast$ is the eigenvector matrix of $\textnormal{tr}(\A_c'\A_c)\cdot p^{-\delta_{r,k_r}} q^{-\delta_{c,1}} \Z_r^{1/2} \bSigma_{A,r} \Z_r^{1/2}$. Similarly, we have \linebreak $\H_c^{\ast} :=  \textnormal{tr}(\A_r'\A_r))^{1/2} \cdot \D_c^{-1/2} (\bGamma_c^\ast)'\Z_c^{1/2}$, with $\D_c := \textnormal{tr}(\A_r'\A_r) \textnormal{diag}\{ \lambda_1(\A_c'\A_c), \ldots, \lambda_{k_c}(\A_c'\A_c)\}$, and $\bGamma_c^\ast$ is the eigenvector matrix of $\textnormal{tr}(\A_r'\A_r) \cdot q^{-\delta_{c,k_c}} p^{-\delta_{r,1}} \Z_c^{1/2} \bSigma_{A,c} \Z_c^{1/2}$.
}

\item[(R2)] We have
\begin{align*}
  &T^{-1} p^{1 + 2\delta_{r,1} -3\delta_{r,k_r}} q^{1 -\delta_{c,1} -\delta_{c,k_c}}, \;\;\;
  p^{1 +\delta_{r,1} -4\delta_{r,k_r}} q^{2 -\delta_{c,1} -\delta_{c,k_c}},\\
  &T^{-1} q^{1 + 2\delta_{c,1} -3\delta_{c,k_c}} p^{1 -\delta_{r,1} -\delta_{r,k_r}}, \;\;\;
  q^{1 +\delta_{c,1} -4\delta_{c,k_c}} p^{2 -\delta_{r,1} -\delta_{r,k_r}} = o(1).
\end{align*}
\end{itemize}
Assumption (AD1) appears in \cite{CenLam2024} as well. This assumption facilitates the proof of the asymptotic normality of each row of $\wh\Q_r$ and $\wh\Q_c$, by asserting that in the decomposition of $\wh\Q_r-\Q_r\H_r$ (resp. $\wh\Q_c - \Q_c\H_c$), certain terms are dominating others even in the lower bound, and hence is truly dominating rather than just having the upper bounds dominating other upper bounds as in the proofs of similar theorems in the broader literature of factor models. Assumption (R2) is needed to make sure that the estimated common component $\wh\C_t$ is consistent element-wise (see Theorem \ref{thm: consistency_F_C}). This is satisfied automatically when all factors are pervasive, for instance.

Let $\Sigma_{\epsilon,ij}$ be the $(i,j)$ entry of $\bSigma_{\epsilon}$ in Assumption (E1).
\begin{theorem}\label{thm:alpha_normality}
Let all assumptions in Theorem \ref{thm: consistency} hold. Assume also for $i\in[p]$ and $j\in[q]$,
\[ \gamma_{\alpha,i}^2 := \lim_{q\rightarrow\infty} \frac{1}{q}\sum_{j=1}^q \Sigma_{\epsilon,ij}^2, \;\;\;  \gamma_{\beta, j}^2 := \lim_{p\rightarrow \infty} \frac{1}{p}\sum_{i=1}^p\Sigma_{\epsilon,ij}^2, \;\;\; \gamma_{\mu}^2 := \lim_{p,q\rightarrow\infty}\frac{1}{pq}\sum_{i\in[p], j\in[q]}\Sigma_{\epsilon,ij}^2.\]
Then for each $t\in[T]$,
\[\sqrt{pq}(\wh\mu_t - \mu_t) \xrightarrow{\c{D}} \c{N}(0, \gamma_{\mu}^2).\]
Take a finite integer $m$, and integers $i_1<i_2<\cdots<i_m$ with $i_\ell\in[p]$. Define $\btheta_{\alpha,t} := (\alpha_{t,i_1},\ldots,\alpha_{t,i_m})'$ and similarly for $\hat{\btheta}_{\alpha,t}$, where $\alpha_{t,i}$ is the $i$-th element of $\balpha_t$. Then for a fixed $t\in[T]$,
\begin{align*}
  \sqrt{q}(\hat{\btheta}_{\alpha, t} - \btheta_{\alpha, t}) \xrightarrow{\c{D}} \c{N}(\0, \text{\em diag}(\gamma_{\alpha,i_1}^2,\ldots,\gamma_{\alpha,i_m}^2)).
\end{align*}
Similarly, take integers $j_1<\cdots<j_m$ where $j_\ell\in[q]$. Define $\btheta_{\beta,t}:= (\beta_{t,j_1},\ldots,\beta_{t,j_m})'$ and similarly for $\wh{\btheta}_{\beta,t}$, where $\beta_{t,j}$ is the $j$-th element of $\bbeta_t$. Then for a fixed $t\in[T]$,
\begin{align*}
  \sqrt{p}(\wh{\btheta}_{\beta,t} - \btheta_{\beta,t}) \xrightarrow{\c{D}} \c{N}(\0, \text{\em diag}(\gamma_{\beta,j_1}^2,\ldots,\gamma_{\beta,j_m}^2)).
\end{align*}
Moreover, for $i\in[p]$ and $j\in[q]$, if the rate for $\wh{C}_{t,ij}-C_{t,ij}$ in Theorem \ref{thm: consistency_F_C} is $o(1)$, then
\begin{align*}
  \wh\gamma_{\alpha,i}^2 := q^{-1}(\wh\E_t\wh\E_t')_{ii}, \;\;\; \wh\gamma_{\beta,j}^2 := p^{-1}(\wh\E_t'\wh\E_t)_{jj}, \;\;\; \wh\gamma_{\mu}^2 := p^{-1}\sum_{i=1}^p\wh\gamma_{\alpha,i}^2 = q^{-1}\sum_{j=1}^q\wh\gamma_{\beta,j}^2
\end{align*}
are consistent estimators for $\gamma_{\alpha,i}^2$, $\gamma_{\beta,j}^2$ and $\gamma_{\mu}^2$ respectively under Assumption (R2), so that
\begin{align*}
  \sqrt{pq}\,\wh\gamma_{\mu}^{-1}(\wh\mu_t - \mu_t) &\xrightarrow{\c{D}} \c{N}(0,1),\\
  \sqrt{q}\, \text{\em diag}(\wh\gamma_{\alpha,i_1}^{-1},\ldots,\wh\gamma_{\alpha,i_m}^{-1})(\wh\btheta_{\alpha,t} - \btheta_{\alpha,t}) &\xrightarrow{\c{D}} \c{N}(\0, \I_m),\\
  \sqrt{p}\, \text{\em diag}(\wh\gamma_{\beta,j_1}^{-1},\ldots,\wh\gamma_{\beta,j_m}^{-1})(\wh\btheta_{\beta,t} - \btheta_{\beta,t}) &\xrightarrow{\c{D}} \c{N}(\0, \I_m).
\end{align*}
\end{theorem}
Recall from Remark \ref{remark: FM_as_MEFM} that FM can be expressed in MEFM, and hence the ability to make inferences on the elements of $\balpha_t$ and $\bbeta_t$ does not facilitate a test for the necessity of MEFM over FM. For such a test, please see Section \ref{subsec: testFMvsMEFM}. Theorem \ref{thm:alpha_normality} gives us the ability to infer on the level of row and column main effects at each time point, which is important if we have target comparisons we want to make for these effects. For instance, if each row represents a country, we can easily compare the main effects at time $t$ for the first country against the average of the second and third simply by considering $\g := (1,-1/2,-1/2)'$, $\btheta_{\alpha,t} := (\alpha_{t,1}, \alpha_{t,2}, \alpha_{t,3})'$ and using Theorem \ref{thm:alpha_normality} to arrive at
\begin{align*}
\sqrt{q}\,(\g'\diag(\wh\gamma_{\alpha,1}^2, \wh\gamma_{\alpha,2}^2, \wh\gamma_{\alpha,3}^2)\g)^{-1/2}\g'(\wh\btheta_{\alpha,t} - \btheta_{\alpha,t}) \xrightarrow{\c{D}} \c{N}(0,1).
\end{align*}

\begin{theorem}\label{thm: asymp_loading}
Let all the assumptions under Theorem \ref{thm: consistency} hold, in addition to (AD1) and (L2). Suppose $k_r$ and $k_c$ are fixed and $p,q,T \to \infty$. If $Tq= o(p^{\delta_{r,1}+ \delta_{r,k_r}})$, we have
\begin{align*}
& (T p^{2\delta_{r,k_r} - \delta_{r,1}} q^{2\delta_{c,1} -1})^{1/2} \cdot (\wh\Q_{r,j\cdot} - \H_r\Q_{r,j\cdot})
\xrightarrow{\c{D}}
\cN\big(\0, T^{-1} p^{2\delta_{r,k_r} - \delta_{r,1}} q^{2\delta_{c,1} -1} \cdot \D_r^{-1} \H_r^\ast \bf{\Xi}_{r,j} (\H_r^\ast)' \D_r^{-1} \big), \\
& \text{where} \;\;\;
\bf{\Xi}_{r,j} := \plim_{p,q,T \to\infty} \textnormal{Var} \big( \sum_{i=1}^p \Q_{r,i\cdot} \sum_{t=1}^T (\C_t \E_t')_{ij} \big) .
\end{align*}
On the other hand, if $Tp= o(q^{\delta_{c,1}+ \delta_{c,k_c}})$, we have
\begin{align*}
& (T q^{2\delta_{c,k_c} - \delta_{c,1}} p^{2\delta_{r,1} -1})^{1/2} \cdot (\wh\Q_{c,j\cdot} - \H_c\Q_{c,j\cdot})
\xrightarrow{\c{D}}
\cN\big(\0, T^{-1} q^{2\delta_{c,k_c} - \delta_{c,1}} p^{2\delta_{r,1} -1} \cdot \D_c^{-1} \H_c^\ast \bf{\Xi}_{c,j} (\H_c^\ast)' \D_c^{-1} \big), \\
& \text{where} \;\;\;
\bf{\Xi}_{c,j} := \plim_{p,q,T \to\infty} \textnormal{Var} \big( \sum_{i=1}^q \Q_{c,i\cdot} \sum_{t=1}^T (\C_t' \E_t)_{ij} \big) .
\end{align*}
\end{theorem}
Theorem \ref{thm: asymp_loading} is essentially Theorem 3 of \cite{CenLam2024} when $K=2$ and $\eta=0$ (no missing values), having the same rate of convergence under potentially weak factors. Hence our MEFM estimation procedure has successfully estimated and removed all time-varying main effects and grand mean, leaving the estimation of the common component exactly the same as in FM.

\subsubsection{Estimation of covariance matrix for factor loading estimators}\label{subsec: cov_est_Qhat}
To practically use Theorem \ref{thm: asymp_loading} for inference, we need to estimate the covariance matrices for $\wh\Q_{r,j\cdot} - \H_r\Q_{r,j\cdot}$ and $\wh\Q_{c,j\cdot} - \H_c\Q_{c,j\cdot}$. We use the heteroscedasticity and autocorrelation consistent (HAC) estimators (\cite{NeweyWest1987}) based on $\{\wh\D_r, \wh\Q_r, \wh\C_t, \wh\E_t \}_{t\in[T]}$ and $\{\wh\D_c, \wh\Q_c, \wh\C_t, \wh\E_t \}_{t\in[T]}$, respectively.

For $\wh\Q_{r,j\cdot} - \H_r\Q_{r,j\cdot}$, with $\eta_r$ such that $\eta_r \to \infty$ and $\eta_r/ (T p^{2\delta_{r,k_r} -\delta_{r,1}} q^{2\delta_{c,1} -1})^{1/4} \to 0$, define an HAC estimator
\begin{align*}
    & \wh\bSigma_{r,j}^{HAC} := \D_{r,0,j} + \sum_{\nu =1}^{\eta_r} \Big(1- \frac{\nu}{1+ \eta_r}\Big) \Big(\D_{r,\nu,j} + \D_{r,\nu,j}'\Big),
    \;\;\; \text{where} \\
    & \D_{r,\nu,j} := \sum_{t=1+\nu}^T \Big\{ \sum_{i=1}^p \Big( T^{-1} \wh\D_r^{-1} \wh\Q_r' \sum_{s=1}^T \wh\C_s \wh\C_{s,i\cdot} \Big) (\wh\C_t \wh\E_t')_{ij}\Big\} \Big\{\sum_{i=1}^p \Big( T^{-1} \wh\D_r^{-1} \wh\Q_r' \sum_{s=1}^T \wh\C_s \wh\C_{s,i\cdot} \Big) (\wh\C_{t-\nu} \wh\E_{t-\nu}')_{ij} \Big\}' .
\end{align*}
Similarly for $\wh\Q_{c,j\cdot} - \H_c\Q_{c,j\cdot}$, with $\eta_c$ such that $\eta_c \to \infty$ and $\eta_c/ (T q^{2\delta_{c,k_c} -\delta_{c,1}} p^{2\delta_{r,1} -1})^{1/4} \to 0$, define
\begin{align*}
    & \wh\bSigma_{c,j}^{HAC} := \D_{c,0,j} + \sum_{\nu =1}^{\eta_c} \Big(1- \frac{\nu}{1+ \eta_c}\Big) \Big(\D_{c,\nu,j} + \D_{c,\nu,j}'\Big),
    \;\;\; \text{where} \\
    & \D_{c,\nu,j} := \sum_{t=1+\nu}^T \Big\{ \sum_{i=1}^q \Big( T^{-1} \wh\D_c^{-1} \wh\Q_c' \sum_{s=1}^T \wh\C_s' \wh\C_{s,\cdot i} \Big) (\wh\C_t' \wh\E_t)_{ij}\Big\} \Big\{ \sum_{i=1}^q \Big( T^{-1} \wh\D_c^{-1} \wh\Q_c' \sum_{s=1}^T \wh\C_s' \wh\C_{s,\cdot i} \Big) (\wh\C_{t-\nu}' \wh\E_{t-\nu} )_{ij}\Big\}' .
\end{align*}

\begin{theorem}\label{thm: HAC}
Let all the assumptions under Theorem \ref{thm: consistency} hold, in addition to (L2), (AD1) and (R2). Suppose $k_r$ and $k_c$ are fixed and $p,q,T \to \infty$. If $Tq= o(p^{\delta_{r,1}+ \delta_{r,k_r}})$, then
\begin{itemize}
    \item [1.] $\wh\D_r^{-1} \wh\bSigma_{r,j}^{HAC} \wh\D_r^{-1}$ is consistent for $\D_r^{-1} \H_r^\ast \bf{\Xi}_{r,j} (\H_r^\ast)' \D_r^{-1}$;
    \item [2.] $T \cdot \big(\wh\bSigma_{r,j}^{HAC} \big)^{-1/2} \wh\D_r (\wh\Q_{r,j\cdot} - \H_r \Q_{r,j\cdot}) \xrightarrow{\c{D}} \cN(\0, \I_{k_r})$.
\end{itemize}
On the other hand, if $Tp= o(q^{\delta_{c,1}+ \delta_{c,k_c}})$, then
\begin{itemize}
    \item [3.] $\wh\D_c^{-1} \wh\bSigma_{c,j}^{HAC} \wh\D_c^{-1}$ is consistent for $\D_c^{-1} \H_c^\ast \bf{\Xi}_{c,j} (\H_c^\ast)' \D_c^{-1}$;
    \item [4.] $T \cdot \big(\wh\bSigma_{c,j}^{HAC} \big)^{-1/2} \wh\D_c (\wh\Q_{c,j\cdot} - \H_c \Q_{c,j\cdot}) \xrightarrow{\c{D}} \cN(\0, \I_{k_c})$.
\end{itemize}
\end{theorem}

\subsection{Testing the sufficiency of FM versus MEFM}\label{subsec: testFMvsMEFM}
In the last section, we introduce how to make inferences on various parameters of MEFM. However, to test if FM is sufficient against our proposed MEFM, simple inferences on the model parameters are not enough in the face of Remark \ref{remark: FM_as_MEFM}. Formally, we want to test, for the time horizon $t\in[T]$,
\begin{align*}
  H_0: \text{FM is sufficient over } \; t\in[T] \;\;\; \longleftrightarrow \;\;\; H_1: \text{MEFM is needed over }\; t\in[T].
\end{align*}
The above problem is complicated by the fact that, in Section \ref{sec: model}, we have seen that MEFM can always be expressed as FM if we are willing to potentially consider a large number of factors. So, how ``large'' an increase in the number of factors do we consider unacceptable?

Remark \ref{remark: FM_as_MEFM} tells us that a special form of MEFM can be expressed back in FM:
\[\Y_t = \mu_t\1_p\1_q' + \balpha_t\1_q' + \1_p\bbeta_t' + \M_p\Acute{\C}_t\M_q + \E_t = \A_r\F_t\A_c' + \E_t, \;\;\; t\in[T],\]
where $\Acute{\C}_t := \A_r\F_t\A_c'$ and
\[\mu_t := (pq)^{-1}\1_p'\Acute{\C}_t\1_q, \;\;\; \balpha_t := q^{-1}\M_p\Acute{\C}_t\1_q, \;\;\;
\bbeta_t := p^{-1}\M_q\Acute{\C}_t'\1_p.\]
If $\A_r$ has rank $k_r$ satisfying Assumption (L1) and $\A_c$ has rank $k_c$, the potential rank of $\M_p\A_r$ is $k_r-1$ (when a column in $\A_r$ is parallel to $\1_p$), and that of $\M_q\A_c$ is $k_c-1$ (when a column in $\A_c$ is parallel to $\1_q$), demonstrating that FM can have an increase in the number of factors, albeit still finite.

Another special example is when both $\balpha_t$ and $\bbeta_t$ are zero, but $\mu_t\neq 0$. Then we can write MEFM as
\begin{align*}
  \Y_t &= \mu_t\1_p\1_q' + \A_r\F_t\A_c' + \E_t = (\A_r,\1_p)\left(
                                                               \begin{array}{cc}
                                                                 \F_t & \0 \\
                                                                 \0' & \mu_t/(pq) \\
                                                               \end{array}
                                                             \right)\left(
                                                                      \begin{array}{c}
                                                                        \A_c' \\
                                                                        \1_q' \\
                                                                      \end{array}
                                                                    \right) + \E_t,
\end{align*}
which is FM with loading matrices $(\A_r,\1_p)$ and $(\A_c,\1_q)$ respectively, and an increase by 1 for both the number of row and column factors.

In light of the above examples, we deem FM sufficient if and only if the number of factors in the FM is still finite and any model variables satisfy the Assumptions in Section \ref{sec: assumption}.

To be able to test $H_0$ against $H_1$, define $\check\E_t$ to be the residual matrix after a fitting of FM (a similar procedure to fitting MEFM but treating $\mu_t$, $\balpha_t$ and $\bbeta_t$ as zero), with
\[\check\E_t := \Y_t - \check\C_t, \;\;\; \text{ where } \; \check\C_t := \check\A_r\check\A_r'\Y_t\check\A_c\check\A_c',\]
with $\check\A_r$ and $\check\A_c$ the $p\times \ell_r$ and $q\times \ell_c$ eigenmatrices of $\sum_{t=1}^T\Y_t\Y_t'$ and $\sum_{t=1}^T\Y_t'\Y_t$ respectively.

\begin{theorem}\label{thm: testingFMvsMEFM}
  Let all assumptions in Theorem \ref{thm: consistency} hold, on top of (R2). Also assume that \linebreak $\wh{C}_{t,ij} - C_{t,ij} = o_P(\min(p^{-1/2},q^{-1/2}))$ in Theorem \ref{thm: consistency_F_C}.  Suppose $k_r,k_c,\ell_r$ and $\ell_c$ are all fixed and known. Then under $H_0$, for each $t\in[T]$, we have
  \begin{align*}
    &\frac{(\wh\E_t\wh\E_t')_{ii} - \sum_{j=1}^q\Sigma_{\epsilon, ij}^2}{\sqrt{\sum_{j=1}^q\var(\epsilon_{t,ij}^2)\Sigma_{\epsilon,ij}^4}}, \; \frac{(\check\E_t\check\E_t')_{ii} - \sum_{j=1}^q\Sigma_{\epsilon, ij}^2}{\sqrt{\sum_{j=1}^q\var(\epsilon_{t,ij}^2)\Sigma_{\epsilon,ij}^4}} \xrightarrow{\c{D}} Z_i
    \xrightarrow{\c{D}} N(0,1) \; \text{ for each } \; i\in[p];\\
    &\frac{(\wh\E_t'\wh\E_t)_{jj} - \sum_{i=1}^p\Sigma_{\epsilon, ij}^2}{\sqrt{\sum_{i=1}^p\var(\epsilon_{t,ij}^2)\Sigma_{\epsilon,ij}^4}}, \; \frac{(\check\E_t'\check\E_t)_{jj} - \sum_{i=1}^p\Sigma_{\epsilon, ij}^2}{\sqrt{\sum_{i=1}^p\var(\epsilon_{t,ij}^2)\Sigma_{\epsilon,ij}^4}} \xrightarrow{\c{D}} W_j \xrightarrow{\c{D}} N(0,1) \; \text{ for each } \; j\in[q],
  \end{align*}
  where $Z_h$ is independent of $Z_\ell$ and $W_h$ is independent of $W_\ell$ for $h\neq \ell$.
  The same asymptotic results hold true under $H_1$ for $(\wh\E_t\wh\E_t')_{ii}$ and $(\wh\E_t'\wh\E_t)_{jj}$ respectively for $i\in[p], j\in[q]$.
\end{theorem}
The assumption $\wh{C}_{t,ij} - C_{t,ij} = o_P(\min(p^{-1/2},q^{-1/2}))$ is satisfied, for instance, when all factors are pervasive and $T,p,q$ are of the same order.
Theorem \ref{thm: testingFMvsMEFM} tells us that for each $t\in[T]$, both
\begin{align}
  x_{\alpha,t} := \max_{i\in[p]}\wh\gamma_{\alpha,i}^2 = \max_{i\in[p]}\{q^{-1}(\wh\E_t\wh\E_t')_{ii}\}, \;\;\;
  y_{\alpha,t} := \max_{i\in[p]}\check\gamma_{\alpha,i}^2 := \max_{i\in[p]}\{q^{-1}(\check\E_t\check\E_t')_{ii}\} \label{eqn: x_alpha}
\end{align}
are distributed approximately the same for large $q$ under $H_0$, and $x_{\alpha,t}$ in particular is distributed the same no matter under $H_0$ or $H_1$. Similarly, define
\begin{align}
  x_{\beta,t} := \max_{j\in[q]}\wh\gamma_{\beta,j}^2 = \max_{j\in[q]}\{p^{-1}(\wh\E_t'\wh\E_t)_{jj}\}, \;\;\;
  y_{\beta,t} := \max_{j\in[q]}\check\gamma_{\beta,j}^2 := \max_{j\in[q]}\{p^{-1}(\check\E_t'\check\E_t)_{jj}\}, \label{eqn: x_beta}
\end{align}
which are distributed approximately the same for large $p$ under $H_0$ from Theorem \ref{thm: testingFMvsMEFM}, and $x_{\beta,t}$ in particular is distributed the same no matter under $H_0$ or $H_1$.

Define $\b{F}_{x,\alpha}$, $\b{F}_{y,\alpha}$, $\b{F}_{x,\beta}$ and $\b{F}_{y,\beta}$ the empirical cumulative distribution functions for $\{x_{\alpha,t}\}_{t\in[T]}$, $\{y_{\alpha,t}\}_{t\in[T]}$, $\{x_{\beta,t}\}_{t\in[T]}$ and $\{y_{\beta,t}\}_{t\in[T]}$ respectively:
\begin{equation}\label{eqn: cdf_x_y}
\begin{split}
&\b{F}_{x,\alpha}(c) := \frac{1}{T}\sum_{t=1}^T\1\{x_{\alpha,t}\leq c\}, \;\;\;
\b{F}_{y,\alpha}(c) := \frac{1}{T}\sum_{t=1}^T\1\{y_{\alpha,t}\leq c\},\\
&\b{F}_{x,\beta}(c) := \frac{1}{T}\sum_{t=1}^T\1\{x_{\beta,t}\leq c\}, \;\;\;
\b{F}_{y,\beta}(c) := \frac{1}{T}\sum_{t=1}^T\1\{y_{\beta,t}\leq c\}.
\end{split}
\end{equation}
\begin{corollary}\label{cor: test_statistic_FMvsMEFM}
 Let all the assumptions in Theorem \ref{thm: testingFMvsMEFM} hold. Define for $0<\theta<1$,
 \begin{align*}
 \wh{q}_{x,\alpha}(\theta) &:=  \inf\{c\;|\;\b{F}_{x,\alpha}(c)\geq \theta\}, \;\;\;
 \wh{q}_{x,\beta}(\theta) :=  \inf\{c\;|\;\b{F}_{x,\beta}(c)\geq \theta\}, \;\;\;
 \end{align*}
 Then under $H_0$, as $T,p,q \rightarrow \infty$, we have for each $t \in [T]$,
 \[\b{P}_{y,\alpha}[y_{\alpha,t} \geq \wh{q}_{x,\alpha}(\theta)] \leq 1-\theta, \;\;\;
 \b{P}_{y,\beta}[y_{\beta,t} \geq \wh{q}_{x,\beta}(\theta)] \leq 1-\theta,\]
 where $\b{P}_{y,\alpha}$ and $\b{P}_{y,\beta}$ are empirical probability measures induced by $\b{F}_{y,\alpha}$ and $\b{F}_{y,\beta}$ respectively.
\end{corollary}
With Corollary \ref{cor: test_statistic_FMvsMEFM}, we can test $H_0$ at significance level $1-\theta$ using the test statistics $y_{\alpha,t}$ and $y_{\beta,t}$, and rejection rules $y_{\alpha,t} \geq \wh{q}_{x,\alpha}(\theta)$ and $y_{\beta,t} \geq \wh{q}_{x,\beta}(\theta)$ respectively. Since we have $y_{\alpha,t}$ and $y_{\beta,t}$ for $t\in[T]$,  we can assess the significance level under $H_0$ by calculating
\[\text{Significance levels} = T^{-1}\sum_{t=1}^T\1\{y_{\alpha,t} \geq \wh{q}_{x,\alpha}(\theta)\},\;\;\; T^{-1}\sum_{t=1}^T\1\{y_{\beta,t} \geq \wh{q}_{x,\beta}(\theta)\},\]
and see if they are close to $1-\theta$. If $H_0$ is not true, then if any $\alpha_{t,i}$ is large, we expect $y_{\alpha,t}$ to be large. Or, if any $\beta_{t,j}$ is large, we expect $y_{\beta,t}$ to be large.

In practice for testing $H_0$ against $H_1$, we estimate $k_r$ and $k_c$, and set $\ell_r = k_r + 1$ and $\ell_c = k_c+1$ in light of the previous argument on how a special form of MEFM can be expressed back in FM. For estimation of $k_r$ and $k_c$, see Section \ref{sec: estimation_rank}.

\subsection{Estimation of the number of factors}\label{sec: estimation_rank}
From (\ref{eqn: Ct_estimator}), we have $T^{-1} \sum_{t=1}^T \wh\L_t\wh\L_t'$ essentially the row sample covariance matrix and $T^{-1} \sum_{t=1}^T \wh\L_t'\wh\L_t$ the column sample covariance matrix. We then propose the eigenvalue-ratio estimators for the number of factors as
\begin{align}
    \wh{k}_r := \argmin_{j} \bigg\{\frac{\lambda_{j+1}\big(T^{-1} \sum_{t=1}^T \wh\L_t\wh\L_t' \big) + \xi_r}{\lambda_{j}\big(T^{-1} \sum_{t=1}^T \wh\L_t\wh\L_t' \big) + \xi_r}, \; j\in[\lfloor p/2 \rfloor] \bigg\}, \;\;\;
    \xi_r \asymp pq\big[(Tq)^{-1/2} + p^{-1/2}\big],
    \label{eqn: hat_k_r} \\
    \wh{k}_c := \argmin_{j} \bigg\{\frac{\lambda_{j+1}\big(T^{-1} \sum_{t=1}^T \wh\L_t'\wh\L_t \big) + \xi_c}{\lambda_{j}\big(T^{-1} \sum_{t=1}^T \wh\L_t'\wh\L_t \big) + \xi_c}, \; j\in[\lfloor q/2 \rfloor] \bigg\}, \;\;\;
    \xi_c \asymp pq\big[(Tp)^{-1/2} + q^{-1/2}\big] .
    \label{eqn: hat_k_c}
\end{align}
Ratio-based estimators are widely studied by researchers. For example, an eigenvalue-ratio estimator is considered in \cite{LamYao2012} and \cite{AhnHorenstein2013}, while a cumulative eigenvalue ratio estimator is proposed by \cite{Zhangetal2023}. Our proposed estimator is similar to the perturbed eigenvalue-ratio estimators as in \cite{Pelger2019}. Technically, we can minimise (\ref{eqn: hat_k_r}) (resp. (\ref{eqn: hat_k_c}) over any $j\in[p]$ (resp. $j\in[q]$), but it is very reasonable to assume $k_r \leq p/2$ and $k_c \leq q/2$ in all applications of factor models. The correction terms $\xi_r$ and $\xi_c$ are added to stabilise the ratio so that consistency follows from the theorem below.

\begin{theorem}\label{thm: rank}
Under Assumptions (IC1), (M1), (F1), (L1), (E1), (E2) and (R1), we have the following.
\begin{itemize}
    \item [1.] $\wh{k}_r$ is a consistent estimator of $k_r$ if
    \[  \left\{
      \begin{array}{ll}
        p^{1-\delta_{r,k_r}} q^{1-\delta_{c,1}} [(Tq)^{-1/2} + p^{-1/2}]  = o(p^{\delta_{r,j+1} - \delta_{r,j}}), & \hbox{$j\in[k_r-1]$ with $k_r\geq 2$;} \\
        p^{1-\delta_{r,1}} q^{1-\delta_{c,1}} [(Tq)^{-1/2} + p^{-1/2}] = o(1), & \hbox{$k_r=1$.}
      \end{array}
    \right.
    \]
    \item [2.] $\wh{k}_c$ is a consistent estimator of $k_c$ if
    \[  \left\{
      \begin{array}{ll}
        q^{1-\delta_{c,k_c}} p^{1-\delta_{r,1}} [(Tp)^{-1/2} + q^{-1/2}]  = o(q^{\delta_{c,j+1} - \delta_{c,j}}), & \hbox{$j\in[k_c-1]$ with $k_c\geq 2$;} \\
        q^{1-\delta_{c,1}} p^{1-\delta_{r,1}} [(Tp)^{-1/2} + q^{-1/2}] = o(1), & \hbox{$k_c=1$.}
      \end{array}
    \right.
    \]
\end{itemize}
\end{theorem}
The extra rate conditions in the theorem are due to existence of potential weak factors and are trivially satisfied for pervasive factors. The theorem is similar to the consistency result in \cite{CenLam2024} for matrix-valued factor models, and this implies that the number of factors in MEFM can be well estimated just as in the case of FM.

\section{Numerical Results}\label{sec: numerical_results}
\subsection{Simulation}\label{subsec: simulation}
We demonstrate the performance of our estimators in this section. We will experiment different settings to assess consistency results as described in Theorem \ref{thm: consistency} and \ref{thm: consistency_F_C}, followed by the asymptotic normality of our estimators in Theorem \ref{thm:alpha_normality} and \ref{thm: asymp_loading}, where the covariance matrices can be constructed by their consistent estimators by Theorem \ref{thm:alpha_normality} and Theorem \ref{thm: HAC}, respectively. We then showcase the results for the rank estimators described in Theorem \ref{thm: rank}. As it is a first to consider matrix factor model with time-varying grand mean and main effects, we unveil the differences between MEFM and FM using numerical results that will illustrate Theorem \ref{thm: testingFMvsMEFM}.

For the data generating process, we use Assumptions (E1), (E2) and (F1) to generate general linear processes for the noise and factor series in model (\ref{eqn: model}). To be precise, the elements in $\F_t$ are independent standardised AR(5) with AR coefficients 0.7, 0.3, -0.4, 0.2, and -0.1. The elements in $\F_{e,t}$ and $\bepsilon_t$ are generated similarly, but their AR coefficients are (-0.7, -0.3, -0.4, 0.2, 0.1) and (0.8, 0.4, -0.4, 0.2, -0.1) respectively. The standard deviation of each element in $\bepsilon_t$ is generated by i.i.d. $|\cN(0,1)|$. To test how robust our method is under heavy-tailed distribution, we consider two distributions for the innovation process in generating $\F_t$, $\F_{e,t}$ and $\bepsilon_t$: 1) i.i.d. $\cN(0,1)$; 2) i.i.d. $t_3$.

The row factor loading matrix $\A_r$ is generated with $\A_r=\M_p\bf{U}_r\B_r$, where each entry of $\bf{U}_r\in\b{R}^{p\times k_r}$ is i.i.d. $\cN(0,1)$, and $\B_r\in\b{R}^{k_r\times k_r}$ is diagonal with the $j$-th diagonal entry being $p^{-\zeta_{r,j}}$, $0\leq \zeta_{r,j}\leq 0.5$. Pervasive (strong) factors have $\zeta_{r,j}=0$, while weak factors have $0<\zeta_{r,j}\leq 0.5$. Note that $\M_p$ is defined in (\ref{eqn: Ct_estimator}) so that (IC1) is satisfied. In a similar way, the column factor loading matrix $\A_c$ is generated independently. Each entry of $\A_{e,r}\in\b{R}^{p\times k_{e,r}}$ is i.i.d. $\cN(0,1)$ and has independent probability of 0.95 being set exactly to 0, and $\A_{e,c}$ is generated similarly. We fix $k_{e,r} = k_{e,c} = 2$ throughout the section.

For any $t\in[T]$, we generate $\mu_t = v_{\mu,t}$, $\balpha_t = \M_p \bf{v}_{\alpha,t}$ and $\bbeta_t = \M_q \bf{v}_{\beta,t}$, where $v_{\mu,t}$ is $\cN(m_\mu, \sigma_\mu^2)$, each element of $\bf{v}_{\alpha,t}$ is i.i.d. $\cN(m_\alpha, \sigma_\alpha^2)$ and that of $\bf{v}_{\beta,t}$ is i.i.d. $\cN(m_\beta, \sigma_\beta^2)$. We set $m_\mu = m_\alpha = m_\beta = 0$ and $\sigma_\mu = \sigma_\alpha = \sigma_\beta = 1$, and every experiment in this section is repeated 1000 times unless specified otherwise.

\subsubsection{Accuracy of various estimators}\label{subsubsec: simulation_consistency}
To assess the accuracy of our estimators, we define the relative mean squared errors (MSE) for $\mu_t$, $\balpha_t$, $\bbeta_t$ and $\C_t$ as the following, respectively,
\begin{align*}
    & \text{relative MSE}_\mu = \frac{\sum_{t=1}^T(\mu_t - \wh{\mu}_t)^2}{\sum_{t=1}^T \mu_t^2} ,
    \,\,\,\,\,
    \text{relative MSE}_{\balpha} = \frac{\sum_{t=1}^T\|\balpha_t - \wh{\balpha}_t\|^2}{\sum_{t=1}^T \|\balpha_t\|^2} ,\\
    & \text{relative MSE}_{\bbeta} = \frac{\sum_{t=1}^T\|\bbeta_t - \wh{\bbeta}_t\|^2}{\sum_{t=1}^T \|\bbeta_t\|^2} ,
    \,\,\,\,\,
    \text{relative MSE}_{\C} = \frac{\sum_{t=1}^T\|\C_t - \wh{\C}_t\|_F^2}{\sum_{t=1}^T \|\C_t\|_F^2} .
\end{align*}
For measuring the accuracy of our factor loading matrix estimators, we use the column space distance,
\[
\c{D}(\Q, \wh\Q) = \Big\|\Q(\Q'\Q)^{-1}\Q' -\wh\Q(\wh\Q'\wh\Q)^{-1}\wh\Q' \Big\|,
\]
for any given $\Q$ and $\wh\Q$, which is a common measure in the literature such as \cite{ChenYangZhang2022} and \cite{ChenFan2023}.

We consider the following settings:
\begin{itemize}
\item[(Ia)] $T=100, \; p=q=40, \; k_r=1, \; k_c=2$. All factors are pervasive with $\zeta_{r,j} = \zeta_{c,j} = 0$. All innovation processes in constructing $\F_t$, $\F_{e,t}$ and $\bepsilon_t$ are i.i.d. standard normal.
\item[(Ib)] Same as (Ia), but one factor is weak with $\zeta_{r,1}=0.2$ and $\zeta_{c,1}=0.2$. Set also $m_\alpha = -2$.
\item[(Ic)] Same as (Ia), but all innovation processes are i.i.d. $t_3$.
\item[(Id)] Same as (Ib), but $T=100, \; p=q=80$ and $\sigma_\alpha = 2$.
\item[(Ie)] Same as (Id), but $T=200$.
\item[(IIa-e)] Same as (Ia) to (Ie) respectively, except that we generate $\F_t$, $\F_{e,t}$ and $\bepsilon_t$ using white noise rather than AR(5).
\end{itemize}
Setting (IIa) to (IIe) are to investigate how temporal dependence in the noise affects our results.

\begin{figure}[t!]
\centering
\begin{minipage}{.45\textwidth}
  \centering
  \includegraphics[width=\columnwidth]{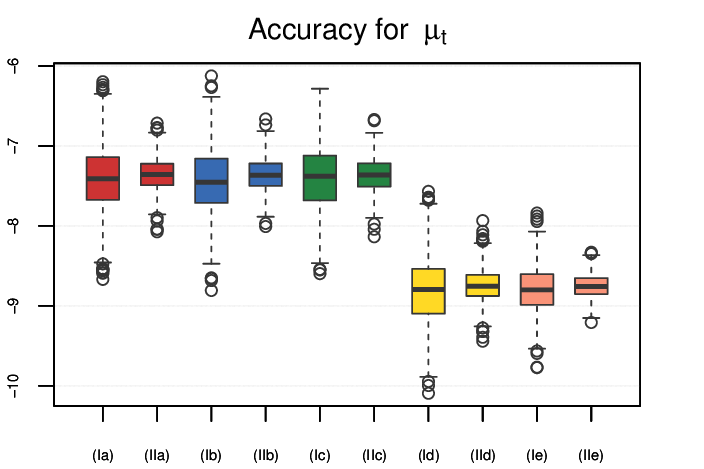}
  \caption{Plot of the relative MSE for $\mu_t$ (in log-scale) from Settings (Ia) to (Ie), comparing with (IIa) to (IIe).}
  \label{fig: consistency_mu}
\end{minipage}%
\hspace{0.05\textwidth}
\begin{minipage}{.45\textwidth}
  \centering
  \includegraphics[width=\columnwidth]{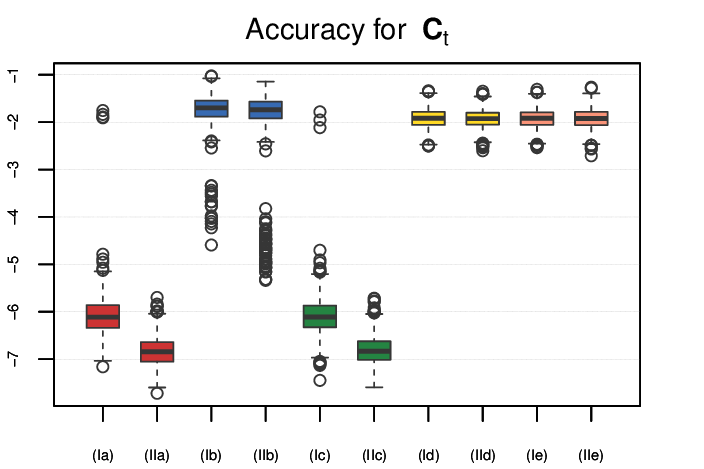}
  \caption{Plot of the relative MSE for $\C_t$ (in log-scale) from Settings (Ia) to (Ie), comparing with (IIa) to (IIe).}
  \label{fig: consistency_C}
\end{minipage}
\end{figure}

\begin{figure}[h!]
\centering
\begin{minipage}{.45\textwidth}
  \centering
  \includegraphics[width=\columnwidth]{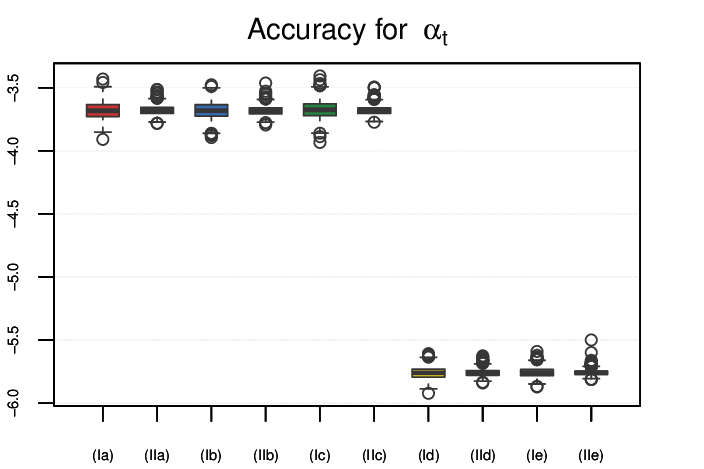}
  \caption{Plot of the relative MSE for $\balpha_t$ (in log-scale) from Settings (Ia) to (Ie), comparing with (IIa) to (IIe).}
  \label{fig: consistency_alpha}
\end{minipage}%
\hspace{0.05\textwidth}
\begin{minipage}{.45\textwidth}
  \centering
  \includegraphics[width=\columnwidth]{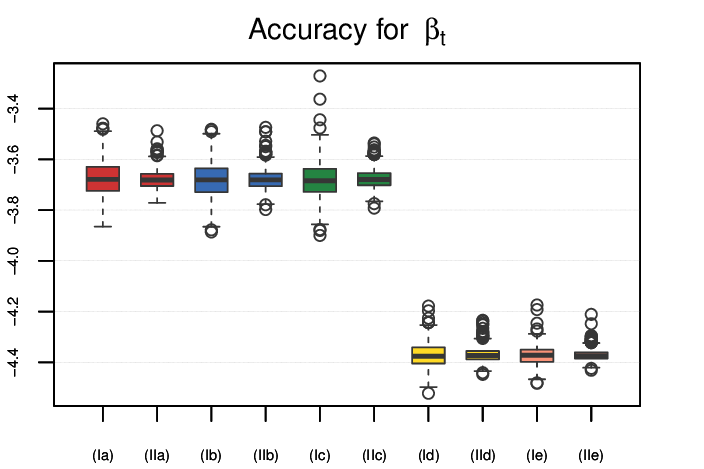}
  \caption{Plot of the relative MSE for $\bbeta_t$ (in log-scale) from Settings (Ia) to (Ie), comparing with (IIa) to (IIe).}
  \label{fig: consistency_beta}
\end{minipage}
\end{figure}

\begin{figure}[h!]
\centering
\begin{minipage}{.45\textwidth}
  \centering
  \includegraphics[width=\columnwidth]{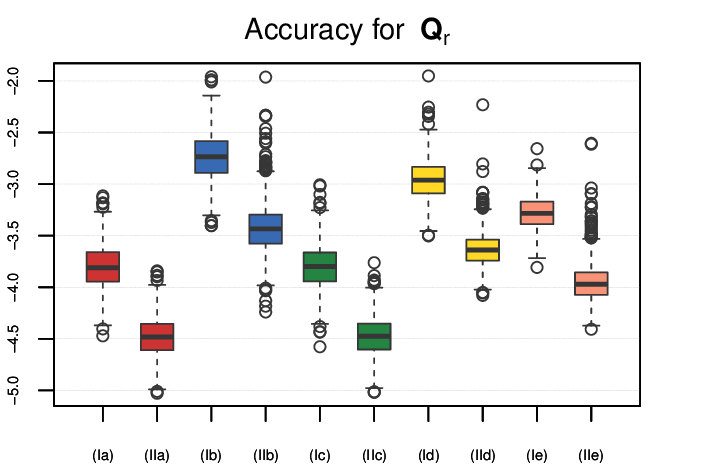}
  \caption{Plot of the row space distance $\c{D}(\Q_r, \wh\Q_r)$ (in log-scale) from Settings (Ia) to (Ie), comparing with (IIa) to (IIe).}
  \label{fig: consistency_Qr}
\end{minipage}%
\hspace{0.05\textwidth}
\begin{minipage}{.45\textwidth}
  \centering
  \includegraphics[width=\columnwidth]{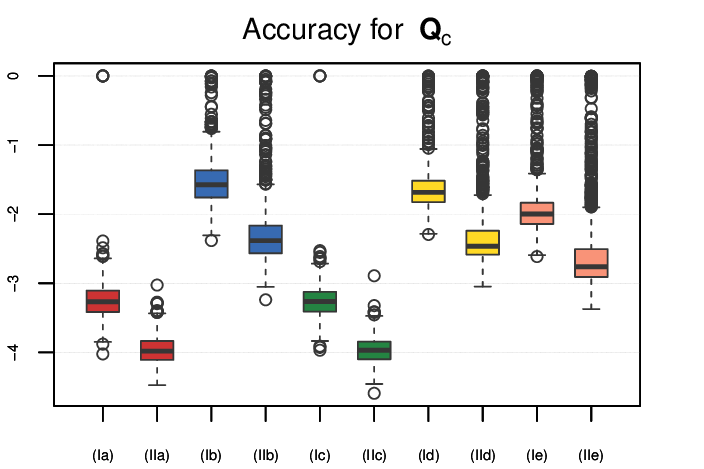}
  \caption{Plot of the column space distance $\c{D}(\Q_c, \wh\Q_c)$ (in log-scale) from Settings (Ia) to (Ie), comparing with (IIa) to (IIe).}
  \label{fig: consistency_Qc}
\end{minipage}
\end{figure}

We report the boxplots of accuracy measures for our estimators from Figure \ref{fig: consistency_mu} to Figure \ref{fig: consistency_Qc}. Note first that stronger temporal dependence leads to larger variance of our estimators in general. The serial dependence mainly undermines the performance of our loading matrix estimators as shown in Figures \ref{fig: consistency_Qr} and \ref{fig: consistency_Qc}, which in turn affects our common component estimator.

Considering the comparisons among (Ia) to (Ie), we see that $\text{relative MSE}_\mu$ can be improved by increasing the spatial dimensions, but is not affected by weak factors. Similar results can be seen from Figure \ref{fig: consistency_alpha} and Figure \ref{fig: consistency_beta} for $\text{relative MSE}_\alpha$ and $\text{relative MSE}_\beta$. The detrimental effects of heavy-tailed innovation processes in Setting (Ic) are most reflected in the corresponding boxplots in Figure \ref{fig: consistency_beta}.

Weak factors can be detrimental to the accuracy of the factor loading matrix estimators, as can be seen by the significant rise in the factor loading space errors from Setting (Ia) to (Ib) in Figure \ref{fig: consistency_Qr} and \ref{fig: consistency_Qc}. In fact, $\wh{k}_c$ barely captures the second factor under Setting (Ib) and (IIb). See Section \ref{subsubsec: simulation_rank} for details. Comparing Setting (Ib) with (Id), Figure \ref{fig: consistency_Qr} and \ref{fig: consistency_Qc} show that increase in data dimensions slightly improves our factor loading matrix estimators, which is consistent to the simulation results in \cite{WangLiuChen2019} for instance.

\subsubsection{Performance for the estimation of the number of factors}\label{subsubsec: simulation_rank}
In this section, we demonstrate the performance of our estimators for the number of factors, as described in Theorem \ref{thm: rank}. First, we set $\xi_r = pq[(Tq)^{-1/2} + p^{-1/2}] /5$ and $\xi_c = pq[(Tp)^{-1/2} + q^{-1/2}] /5$, so that the conditions for $\xi_r$ and $\xi_c$ in (\ref{eqn: hat_k_r}) and (\ref{eqn: hat_k_c}) are respectively satisfied. A wide range of values other than $1/5$ for $\xi_r$ and $\xi_c$ are experimented, but $1/5$ is working the best in vast majority of settings, and hence we do not recommend treating this as a tuning parameter.

We present the results for each of the following settings:
\begin{itemize}
\item[(IIIa)] $k_r = k_c = 3$. All factors are pervasive with $\zeta_{r,j} = \zeta_{c,j} = 0$ for all $j\in[3]$. All innovation processes involved are i.i.d. standard normal.
\item[(IIIb)] Same as (IIIa), but some factors are weak with $\zeta_{r,1} = \zeta_{c,1} = \zeta_{c,2} =0.2$.
\item[(IIIc)] Same as (IIIa), but all factors are weak with $\zeta_{r,j} = \zeta_{c,j}= 0.2$ for all $j\in[3]$.
\end{itemize}
We experiment the above settings with various $(p,q)$ pairs among $(10, 10)$, $(10, 20)$ and $(20, 20)$, with the choice $T= 0.5 \cdot pq$ or $T = pq$. The setup is similar to \cite{WangLiuChen2019} and \cite{ChenFan2023}, but we use smaller sets of dimensions since the accuracy of our estimators are approaching 1 with larger dimensions, which reveal little intricacies among different settings.

\begin{table}[h!]
\begin{center}
\begin{tabular}{|c|cccccc|}
\hline
 & \multicolumn{2}{c|}{$p,q=10,10$}                                      & \multicolumn{2}{c|}{$p,q=10,20$}                                      & \multicolumn{2}{c|}{$p,q=20,20$}                 \\ \cline{2-7}
    $(\wh{k}_r, \wh{k}_c)$      & \multicolumn{1}{c|}{$T=.5pq$}      & \multicolumn{1}{c|}{$T=pq$}        & \multicolumn{1}{c|}{$T=.5pq$}      & \multicolumn{1}{c|}{$T=pq$}        & \multicolumn{1}{c|}{$T=.5pq$}      & $T=pq$        \\ \hline
        & \multicolumn{6}{c|}{Setting (IIIa)}                \\ \hline
$(2,3)$           & \multicolumn{1}{c|}{0.121}           & \multicolumn{1}{c|}{0.112}           & \multicolumn{1}{c|}{0.128}           & \multicolumn{1}{c|}{0.11}           & \multicolumn{1}{c|}{0}           & 0.004          \\ \hline
$(3,2)$             & \multicolumn{1}{c|}{0.124}           & \multicolumn{1}{c|}{0.111}           & \multicolumn{1}{c|}{0.004}           & \multicolumn{1}{c|}{0.003}          & \multicolumn{1}{c|}{0.001}          & 0.001          \\ \hline
$(3,3)$             & \multicolumn{1}{c|}{\textbf{0.583}} & \multicolumn{1}{c|}{\textbf{0.659}} & \multicolumn{1}{c|}{\textbf{0.833}} & \multicolumn{1}{c|}{\textbf{0.855}} & \multicolumn{1}{c|}{\textbf{0.999}} & \textbf{0.995} \\ \hline
other             & \multicolumn{1}{c|}{0.172}          & \multicolumn{1}{c|}{0.118}          & \multicolumn{1}{c|}{0.035}          & \multicolumn{1}{c|}{0.032}          & \multicolumn{1}{c|}{0}          & 0          \\ \hline
           & \multicolumn{6}{c|}{Setting (IIIb)}           \\ \hline
$(2,3)$             & \multicolumn{1}{c|}{0.135}          & \multicolumn{1}{c|}{0.13}          & \multicolumn{1}{c|}{0.23}          & \multicolumn{1}{c|}{0.257}          & \multicolumn{1}{c|}{0.228}          & 0.149          \\ \hline
$(3,2)$             & \multicolumn{1}{c|}{0.079}          & \multicolumn{1}{c|}{0.096}          & \multicolumn{1}{c|}{0.024}          & \multicolumn{1}{c|}{0.017}          & \multicolumn{1}{c|}{0.022}          & 0.02          \\ \hline
$(3,3)$             & \multicolumn{1}{c|}{\textbf{0.136}} & \multicolumn{1}{c|}{\textbf{0.17}} & \multicolumn{1}{c|}{\textbf{0.289}} & \multicolumn{1}{c|}{\textbf{0.347}} & \multicolumn{1}{c|}{\textbf{0.556}} & \textbf{0.637} \\ \hline
other             & \multicolumn{1}{c|}{0.65}          & \multicolumn{1}{c|}{0.604}          & \multicolumn{1}{c|}{0.457}          & \multicolumn{1}{c|}{0.379}          & \multicolumn{1}{c|}{0.194}          & 0.194          \\ \hline
          & \multicolumn{6}{c|}{Setting (IIIc)}                    \\ \hline
$(2,3)$             & \multicolumn{1}{c|}{0.082}          & \multicolumn{1}{c|}{0.085}          & \multicolumn{1}{c|}{0.218}          & \multicolumn{1}{c|}{0.254}          & \multicolumn{1}{c|}{0.089}          & 0.096          \\ \hline
$(3,2)$             & \multicolumn{1}{c|}{0.075}          & \multicolumn{1}{c|}{0.124}          & \multicolumn{1}{c|}{0.04}          & \multicolumn{1}{c|}{0.035}          & \multicolumn{1}{c|}{0.088}          & 0.089          \\ \hline
$(3,3)$             & \multicolumn{1}{c|}{\textbf{0.073}} & \multicolumn{1}{c|}{\textbf{0.096}} & \multicolumn{1}{c|}{\textbf{0.209}} & \multicolumn{1}{c|}{\textbf{0.257}} & \multicolumn{1}{c|}{\textbf{0.614}} & \textbf{0.646} \\ \hline
other             & \multicolumn{1}{c|}{0.77}          & \multicolumn{1}{c|}{0.695}          & \multicolumn{1}{c|}{0.533}          & \multicolumn{1}{c|}{0.454}          & \multicolumn{1}{c|}{0.209}          & 0.169          \\ \hline
\end{tabular}
\end{center}
\caption {Results for Setting (IIIa) to (IIIc). Each cell reports the frequency of $(\wh{k}_r, \wh{k}_c)$ under the setting in the corresponding column. The true number of factors is $(k_r, k_c) = (3,3)$, and the cells corresponding to correct estimations are bolded.}
\label{tab: simulation_rank}
\end{table}

From the results in Table \ref{tab: simulation_rank}, our eigenvalue-ratio estimators is working well with MEFM. The accuracy of $\wh{k}_r$ and $\wh{k}_c$ suffers from the existence of weak factors, which is also seen in traditional FM  (see for instance \cite{ChenLam2024} and \cite{CenLam2024}). In particular, the accuracy of our estimators drops significantly as we move from Setting (IIIa) to (IIIc), and in general large dimensions are beneficial to our estimation. Lastly, note that although we have two weak factors in the column loading matrix while there is only one weak factor in the row loading matrix, the correct proportion of $\wh{k}_c$ is much larger than that of $\wh{k}_r$ for $(p,q) = (10,20)$. This hints at the importance of data dimensions over factor strength, which can also be seen from the fact that the results for $(p,q) = (20,20)$ under Setting (IIIc) are comparable with those for $(p,q) = (10,10)$ under Setting (IIIa).

\subsubsection{Asymptotic normality}\label{subsubsec: simulation_asymp}
We numerically demonstrate the asymptotic normality results in Theorems \ref{thm:alpha_normality} and \ref{thm: asymp_loading} in this section. For the ease of demonstration, we consider $t=10$ only for the asymptotic distribution of $\wh\mu_t$, $\wh\btheta_{\alpha, t}=(\wh\alpha_{t,1}, \wh\alpha_{t,2}, \wh\alpha_{t,3})'$ and $\wh\btheta_{\beta, t}=(\wh\beta_{t,1}, \wh\beta_{t,2}, \wh\beta_{t,3})'$, and for $\wh\btheta_{\alpha, t}$ and $\wh\btheta_{\beta, t}$ we will only report results for the third component. We will also demonstrate the asymptotic normality for $(\wh\Q_c)_{1\cdot}$ and present the results for $(\wh\Q_c)_{11}$, i.e., the first entry of the first row in the column loading matrix estimator. To consistently estimate its covariance matrix, we use Theorem \ref{thm: HAC} with $\eta_c = \lfloor (Tpq)^{1/4} /5 \rfloor$.

\begin{figure}[t!]
\begin{center}
\centerline{\includegraphics[width=0.9\columnwidth,scale=0.9]{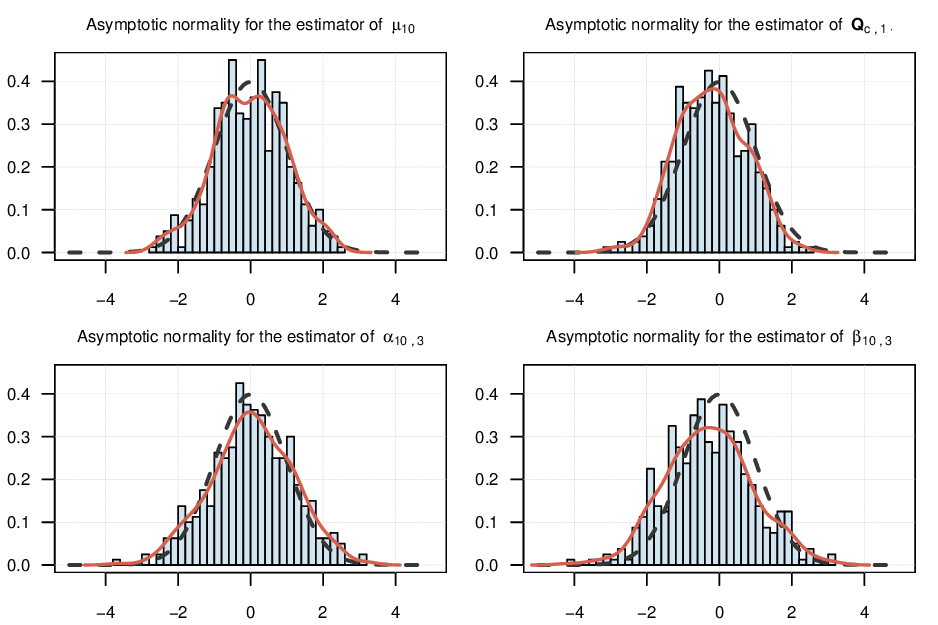}}
\caption{Histograms of $\sqrt{pq}\,\wh\gamma_{\mu}^{-1}(\wh\mu_{10} - \mu_{10})$ (top-left), $[T\, (\wh\bSigma_{c,1}^{HAC})^{-1/2} \wh\D_c (\wh\Q_{c,1\cdot}-\H_1^a\Q_{c,1\cdot})]_1$ (top-right), $\sqrt{q}\, [\text{diag}(\wh\gamma_{\alpha,1}^{-1}, \wh\gamma_{\alpha,2}^{-1}, \wh\gamma_{\alpha,3}^{-1})(\wh\btheta_{\alpha,10} - \btheta_{\alpha,10})]_3$ (bottom-left), and $\sqrt{p}\, [\text{diag}(\wh\gamma_{\beta,1}^{-1}, \wh\gamma_{\beta,2}^{-1}, \wh\gamma_{\beta,3}^{-1})(\wh\btheta_{\beta,10} - \btheta_{\beta,10})]_3$ (bottom-right). In each panel, the curve (in red) is the empirical density, and the density curve for $\cN(0,1)$ (in black, dotted) is also superimposed.}
\label{fig: asymp_hist}
\end{center}
\end{figure}

We use heavy-tailed innovations to investigate the robustness of our results, hence Setting (Ic) is adapted except that we generate $\F_t$, $\F_{e,t}$ and $\bepsilon_t$ using AR(1) with coefficient $-0.2$. Due to the different rates of convergence in Theorems \ref{thm:alpha_normality} and \ref{thm: asymp_loading}, we specify different dimensions $(T,p,q)$ in the following settings:
\begin{align*}
    & \wh\mu_t: (80,100,100), \;\;\;
    \wh\btheta_{\alpha, t}: (60,60,300), \;\;\;
    \wh\btheta_{\beta, t}: (60,300,60), \;\;\;
    (\wh\Q_c)_{1\cdot}: (60,60,300),
\end{align*}
where the dimension setting for $(\wh\Q_c)_{1\cdot}$ is to align with the rate conditions in Theorem \ref{thm: asymp_loading} that $Tp/q^2 \to 0$ under pervasive factors. Each setting is repeated 400 times, and we present the histograms of our four estimators in Figure \ref{fig: asymp_hist}.

Our plots stand as empirical evidence of Theorem \ref{thm:alpha_normality}, \ref{thm: asymp_loading} and \ref{thm: HAC}. It might worth noting that the spread of the normalised empirical density for $\wh{\beta}_{10,3}$ is slightly larger than expected by comparing with the superimposed standard normal. The same problem is not seen in the histogram for $\wh{\alpha}_{10,3}$. With true $(k_r,k_c)=(1,2)$, the common component estimation using $(p,q)=(300,60)$ is worse than that using $(p,q)=(60,300)$ due to insufficient column dimension relative to $k_c$. Hence it leads to worse estimators for errors and $(\wh\gamma_{\beta,1}^{-1}, \wh\gamma_{\beta,2}^{-1}, \wh\gamma_{\beta,3}^{-1})$ under $(p,q)=(300,60)$. Hence inference performances on the time-varying row and column effect estimators are affected by the latent number of factors.

\subsubsection{Testing MEFM versus FM}\label{subsubsec: simulation_test}
We demonstrate numerical results for Corollary \ref{cor: test_statistic_FMvsMEFM} in this section. We consider the two scenarios below:
\begin{itemize}
\item [1.] \textit{(Global effect.)} The entries of at least one of $\balpha_t$ and $\bbeta_t$ are in general non-zero for each $t$.
\item [2.] \textit{(Local effect.)} The entries of at least one of $\balpha_t$ and $\bbeta_t$ are sparse for each $t$, i.e., given any $t$, there are some non-zero entries in at least one of $\balpha_t$ and $\bbeta_t$ with all other entries zero.
\end{itemize}

Throughout this section, we generate the time-varying grand mean and main effects using Rademacher random variables such that $v_{\mu,t}$ is i.i.d. Rademacher multiplied by some $u_\mu$ and each entry of $\v_{\alpha,t}, \, \v_{\beta,t}$ is i.i.d. Rademacher multiplied by some $u_\alpha, \, u_\beta$ respectively, recalling that $\mu_t = v_{\mu,t}, \; \balpha_t = \M_p\v_{\alpha,t}$ and $\bbeta_t = \M_q\v_{\beta,t}$. Hence, setting $u_\mu = u_\alpha = u_\beta = 0$ corresponds to generating a traditional FM. We set $k_r = k_c = 2$, and consider the following settings:
\begin{itemize}
    \item [(IVa)] $T=p=q=40$. All factors are pervasive with $\zeta_{r,j} = \zeta_{c,j} = 0$. All innovation processes in constructing $\F_t$, $\F_{e,t}$ and $\bepsilon_t$ are i.i.d. standard normal. Set $u_\mu = u_\beta = 0$, and we select $u_\alpha$ from $0.1,\, 0.5,\, 1$.
    \item [(IVb)] Same as (IVa), but fix $u_\alpha = 0.1$ and select $u_\beta$ from $0.1,\, 0.5,\, 1$.
    \item [(IVc)] Same as (IVa), except that $u_\alpha = 1$, and when generating $\balpha_t = \M_p \bf{v}_{\alpha,t}$ as specified previously, we only keep the first $u_{local}$ entries of $\bf{v}_{\alpha,t}$ as non-zero where $u_{local}$ is selected from $2, \, 5, \, 10$.
\end{itemize}

Setting (IVa) and (IVb) are designed for testing global effects, and Setting (IVc) for local effects. For each setting, we construct $y_{\alpha,t}$, $y_{\beta,t}$ and use $\theta=0.95$ in Corollary \ref{cor: test_statistic_FMvsMEFM}. Each experiment is repeated 400 times and we report both $\textbf{reject}_\alpha := T^{-1} \sum_{t=1}^T \b{1}\{y_{\alpha,t} \geq \wh{q}_{x,\alpha}(0.95)\}$ and $\textbf{reject}_\beta := T^{-1} \sum_{t=1}^T \b{1}\{y_{\beta,t} \geq \wh{q}_{x,\beta}(0.95)\}$.

As explained under Corollary \ref{cor: test_statistic_FMvsMEFM}, we expect $\textbf{reject}_\alpha$ and $\textbf{reject}_\beta$ to be close to $1-\theta = 0.05$ if FM is sufficient. From Table \ref{tab: simulation_testFM}, our proposed test works well since it suggests FM is insufficient as we strengthen $\balpha_t$ or $\bbeta_t$. In particular, even if the signal of $\balpha_t$ is not strong enough such as $u_\alpha = 0.1$, Setting (IVb) shows that additional signals from $\bbeta_t$ allows us to reject the use of FM. The comparison between $\textbf{reject}_\alpha$ and $\textbf{reject}_\beta$ is indicative of which effect is stronger. According to the results for (IVc) in the table, our test is capable of detecting local effect such that $\textbf{reject}_\alpha$ is far from $0.05$ even when only two entries in $\balpha_t$ are non-zero.

\begin{table}[t!]\centering
\small
\ra{1.3}
\begin{tabular}{@{}rrcrrrcrrrcrrr@{}}\toprule
& \multicolumn{1}{c}{Size} & \phantom{abc}
& \multicolumn{3}{c}{ Setting (IVa)} & \phantom{abc} & \multicolumn{3}{c}{ Setting (IVb)} & \phantom{abc} & \multicolumn{3}{c}{ Setting (IVc)} \\
\cmidrule{4-6} \cmidrule{8-10} \cmidrule{12-14}
Parameter
& $0$  &&  $0.1$ & $0.5$ & $1$ && $0.1$ & $0.5$ & $1$ && $2$ & $5$ & $10$ \\ \midrule
$\textbf{reject}_\alpha$ & $5_{(4)}$ && $11_{(7)}$ & $63_{(31)}$ & $96_{(15)}$ && $13_{(8)}$ & $53_{(30)}$ & $86_{(23)}$ && $37_{(17)}$ & $77_{(24)}$ & $85_{(27)}$\\
$\textbf{reject}_\beta$ & $5_{(4)}$ && $11_{(7)}$ & $52_{(28)}$ & $87_{(22)}$ && $13_{(8)}$ & $62_{(32)}$ & $96_{(16)}$ && $14_{(8)}$ & $28_{(16)}$ & $48_{(26)}$ \\
\bottomrule
\end{tabular}
\caption {Results for Setting (IVa) to (IVd). Each cell reports the mean and SD (subscripted, in bracket), both multiplied by $100$. The parameters for Settings (IVa), (IVb) and (IVc) are $u_\alpha$, $u_\beta$ and $u_{local}$, respectively. Setting (IVa) with $u_\alpha = 0$ is reported in the first column, representing the size of the test.}
\label{tab: simulation_testFM}
\end{table}

\begin{figure}[t!]
\centering
\begin{minipage}{.45\textwidth}
  \centering
  \includegraphics[width=\columnwidth]{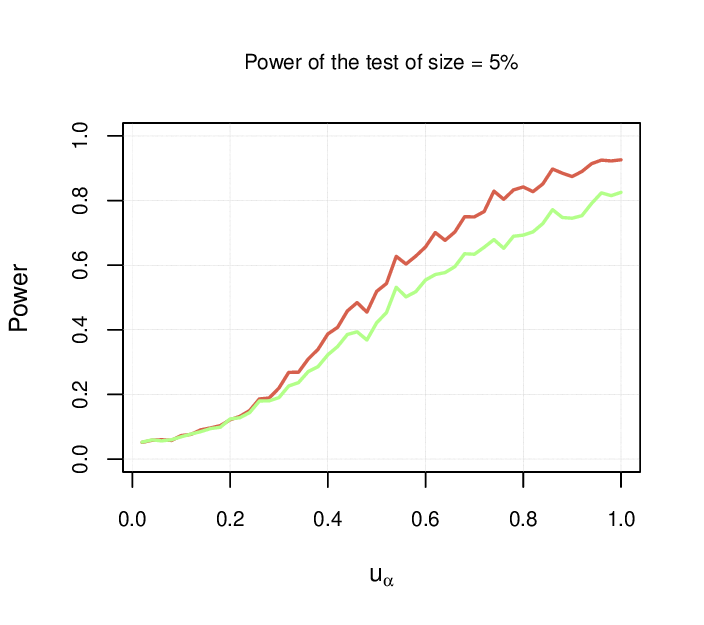}
  \caption{Statistical power curve of testing the null hypothesis that FM is sufficient for the given series, against the alternative that MEFM is necessary. Each power value is computed as the average over 400 runs of $\textbf{reject}_\alpha$ (in red) and $\textbf{reject}_\beta$ (in green) under Setting (IVa) except that $(T,\, p,\, q) = (60, \, 80, \, 80)$.}
  \label{fig: testFM_power_global}
\end{minipage}%
\hspace{0.05\textwidth}
\begin{minipage}{.45\textwidth}
  \centering
  \vspace{45pt}
  \includegraphics[width=\columnwidth]{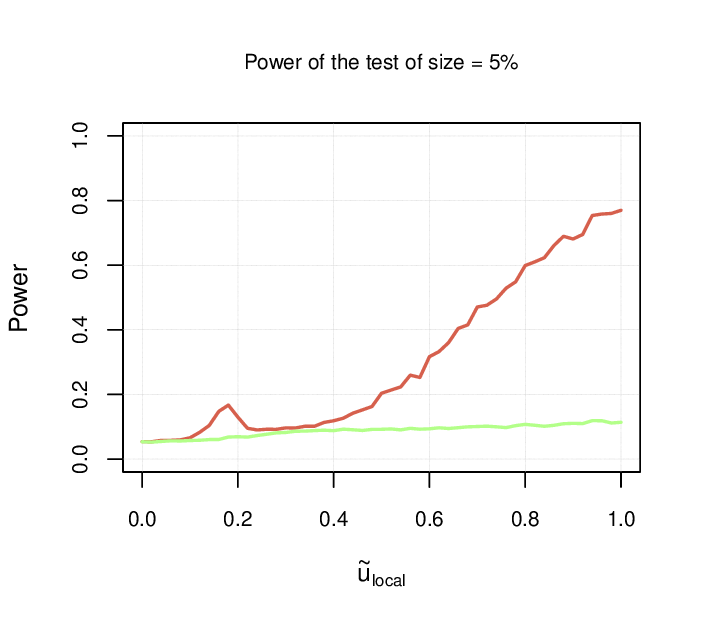}
  \caption{Statistical power curve of testing the null hypothesis that FM is sufficient for the given series, against the alternative that MEFM is necessary. Refer to Figure \ref{fig: testFM_power_global} for how the power is computed. The data is generated under Setting (IVc) except that $(T,\, p,\, q) = (60, \, 80, \, 80)$ and $\balpha_t$ is generated as $\balpha_1 = \wt{u}_{local} \, (1, \, 1, \, -2, \, 0,\, \dots,\, 0)'$, $\balpha_2 = \wt{u}_{local} \, (1, \, 2, \, -3, \, 0, \, \dots, \, 0)'$,  $\balpha_3 = \wt{u}_{local} \, (2, \, -5, \, 3, \, 0, \, \dots, \, 0)'$ and $\balpha_{3\ell + i} = \balpha_i$ for $\ell$ a positive integer and $i=1,2,3$, so that each $\balpha_t$ has non-zero entries only in the first three indices.}
  \label{fig: testFM_power_local}
\end{minipage}
\end{figure}

Extensive experiments on different dimensions, factor strengths or grand mean magnitudes are performed. All indicate similar interpretation as the above settings and hence the results are not shown here. The power curve for Setting (IVa) is also presented in Figure \ref{fig: testFM_power_global} to support the use of our test, with $(T,\, p,\, q) = (60, \, 80, \, 80)$ and $u_\alpha$ ranging from $0.02$ to $1$. Besides, we also show the power curve for local effect in Figure \ref{fig: testFM_power_local}, for Setting (IVc) except that $(T,\, p,\, q) = (60, \, 80, \, 80)$ and we generate $\balpha_t$ as described in the caption. Both power curves show that the test is able to reject the use of FM if signals from the time-varying main effects are large, either globally or locally. And in both figures, when $u_{\alpha}$ is close to 0.02 or $\wt{u}_{local}$ close to 0, the value of the power curves are all very close to 0.05, which is exactly what we want for the size of the tests.

%
%

\subsection{Real data analysis}\label{subsec: real_data}

\subsubsection{NYC taxi traffic}\label{subsubsec: taxi}
We analyse a set of taxi traffic data in New York city in this example. The data includes all individual taxi rides operated by Yellow Taxi in New York City, published at

\url{https://www1.nyc.gov/site/tlc/about/tlc-trip-record-data.page}.

For simplicity, we only consider the rides within Manhattan Island, which comprises most of the data. The dataset contains 842 million trip records within the period of January 1, 2013 to December 31, 2022. Each trip record includes features such as pick-up and drop-off dates/times, pick-up and drop-off locations, trip distances, itemized fares, rate types, payment types, and driver-reported passenger counts. Our example here focuses on the drop-off dates/times and locations.

To classify the drop-off locations in Manhattan, they are coded according to 69 predefined zones in the dataset. Moreover, each day is divided into 24 hourly periods to represent the drop-off times each day, with the first hourly period from 0 a.m. to 1 a.m. The total number of rides moving among the zones within each hour are recorded, yielding data $\Y_t\in\b{R}^{69\times 24}$ each day, where $y_{i_1,i_2,t}$ is the number of trips to zone $i_1$ and the pick-up time is within the $i_2$-th hourly period on day $t$.

We consider the non-business-day series which is 1,133 days long, within the period of January 1, 2013 to December 31, 2022. Using MEFM, the estimated rank of the core factors is $(2,\, 2)$ according to our proposed eigenvalue ratio estimator.
As mentioned in Section \ref{subsec: testFMvsMEFM}, we therefore use $(3, \, 3)$ as the number of factors to estimate FM and test if FM is sufficient. We compute $\textbf{reject}_\alpha = 0.064$ and $\textbf{reject}_\beta = 0.133$ which are defined in Section \ref{subsubsec: simulation_test}. They should be close to $1-\theta = 0.05$ according to Corollary \ref{cor: test_statistic_FMvsMEFM} if FM is sufficient.
Hence we reject the use of traditional FM due to the signals in $\wh\bbeta_t$.

\begin{figure}[h!]
\begin{center}
\centerline{\includegraphics[width=0.9\columnwidth,scale=0.9]{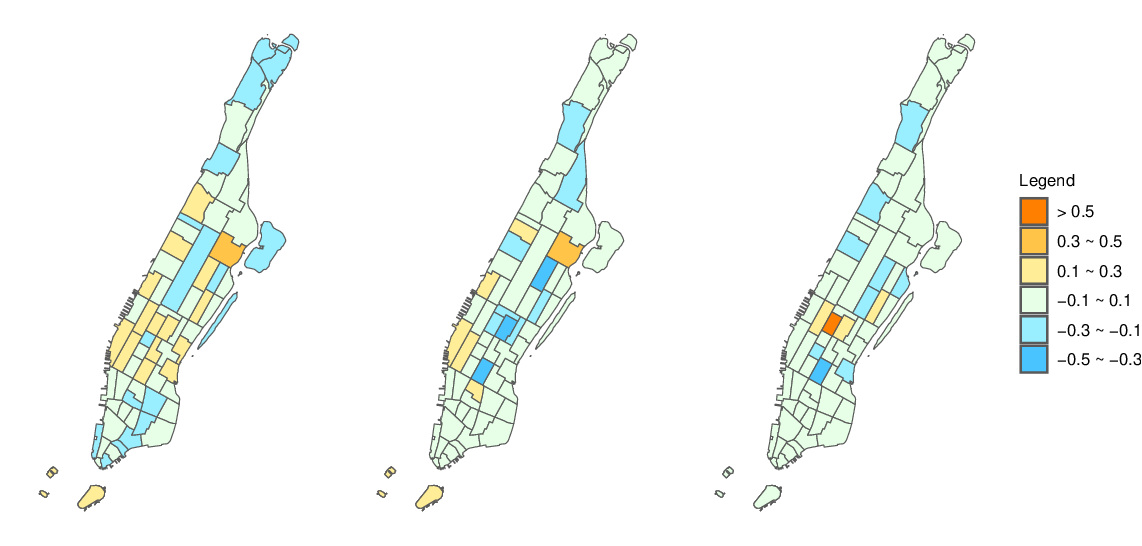}}
\caption{Estimated loading on three dropoff factors using MEFM, i.e. $\wh\Q_{1,\cdot 1}$ (left), $\wh\Q_{1,\cdot 2}$ (middle) and $\wh\Q_{1,\cdot 3}$ (right).}
\label{fig: real_data_Q1_MEFM}
\end{center}
\end{figure}

\begin{figure}[h!]
\begin{center}
\centerline{\includegraphics[width=0.9\columnwidth,scale=0.9]{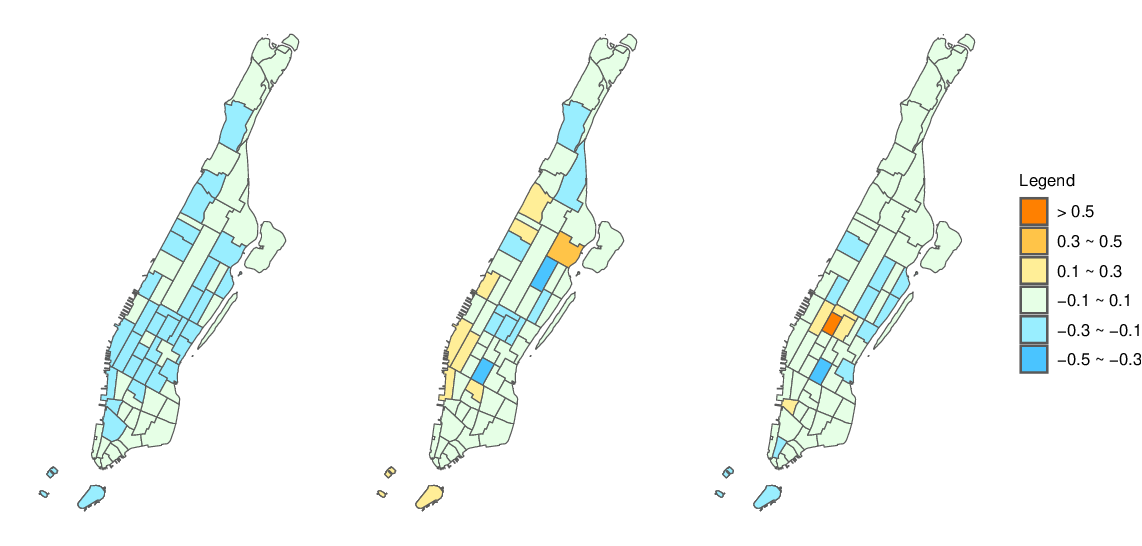}}
\caption{Estimated loading on three dropoff factors using FM, similar to Figure \ref{fig: real_data_Q1_MEFM}.}
\label{fig: real_data_Q1_FM}
\end{center}
\end{figure}

{
To compare MEFM with FM, we use core rank $(3,\, 3)$ to estimate MEFM for the rest of this section. Figure \ref{fig: real_data_Q1_MEFM} and \ref{fig: real_data_Q1_FM} illustrate the heatmaps of the estimated loading columns on the three dropoff factors using MEFM and FM, respectively. From both heatmaps, we can identify the first factor as active areas, the second as dining and sports areas and the third as downtown areas. The three factors are similar to their corresponding counterparts, except that the first factor estimated using MEFM is more indicative on the active areas to taxi traffic in Manhattan by its emphasised orange zone which corresponds to East Harlem.
}

\begin{figure}[h!]
\begin{center}
\centerline{\includegraphics[width=0.8\columnwidth,scale=0.8]{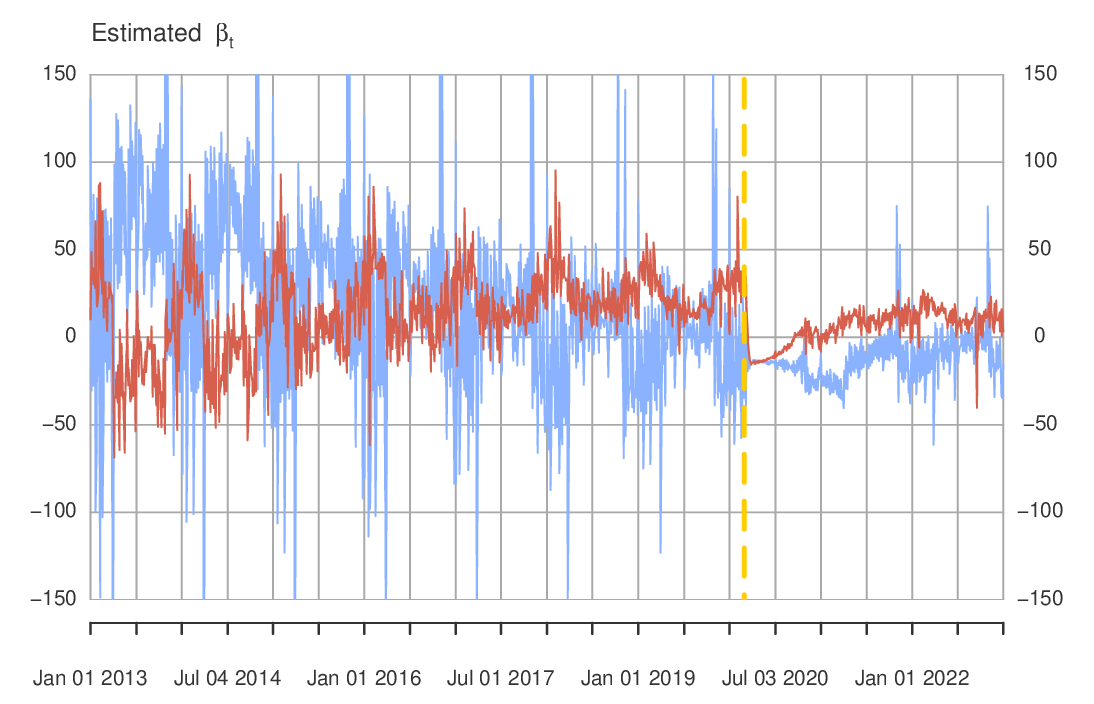}}
\caption{Plot of the estimated hour effects for periods from 1 a.m. to 2 a.m. (in blue) and from 5 p.m. to 6.p.m. (in red). The date for the first confirmed case of COVID-19 in New York is also shown (dotted  yellow vertical line).}
\label{fig: real_data_beta}
\end{center}
\end{figure}

To gain further understanding on the taxi traffic, we show the scaled $\wh\Q_2$ by MEFM and FM in Tables \ref{tab: real_data_Q2_MEFM} and \ref{tab: real_data_Q2_FM}, respectively. We can see that for the rush hours between 6 p.m. to 11 p.m., the estimated loadings almost vanish for MEFM, which is consistent with the fact that $\wh\bbeta_t$ captures the common hour effect on Manhattan life style. This also provides an intuition why the time-varying column/hour effect is strong, since in non-business days, the way that daily hours affecting the taxi traffic can change drastically over time as compared to the same when Manhattan zones are considered.
For demonstration purpose, we plot both $\wh{\beta}_{t,2}$ and $\wh{\beta}_{t,18}$ in Figure \ref{fig: real_data_beta}, where the former series features the mid-night effects and the latter features the night-life effects. Both series demonstrate obvious seasonality before COVID-19 as indicated on the plot.

The business-day series is also analysed, but since both $\textbf{reject}_\alpha$ and $\textbf{reject}_\beta$ are not significant, the estimated model is not shown here. The fact that the time-varying hour effect is not strong for business days is probably due to a rather routine working hours. Thus the hour effect is hardly changing and can be absorbed into a fixed mean, so that FM would be sufficient.

\begin{table}[h!]
\footnotesize
\setlength{\tabcolsep}{4pt}
\begin{center}
\begin{tabular}{l||ccccccccccccccccccccccccccc}
\hline $0_\text{am}$ & & 2 & & 4 & & 6 & & 8 & & 10 & & $12_\text{pm}$ & & 2 & & 4 & & 6 & & 8 & & 10 & & $12_\text{am}$ \\
\hline 1 & -2 & -5 & -6 & \red{-7} & \red{-7} & \red{-7} & -6 & -5 & -3 & 0 & 3 & 5 & 6 & 6 & 5 & 5 & 4 & 4 & 5 & 5 & 2 & 0 & 0 & -1 \\
\hline 2 & 6 & 5 & 3 & 1 & -1 & -4 & -5 & -6 & \red{-7} & \red{-7} & -6 & -5 & -3 & -2 & -1 & -2 & -1 & -1 & 2 & 5 & \red{8} & 6 & 6 & \red{7} \\
\hline 3 & -1 & \red{-13} & \red{-9} & -6 & -2 & 2 & 4 & 5 & 6 & 4 & 2 & -2 & -4 & -5 & -4 & -3 & -2 & -1 & 0 & 2 & 5 & 4 & 7 & \red{9} \\
\hline
\end{tabular}
\end{center}
\caption{Estimated loading matrix $\wh\Q_2$ using MEFM, after scaling. Magnitudes larger than 6 are highlighted in red.}
\label{tab: real_data_Q2_MEFM}
\end{table}

\begin{table}[h!]
\footnotesize
\setlength{\tabcolsep}{4pt}
\begin{center}
\begin{tabular}{l||ccccccccccccccccccccccccccc}
\hline $0_\text{am}$ & & 2 & & 4 & & 6 & & 8 & & 10 & & $12_\text{pm}$ & & 2 & & 4 & & 6 & & 8 & & 10 & & $12_\text{am}$ \\
\hline 1 & -5 & -5 & -4 & -3 & -2 & -1 & -1 & -1 & -2 & -3 & -4 & -5 & -5 & -6 & -5 & -5 & -5 & -5 & \red{-6} & \red{-6} & \red{-6} & -5 & -5 & -5 \\
\hline 2 & 5 & \red{7} & \red{7} & 5 & 4 & 2 & 0 & -2 & -4 & \red{-6} & \red{-6} & \red{-6} & \red{-6} & -5 & -4 & -4 & -3 & -3 & -2 & 1 & 4 & 4 & 5 & \red{6} \\
\hline 3 & 1 & \red{-13} & \red{-10} & \red{-9} & \red{-6} & -3 & -1 & 0 & 1 & 0 & -1 & -3 & -3 & -3 & -3 & -2 & -1 & 0 & 3 & \red{6} & \red{6} & \red{6} & \red{8} & \red{11} \\
\hline
\end{tabular}
\end{center}
\caption{Estimated loading matrix $\wh\Q_2$ using FM, after scaling. Magnitudes larger than 5 are highlighted in red.}
\label{tab: real_data_Q2_FM}
\end{table}

\section{Appendix: Proof of all theorems and lemmas}\label{sec:Appendix}

\subsection{Proof of Theorems}

\textbf{\textit{Proof of Theorem \ref{thm: identification}.}} Suppose we have another set of parameters, $(\wt{\mu}_t, \wt{\balpha}_t, \wt{\bbeta}_t, \wt{\Q}_r, \wt{\Q}_c, \wt{\F}_{Z,t})$ for $t\in[T]$, also satisfying (\ref{eqn: model}). For each $t\in[T]$, left-multiplying by $\1_p'$ and right-multiplying by $\1_q$ on (\ref{eqn: model}), we would arrive at $pq \wt{\mu}_t=pq \mu_t$ according to (IC1), so $\mu_t$ is identified. Similarly, $\wt\balpha_t = \balpha_t$ and $\wt\bbeta = \bbeta_t$, by separately left-multiplying $\1_p'$ and right-multiplying $\1_q$ on $\Y_t$.

We hence have $\wt\Q_r\wt\F_{Z,t}\wt\Q_c' = \Q_r\F_{Z,t}\Q_c'$, but the factor loading matrices and factor series require further identification due to the multiplicative form. Without loss of generality, write
\[
\wt\Q_r = \Q_r \M_r + \bGamma_r, \;\;\;
\text{where } \bGamma_r'\Q_r = \0,
\]
with $\M_r\in\b{R}^{k_r\times k_r}$ and $\bGamma_r\in\b{R}^{p\times k_r}$, but can have zero columns. Then we have
\[
\0=\bGamma_r'\Q_r\F_{Z,t}\Q_c'\wt\Q_c =
\bGamma_r'\wt\Q_r\wt\F_{Z,t}\wt\Q_c'\wt\Q_c =
\bGamma_r'\bGamma_r\wt\F_{Z,t}\wt\Q_c'\wt\Q_c,
\]
which can only be true in general if $\bGamma_r=\0$ since $\wt\F_t$ is random and $\wt\Q_c'\wt\Q_c\to\bSigma_{A,c}$ due to (L1). Using (L1), $\M_r$ is of full rank and hence $\wt\Q_r$ and $\Q_r$ share the same column space. Similarly, the factor loading space of $\Q_c$ is identified, and $\F_{Z,t}$ is hence identified once $\Q_r$ and $\Q_c$ are given correspondingly.
$\square$

\textbf{\textit{Proof of Theorem \ref{thm: consistency}.}} By Assumption (IC1), we have $\mu_t=\1_p '(\Y_t-\E_t) \1_q/pq$ and hence
\begin{equation}
\label{eqn: rate_mu}
(\wh\mu_t - \mu_t)^2 = \frac{1}{p^2q^2} \Big(\1_p' \E_t \1_q \Big)^2 = \frac{1}{p^2q^2} \Big(\sum_{i=1}^p\sum_{j=1}^q E_{t,ij}\Big)^2 .
\end{equation}
Assumption (E1) implies each entry of $\E_t$ has zero mean and bounded fourth moment, and we have
\begin{equation}
\label{eqn: 1E1_bound}
    \b{E}\bigg[\Big(\sum_{i=1}^p\sum_{j=1}^q E_{t,ij}\Big)^2 \bigg] =
    \text{Var}\Big(\sum_{i=1}^p\sum_{j=1}^q E_{t,ij}\Big) =
    \sum_{i=1}^p \sum_{l=1}^p \sum_{j=1}^q \sum_{h=1}^q\text{cov}(E_{t,ij}, E_{t,lh})
    = O(pq) ,
\end{equation}
where we used Lemma \ref{lemma: correlation_Et_Ft} in the last equality. Thus with (\ref{eqn: rate_mu}), $(\wh\mu_t - \mu_t)^2 = O_P(p^{-1} q^{-1})$.

Similar to the rate for $\wh\mu_t$, by again (IC1) we have $\balpha_t = q^{-1}\Y_t\1_q - \mu_t \1_p - q^{-1}\E_t\1_q$ and $\bbeta_t = p^{-1}\Y_t'\1_p - \mu_t \1_q - p^{-1}\E_t'\1_p$. Then we have
\begin{align}
    \frac{1}{p}\cdot\|\wh\balpha_t - \balpha_t\|^2 &=
    \frac{1}{p}\cdot\Big\|(\mu_t - \wh\mu_t) \1_p +
    q^{-1}\E_t\1_q \Big\|^2 ,
    \label{eqn: rate_alpha} \\
    \frac{1}{q}\cdot\|\wh\bbeta_t - \bbeta_t\|^2 &=
    \frac{1}{q}\cdot \Big\|(\mu_t - \wh\mu_t) \1_q + p^{-1}\E_t'\1_p \Big\|^2 .
    \label{eqn: rate_beta}
\end{align}
From (\ref{eqn: rate_mu}) we have $\|(\mu_t - \wh\mu_t) \1_p\|^2 = O_P(q^{-1})$ and $\|(\mu_t - \wh\mu_t) \1_q\|^2 = O_P(p^{-1})$. Furthermore, by Lemma \ref{lemma: correlation_Et_Ft},
\[
\b{E}\Big[\|q^{-1}\E_t\1_q\|^2 \Big] = q^{-2} \cdot \sum_{i=1}^p \text{Var}\Big(\sum_{j=1}^q E_{t,ij}\Big) = O(pq^{-1}).
\]
Similarly, $\|p^{-1}\E_t'\1_p\|^2 = O_P(qp^{-1})$. Then we have (\ref{eqn: rate_alpha}) and (\ref{eqn: rate_beta}) as
\begin{equation*}
    \frac{1}{p}\cdot\|\wh\balpha_t - \balpha_t\|^2 = O_P(q^{-1}) ,
    \;\;\;
    \frac{1}{q}\cdot\|\wh\bbeta_t - \bbeta_t\|^2 =
    O_P(p^{-1}) .
\end{equation*}

In the rest of the proof, we show consistency of the factor loading estimators. From (\ref{eqn: Ct_estimator}), we have
\begin{align*}
    \wh\L_t &= \Y_t + (pq)^{-1} \1_p'\Y_t\1_q \1_p\1_q' - q^{-1}\Y_t\1_q \1_q' -  p^{-1}\1_p\1_p'\Y_t, \\
    \wh\L_t' &= \Y_t' + (pq)^{-1} \1_q\1_p'\1_q'\Y_t' \1_p - q^{-1} \1_q\1_q'\Y_t' -  p^{-1}\Y_t'\1_p\1_p' ,
\end{align*}
and hence the following decomposition
\begin{equation}
\label{eqn: sample_covariance_row}
\begin{split}
    \wh\L_t\wh\L_t' &= \Y_t\Y_t' +
    (pq)^{-1}\Y_t \1_q\1_p'\1_q'\Y_t' \1_p -
    q^{-1}\Y_t\1_q\1_q'\Y_t' -
    p^{-1}\Y_t\Y_t'\1_p\1_p' \\
    &+
    (pq)^{-1}\1_p'\Y_t\1_q \1_p\1_q'\Y_t' +
    (pq)^{-2}\1_p'\Y_t\1_q \1_p\1_q'\1_q\1_p'\1_q'\Y_t' \1_p\\
    &-
    q^{-1} (pq)^{-1}\1_p'\Y_t\1_q \1_p\1_q'
    \1_q\1_q'\Y_t' -
    p^{-1}(pq)^{-1}\1_p'\Y_t\1_q \1_p\1_q'
    \Y_t'\1_p\1_p'\\
    &-
    q^{-1}\Y_t\1_q \1_q'\Y_t' -
    q^{-1}(pq)^{-1}\Y_t\1_q\1_q'\1_q\1_p'\1_q'\Y_t'\1_p+
    q^{-2}\Y_t\1_q\1_q'\1_q \1_q'\Y_t' +
    q^{-1}p^{-1}\Y_t\1_q \1_q' \Y_t'\1_p\1_p' \\
    &-
    p^{-1}\1_p\1_p'\Y_t \Y_t' -
    (pq)^{-1}p^{-1}\1_p\1_p'\Y_t
    \1_q\1_p'\1_q'\Y_t' \1_p +
    q^{-1}p^{-1}\1_p\1_p'\Y_t \1_q\1_q'\Y_t' +
    p^{-2}\1_p\1_p'\Y_t\Y_t'\1_p\1_p' \\
    &=
    \Y_t\Y_t' +
    (pq)^{-1}\1_q'\Y_t' \1_p\Y_t \1_q\1_p' -
    p^{-1}\Y_t\Y_t'\1_p\1_p' -
    q^{-1}\Y_t\1_q \1_q'\Y_t' -
    p^{-1}\1_p\1_p'\Y_t \Y_t' \\
    &-
    (pq)^{-1}p^{-1}\1_p'\Y_t
    \1_q\1_q'\Y_t' \1_p\1_p\1_p' +
    (pq)^{-1}\1_p'\Y_t \1_q\1_p\1_q'\Y_t' +
    p^{-2}\1_p'\Y_t\Y_t'\1_p\1_p\1_p' \\
    &=
    \Y_t\Y_t' + \cQ_1 - \cQ_2
    - \cQ_3 - \cQ_4 - \cQ_5 + \cQ_6 + \cQ_7 ,
    \;\;\; \text{where} \\
    \cQ_1 &:= (pq)^{-1}\1_q'\Y_t' \1_p\Y_t \1_q\1_p',
    \;\;\;
    \cQ_2 := p^{-1}\Y_t\Y_t'\1_p\1_p' ,
    \;\;\;
    \cQ_3 := q^{-1}\Y_t\1_q \1_q'\Y_t' ,
    \;\;\;
    \cQ_4 := p^{-1}\1_p\1_p'\Y_t \Y_t', \\
    \cQ_5 &:= (pq)^{-1}p^{-1}\1_p'\Y_t
    \1_q\1_q'\Y_t' \1_p\1_p\1_p',
    \;\;\;
    \cQ_6 := (pq)^{-1}\1_p'\Y_t \1_q\1_p\1_q'\Y_t',
    \;\;\;
    \cQ_7 := p^{-2}\1_p'\Y_t\Y_t'\1_p\1_p\1_p'.
\end{split}
\end{equation}
From (\ref{eqn: model}), we have $\1_q'\Y_t'\1_p = pq\mu_t +\1_q'\E_t'\1_p$ and $\Y_t\1_q\1_p' = q\mu_t \1_p\1_p' +q\balpha_t\1_p' +\E_t\1_q\1_p'$. Thus,
\begin{equation}
\label{eqn: cQ1}
\begin{split}
    \cQ_1 &=
     (pq)^{-1}(pq^2\mu_t^2\1_p\1_p' +
     pq^2\mu_t\balpha_t\1_p' +
     pq\mu_t\E_t\1_q\1_p'
     + q\mu_t\1_q'\E_t'\1_p\1_p\1_p' +
     q\1_q'\E_t'\1_p\balpha_t\1_p' +
     \1_q'\E_t'\1_p\E_t\1_q\1_p') \\
     &= q\mu_t^2\1_p\1_p' +
     q\mu_t\balpha_t\1_p' +
     \mu_t\E_t\1_q\1_p'
     + p^{-1}\mu_t\1_q'\E_t'\1_p\1_p\1_p' +
     p^{-1}\1_q'\E_t'\1_p\balpha_t\1_p' +
     (pq)^{-1}\1_q'\E_t'\1_p\E_t\1_q\1_p' .
\end{split}
\end{equation}
Similarly, we have $\Y_t \Y_t' = (\mu_t \1_p\1_q' + \balpha_t\1_q' + \1_p\bbeta_t'  +\C_t +\E_t)(\mu_t \1_q\1_p' + \1_q\balpha_t' + \bbeta_t\1_p'  +\C_t' +\E_t')$. Further with Assumption (IC1),
\begin{equation}
\label{eqn: cQ2}
\begin{split}
    \cQ_2 &= p^{-1}(
     \mu_t^2\1_p\1_q'\1_q\1_p' + \mu_t\1_p\1_q'\1_q\balpha_t'
     + \mu_t\1_p\1_q'\bbeta_t\1_p' + \mu_t\1_p\1_q'\C_t'
     + \mu_t\1_p\1_q'\E_t' \\
     &+ \mu_t\balpha_t\1_q'\1_q\1_p' + \balpha_t\1_q'\1_q\balpha_t' + \balpha_t\1_q'\bbeta_t\1_p' + \balpha_t\1_q'\C_t' +
     \balpha_t\1_q'\E_t' \\
     &+ \mu_t\1_p\bbeta_t'\1_q\1_p' + \1_p\bbeta'\1_q\balpha_t' + \1_p\bbeta_t'\bbeta_t\1_p'
     + \1_p\bbeta_t'\C_t' + \1_p\bbeta_t'\E_t' \\
     &+ \mu_t\C_t\1_q\1_p' + \C_t\1_q\balpha_t' + \C_t\bbeta_t \1_p' + \C_t\C_t' + \C_t\E_t' \\
     &+ \mu_t\E_t\1_q\1_p' + \E_t\1_q\balpha_t' +
     \E_t\bbeta_t\1_p' + \E_t\C_t' + \E_t\E_t'
     )\1_p\1_p' \\
     &=
     q\mu_t^2\1_p\1_p'
     + p^{-1}\mu_t\1_q'\E_t'\1_p\1_p\1_p'
     + q\mu_t\balpha_t\1_p'
     + p^{-1}\1_q'\E_t'\1_p\balpha_t\1_p'
     +\1_p\bbeta_t'\bbeta_t\1_p' + p^{-1}\bbeta_t'\E_t' \1_p\1_p\1_p'   \\
     &+ \C_t\bbeta_t \1_p' + p^{-1}\C_t\E_t'\1_p\1_p'
     + \mu_t\E_t\1_q\1_p' + \E_t\bbeta_t\1_p' + p^{-1}\E_t\E_t'\1_p\1_p' .
\end{split}
\end{equation}
Since $\Y_t\1_q = q\mu_t \1_p + q\balpha_t +\E_t\1_q$, we have
\begin{equation}
\label{eqn: cQ3}
\begin{split}
    \cQ_3 &= q^{-1}(q\mu_t \1_p + q\balpha_t +\E_t\1_q)
    (q\mu_t \1_p' + q\balpha_t' +\1_q'\E_t') \\
    &=
    q\mu_t^2\1_p\1_p' + q\mu_t\1_p\balpha_t' +
    \mu_t\1_p\1_q'\E_t' + q\mu_t\balpha_t\1_p' +
    q\balpha_t\balpha_t' + \balpha_t\1_q'\E_t' +
    \mu_t\E_t\1_q\1_p' + \E_t\1_q\balpha_t' +
    q^{-1}\E_t\1_q\1_q'\E_t' .
\end{split}
\end{equation}
Similar to (\ref{eqn: cQ2}), we have
\begin{equation}
\label{eqn: cQ4}
\begin{split}
    \cQ_4 &=
     q\mu_t^2\1_p\1_p' + q\mu_t\1_p\balpha_t'
     + \mu_t\1_p\1_q'\bbeta_t\1_p' + \mu_t\1_p\1_q'\C_t'
     + \mu_t\1_p\1_q'\E_t' + \mu_t\1_p\1_p'\balpha_t\1_q'\1_q\1_p' \\
     &+ p^{-1}\1_p\1_p'\balpha_t\1_q'\1_q\balpha_t' + p^{-1}\1_p\1_p'\balpha_t\1_q'\bbeta_t\1_p' + p^{-1}\1_p\1_p'\balpha_t\1_q'\C_t' +
     p^{-1}\1_p\1_p'\balpha_t\1_q'\E_t' \\
     &+ \mu_t\1_p\bbeta_t'\1_q\1_p' + \1_p\bbeta'\1_q\balpha_t' + \1_p\bbeta_t'\bbeta_t\1_p'
     + \1_p\bbeta_t'\C_t' + \1_p\bbeta_t'\E_t'
     + p^{-1}\mu_t\1_p\1_p'\C_t\1_q\1_p' \\
     &+ p^{-1}\1_p\1_p'\C_t\1_q\balpha_t' + p^{-1}\1_p\1_p'\C_t\bbeta_t \1_p' + p^{-1}\1_p\1_p'\C_t\C_t' + p^{-1}\1_p\1_p'\C_t\E_t' \\
     &+ p^{-1}\mu_t\1_p\1_p'\E_t\1_q\1_p' + p^{-1}\1_p\1_p'\E_t\1_q\balpha_t' +
     p^{-1}\1_p\1_p'\E_t\bbeta_t\1_p' + p^{-1}\1_p\1_p'\E_t\C_t' + p^{-1}\1_p\1_p'\E_t\E_t' \\
     &=
     q\mu_t^2\1_p\1_p' + q\mu_t\1_p\balpha_t'
     + \mu_t\1_p\1_q'\E_t' + \bbeta_t'\bbeta_t\1_p\1_p'
     + \1_p\bbeta_t'\C_t' + \1_p\bbeta_t'\E_t' \\
     &+ p^{-1}\mu_t\1_p'\E_t\1_q\1_p\1_p' + p^{-1}\1_p'\E_t\1_q\1_p\balpha_t' +
     p^{-1}\1_p'\E_t\bbeta_t\1_p\1_p' + p^{-1}\1_p\1_p'\E_t\C_t' + p^{-1}\1_p\1_p'\E_t\E_t' ,
\end{split}
\end{equation}
where the last equality used Assumption (IC1). For $\cQ_5$ and $\cQ_6$, we have
\begin{equation}
\label{eqn: cQ5}
\begin{split}
    \cQ_5 &=
     (pq)^{-1}p^{-1}(pq\mu_t +\1_p'\E_t\1_q)(pq\mu_t +\1_q'\E_t'\1_p)\1_p\1_p' \\
     &= q\mu_t^2\1_p\1_p' + p^{-1}\mu_t\1_q'\E_t'\1_p\1_p\1_p' +
     p^{-1}\mu_t\1_p'\E_t\1_q\1_p\1_p' + (pq)^{-1}p^{-1}\1_p'\E_t\1_q\1_q'\E_t'\1_p\1_p\1_p' ,
\end{split}
\end{equation}
and
\begin{equation}
\label{eqn: cQ6}
\begin{split}
    \cQ_6 &=
     (\mu_t + (pq)^{-1}\1_q'\E_t'\1_p)
     (q\mu_t \1_p\1_p' + q\1_p\balpha_t' + \1_p\1_q'\bbeta_t\1_p'  +\1_p\1_q'\C_t' +\1_p\1_q'\E_t') \\
     &= q\mu_t^2\1_p\1_p' + q\mu_t\1_p\balpha_t' +
     \mu_t\1_p\1_q'\E_t' + p^{-1}\mu_t\1_q'\E_t'\1_p\1_p\1_p'
     + p^{-1}\1_q'\E_t'\1_p\1_p\balpha_t' +
     (pq)^{-1}\1_q'\E_t'\1_p\1_p\1_q'\E_t' .
\end{split}
\end{equation}
Lastly for $\cQ_7$, we have similar to (\ref{eqn: cQ2}) that
\begin{equation}
\label{eqn: cQ7}
\begin{split}
    \cQ_7 &= p^{-2}\1_p'(
     \mu_t^2\1_p\1_q'\1_q\1_p' + \mu_t\1_p\1_q'\1_q\balpha_t'
     + \mu_t\1_p\1_q'\bbeta_t\1_p' + \mu_t\1_p\1_q'\C_t'
     + \mu_t\1_p\1_q'\E_t' \\
     &+ \mu_t\balpha_t\1_q'\1_q\1_p' + \balpha_t\1_q'\1_q\balpha_t' + \balpha_t\1_q'\bbeta_t\1_p' + \balpha_t\1_q'\C_t' +
     \balpha_t\1_q'\E_t' \\
     &+ \mu_t\1_p\bbeta_t'\1_q\1_p' + \1_p\bbeta'\1_q\balpha_t' + \1_p\bbeta_t'\bbeta_t\1_p'
     + \1_p\bbeta_t'\C_t' + \1_p\bbeta_t'\E_t' \\
     &+ \mu_t\C_t\1_q\1_p' + \C_t\1_q\balpha_t' + \C_t\bbeta_t \1_p' + \C_t\C_t' + \C_t\E_t' \\
     &+ \mu_t\E_t\1_q\1_p' + \E_t\1_q\balpha_t' +
     \E_t\bbeta_t\1_p' + \E_t\C_t' + \E_t\E_t'
     )\1_p\1_p\1_p' \\
     &=
     p^{-2}(
     pq\mu_t^2\1_p' + pq\mu_t\balpha_t'
     + p\mu_t\1_q'\E_t' + p\bbeta_t'\bbeta_t\1_p'
     + p\bbeta_t'\C_t' + p\bbeta_t'\E_t' \\
     &+ \mu_t\1_p'\E_t\1_q\1_p' + \1_p'\E_t\1_q\balpha_t' +
     \1_p'\E_t\bbeta_t\1_p' + \1_p'\E_t\C_t' + \1_p'\E_t\E_t'
     )\1_p\1_p\1_p' \\
     &=
     q\mu_t^2\1_p\1_p'
     + p^{-1}\mu_t\1_q'\E_t'\1_p\1_p\1_p' + \bbeta_t'\bbeta_t\1_p\1_p'
    + p^{-1}\bbeta_t'\E_t'\1_p\1_p\1_p' \\
     &+ p^{-1}\mu_t\1_p'\E_t\1_q\1_p\1_p' +
     p^{-1}\1_p'\E_t\bbeta_t\1_p\1_p' + p^{-2}\1_p'\E_t\E_t'\1_p\1_p\1_p' .
\end{split}
\end{equation}

With (\ref{eqn: cQ1}), (\ref{eqn: cQ2}), (\ref{eqn: cQ3}), (\ref{eqn: cQ4}), (\ref{eqn: cQ5}), (\ref{eqn: cQ6}), and (\ref{eqn: cQ7}), we have
\begin{equation}
\label{eqn: sample_covariance_row_final}
\begin{split}
    \wh\L_t\wh\L_t' &= \Y_t\Y_t' + \cQ_1 - \cQ_2
    - \cQ_3 - \cQ_4 - \cQ_5 + \cQ_6 + \cQ_7 \\
    &=
    q\mu_t^2\1_p\1_p' + q\mu_t\1_p\balpha_t'
     + \mu_t\1_p\1_q'\E_t'
     + q\mu_t\balpha_t\1_p' + q\balpha_t\balpha_t' +
     \balpha_t\1_q'\E_t' + \1_p\bbeta_t'\bbeta_t\1_p'
     + \1_p\bbeta_t'\C_t' \\
     &+ \1_p\bbeta_t'\E_t'
     + \C_t\bbeta_t \1_p' + \C_t\C_t' + \C_t\E_t'
     + \mu_t\E_t\1_q\1_p' + \E_t\1_q\balpha_t' +
     \E_t\bbeta_t\1_p' + \E_t\C_t' + \E_t\E_t' \\
     &+ q\mu_t^2\1_p\1_p' +
     q\mu_t\balpha_t\1_p' +
     \mu_t\E_t\1_q\1_p' +
     p^{-1}\mu_t\1_q'\E_t'\1_p\1_p\1_p' +
     p^{-1}\1_q'\E_t'\1_p\balpha_t\1_p' +
     (pq)^{-1}\1_q'\E_t'\1_p\E_t\1_q\1_p' \\
     &- q\mu_t^2\1_p\1_p'
     - p^{-1}\mu_t\1_q'\E_t'\1_p\1_p\1_p'
     - q\mu_t\balpha_t\1_p'
     - p^{-1}\1_q'\E_t'\1_p\balpha_t\1_p'
     -\1_p\bbeta_t'\bbeta_t\1_p' - p^{-1}\bbeta_t'\E_t' \1_p\1_p\1_p' \\
     &- \C_t\bbeta_t \1_p'
     - p^{-1}\C_t\E_t'\1_p\1_p'
     - \mu_t\E_t\1_q\1_p' - \E_t\bbeta_t\1_p' - p^{-1}\E_t\E_t'\1_p\1_p'
     - q\mu_t^2\1_p\1_p' - q\mu_t\1_p\balpha_t' \\
     &- \mu_t\1_p\1_q'\E_t' - q\mu_t\balpha_t\1_p' -
    q\balpha_t\balpha_t' - \balpha_t\1_q'\E_t' -
    \mu_t\E_t\1_q\1_p' - \E_t\1_q\balpha_t' -
    q^{-1}\E_t\1_q\1_q'\E_t' \\
     &- q\mu_t^2\1_p\1_p' - q\mu_t\1_p\balpha_t'
     - \mu_t\1_p\1_q'\E_t' - \bbeta_t'\bbeta_t\1_p\1_p'
     - \1_p\bbeta_t'\C_t' - \1_p\bbeta_t'\E_t'
     - p^{-1}\mu_t\1_p'\E_t\1_q\1_p\1_p' \\
     &- p^{-1}\1_p'\E_t\1_q\1_p\balpha_t' -
     p^{-1}\1_p'\E_t\bbeta_t\1_p\1_p' - p^{-1}\1_p\1_p'\E_t\C_t' - p^{-1}\1_p\1_p'\E_t\E_t'
     - q\mu_t^2\1_p\1_p' \\
     &- p^{-1}\mu_t\1_q'\E_t'\1_p\1_p\1_p' -
     p^{-1}\mu_t\1_p'\E_t\1_q\1_p\1_p' - (pq)^{-1}p^{-1}\1_p'\E_t\1_q\1_q'\E_t'\1_p\1_p\1_p'
     + q\mu_t^2\1_p\1_p' \\
     &+ q\mu_t\1_p\balpha_t' +
     \mu_t\1_p\1_q'\E_t' + p^{-1}\mu_t\1_q'\E_t'\1_p\1_p\1_p'
     + p^{-1}\1_q'\E_t'\1_p\1_p\balpha_t' +
     (pq)^{-1}\1_q'\E_t'\1_p\1_p\1_q'\E_t' \\
     &+ q\mu_t^2\1_p\1_p'
     + p^{-1}\mu_t\1_q'\E_t'\1_p\1_p\1_p' + \bbeta_t'\bbeta_t\1_p\1_p'
     + p^{-1}\bbeta_t'\E_t'\1_p\1_p\1_p'
     + p^{-1}\mu_t\1_p'\E_t\1_q\1_p\1_p' \\
     &+
     p^{-1}\1_p'\E_t\bbeta_t\1_p\1_p' + p^{-2}\1_p'\E_t\E_t'\1_p\1_p\1_p' \\
     &=
     \C_t\C_t' + \C_t\E_t' + \E_t\C_t' + \E_t\E_t' + (pq)^{-1}\1_q'\E_t'\1_p\E_t\1_q\1_p'
     + (pq)^{-1}\1_q'\E_t'\1_p\1_p\1_q'\E_t' \\
     & - p^{-1}\C_t\E_t'\1_p\1_p' - p^{-1}\1_p\1_p'\E_t\C_t'
     - p^{-1}\E_t\E_t'\1_p\1_p' - p^{-1}\1_p\1_p'\E_t\E_t' - p^{-1}\E_t\1_q\1_q'\E_t' \\
     &- (pq)^{-1}p^{-1} (\1_q'\E_t'\1_p)^2\1_p\1_p' + p^{-2}\1_p'\E_t\E_t'\1_p\1_p\1_p' .
\end{split}
\end{equation}
Swapping the roles of row and column factor loadings, we can arrive at similarly
\begin{equation}
\label{eqn: sample_covariance_column_final}
\begin{split}
    \wh\L_t'\wh\L_t &=
     \C_t'\C_t + \C_t'\E_t + \E_t'\C_t + \E_t'\E_t + (pq)^{-1}\1_p'\E_t\1_q\E_t'\1_p\1_q'
     + (pq)^{-1}\1_p'\E_t\1_q\1_q\1_p'\E_t \\
     & - q^{-1}\C_t'\E_t\1_q\1_q' - q^{-1}\1_q\1_q'\E_t'\C_t
     - q^{-1}\E_t'\E_t\1_q\1_q' - q^{-1}\1_q\1_q'\E_t'\E_t - q^{-1}\E_t'\1_p\1_p'\E_t \\
     &- (pq)^{-1}q^{-1} (\1_p'\E_t\1_q)^2\1_q\1_q' + q^{-2}\1_q'\E_t'\E_t\1_q\1_q\1_q' .
\end{split}
\end{equation}

For ease of notation, we define
\begin{align}
    \R_{r,t} &:= \C_t\E_t' + \E_t\C_t' + \E_t\E_t' + (pq)^{-1}\1_q'\E_t'\1_p\E_t\1_q\1_p'
     + (pq)^{-1}\1_q'\E_t'\1_p\1_p\1_q'\E_t'  \notag \\
     & - p^{-1}\C_t\E_t'\1_p\1_p' - p^{-1}\1_p\1_p'\E_t\C_t'
     - p^{-1}\E_t\E_t'\1_p\1_p' - p^{-1}\1_p\1_p'\E_t\E_t' - p^{-1}\E_t\1_q\1_q'\E_t' \notag \\
     &- (pq)^{-1}p^{-1} (\1_q'\E_t'\1_p)^2\1_p\1_p' + p^{-2}\1_p'\E_t\E_t'\1_p\1_p\1_p'
     , \label{eqn: R_rt} \\
     \R_{c,t} &:= \C_t'\E_t + \E_t'\C_t + \E_t'\E_t + (pq)^{-1}\1_p'\E_t\1_q\E_t'\1_p\1_q'
     + (pq)^{-1}\1_p'\E_t\1_q\1_q\1_p'\E_t \notag \\
     & - q^{-1}\C_t'\E_t\1_q\1_q' - q^{-1}\1_q\1_q'\E_t'\C_t
     - q^{-1}\E_t'\E_t\1_q\1_q' - q^{-1}\1_q\1_q'\E_t'\E_t - q^{-1}\E_t'\1_p\1_p'\E_t \notag \\
     &- (pq)^{-1}q^{-1} (\1_p'\E_t\1_q)^2\1_q\1_q' + q^{-2}\1_q'\E_t'\E_t\1_q\1_q\1_q'
     ,\label{eqn: R_ct}
\end{align}
so that from (\ref{eqn: sample_covariance_row_final}) and (\ref{eqn: sample_covariance_column_final}), we can write
\[
\wh\L_t\wh\L_t' = \C_t\C_t' + \R_{r,t} ,
\;\;\;
\wh\L_t'\wh\L_t = \C_t'\C_t + \R_{c,t} .
\]
Recall that we denote by $\wh\D_r$ the $k_r \times k_r$ diagonal matrix with the first largest $k_r$ eigenvalues of $T^{-1} \sum_{t=1}^T \wh\L_t\wh\L_t'$ on the main diagonal, and since $\wh\Q_r$ consists of the corresponding eigenvectors, we have
\begin{equation}
\label{eqn: row_svd}
    \wh\Q_r\wh\D_r =
    T^{-1}\sum_{t=1}^T \wh\L_t\wh\L_t'\wh\Q_r .
\end{equation}
With (\ref{eqn: sample_covariance_row_final}) and $\C_t\C_t' = \Q_r\F_{Z,t}\Q_c'\Q_c\F_{Z,t}'\Q_r'$, we can write the $j$-th row of estimated row factor loading as
\begin{equation*}
\begin{split}
    \wh\Q_{r,j\cdot} &= T^{-1}\wh\D_r^{-1} \sum_{i=1}^p\wh\Q_{r,i\cdot} \sum_{t=1}^T (\wh\L_t\wh\L_t')_{ij} \\
    &= T^{-1}\wh\D_r^{-1} \sum_{i=1}^p\wh\Q_{r,i\cdot} \Q_{r,i\cdot}' \sum_{t=1}^T (\F_{Z,t}\Q_c'\Q_c\F_{Z,t}')\Q_{r,j\cdot}
    + T^{-1}\wh\D_r^{-1} \sum_{i=1}^p \wh\Q_{r,i\cdot} \sum_{t=1}^T ( \R_{r,t} )_{ij}.
\end{split}
\end{equation*}
Thus with the definition $\H_r = T^{-1}\wh\D_r^{-1} \wh\Q_r' \Q_r \sum_{t=1}^T (\F_{Z,t} \Q_c' \Q_c \F_{Z,t}')$, we have
\begin{equation*}
    \wh\Q_{r,j\cdot} - \H_r\Q_{r,j\cdot} =
    T^{-1}\wh\D_r^{-1} \sum_{i=1}^p
    \wh\Q_{r,i\cdot} \sum_{t=1}^T (\R_{r,t})_{ij} ,
\end{equation*}
and hence we have
\begin{equation*}
\begin{split}
    &\hspace{5mm}
    \|\wh\Q_r - \Q_r\H_r' \|_F^2 =
    \sum_{j=1}^p \|\wh\Q_{r,j\cdot} - \H_r\Q_{r,j\cdot} \|^2 =
    \sum_{j=1}^p \Big\|T^{-1}\wh\D_r^{-1} \wh\Q_r' \Big(\sum_{t=1}^T \R_{r,t}\Big)_{\cdot j} \Big\|^2 \\
    &\leq
    T^{-2} \cdot \|\wh\D_r^{-1}\|_F^2 \cdot
    \|\wh\Q_r\|_F^2 \cdot \Big\| \sum_{t=1}^T \R_{r,t} \Big\|_F^2
    = O_P\Big(T^{-1}p^{2(1-\delta_{r,k_r})}q^{1-2\delta_{c,1}} + p^{1-2\delta_{r,k_r}}q^{2(1-\delta_{c,1})} \Big),
\end{split}
\end{equation*}
where the last equality used Lemma \ref{lemma: rate_R_rt} and Lemma \ref{lemma: norm_hat_D}. The consistency of $\wh\Q_c$ can be similarly shown (omitted here). This completes the proof of Theorem \ref{thm: consistency}. $\square$

\textbf{\textit{Proof of Theorem \ref{thm: consistency_F_C}.}}
From (\ref{eqn: Ct_estimator}), we can first write
\begin{equation*}
\begin{split}
    \wh\L_t &= \Y_t + (pq)^{-1} \1_p'\Y_t\1_q \1_p\1_q' - q^{-1}\Y_t\1_q \1_q' -  p^{-1}\1_p\1_p'\Y_t \\
    &=
    \mu_t \1_p\1_q' + \balpha_t\1_q' + \1_p\bbeta_t'  +\Q_r\F_{Z,t}\Q_c' +\E_t \\
    &+ (pq)^{-1} \1_p'(\mu_t \1_p\1_q' + \balpha_t\1_q' + \1_p\bbeta_t'  +\Q_r\F_{Z,t}\Q_c' +\E_t )\1_q \1_p\1_q' \\
    &- q^{-1}(\mu_t \1_p\1_q' + \balpha_t\1_q' + \1_p\bbeta_t'  +\Q_r\F_{Z,t}\Q_c' +\E_t )\1_q \1_q' \\
    &- p^{-1}\1_p\1_p'(\mu_t \1_p\1_q' + \balpha_t\1_q' + \1_p\bbeta_t'  +\Q_r\F_{Z,t}\Q_c' +\E_t ) \\
    &=
    \mu_t \1_p\1_q' + \balpha_t\1_q' + \1_p\bbeta_t'  +\Q_r\F_{Z,t}\Q_c' +\E_t + (pq)^{-1} \1_p'\mu_t \1_p\1_q'\1_q \1_p\1_q' + (pq)^{-1} \1_p' \balpha_t\1_q'\1_q \1_p\1_q' \\
    &+ (pq)^{-1} \1_p'\1_p\bbeta_t'\1_q \1_p\1_q' + (pq)^{-1} \1_p' \Q_r\F_{Z,t}\Q_c'\1_q \1_p\1_q' + (pq)^{-1} \1_p'\E_t \1_q \1_p\1_q' \\
    &- q^{-1}\mu_t \1_p\1_q'\1_q \1_q' - q^{-1}\balpha_t\1_q'\1_q \1_q' - q^{-1}\1_p\bbeta_t'\1_q \1_q' - q^{-1} \Q_r\F_{Z,t}\Q_c' \1_q \1_q' - q^{-1}\E_t\1_q \1_q' \\
    &- p^{-1}\1_p\1_p'\mu_t \1_p\1_q' - p^{-1}\1_p\1_p'\balpha_t\1_q' -p^{-1}\1_p\1_p'\1_p\bbeta_t' -p^{-1}\1_p\1_p'\Q_r\F_{Z,t}\Q_c' -p^{-1}\1_p\1_p'\E_t \\
    &=
    \Q_r\F_{Z,t}\Q_c' +\E_t + (pq)^{-1} \1_p'\E_t \1_q \1_p\1_q' - q^{-1}\E_t\1_q \1_q' -p^{-1}\1_p\1_p'\E_t,
\end{split}
\end{equation*}
where the last equality used Assumption (IC1). Thus, we have
\begin{equation}
\label{eqn: F_Zt_decomp}
\begin{split}
    &\hspace{5mm}
    \wh\F_{Z,t} - (\H_r^{-1})' \F_{Z,t} \H_c^{-1} =
    \wh\Q_r' \wh\L_t \wh\Q_c - (\H_r^{-1})' \F_{Z,t} \H_c^{-1} \\
    &=
    \wh\Q_r'(\Q_r\H_r')(\H_r^{-1})' \F_{Z,t} \H_c^{-1} (\Q_c\H_c')'\wh\Q_c - (\H_r^{-1})' \F_{Z,t} \H_c^{-1} + \wh\Q_r'\E_t\wh\Q_c \\
    & + (pq)^{-1} \wh\Q_r'\1_p'\E_t \1_q \1_p\1_q'\wh\Q_c - q^{-1}\wh\Q_r'\E_t\1_q \1_q'\wh\Q_c - p^{-1}\wh\Q_r'\1_p\1_p'\E_t\wh\Q_c \\
    &=
    \wh\Q_r'(\Q_r\H_r' - \wh\Q_r)(\H_r^{-1})' \F_{Z,t} \H_c^{-1} (\Q_c\H_c' - \wh\Q_c)'\wh\Q_c + \wh\Q_r'(\Q_r\H_r' - \wh\Q_r)(\H_r^{-1})' \F_{Z,t} \H_c^{-1} \\
    & + (\H_r^{-1})' \F_{Z,t} \H_c^{-1} (\Q_c\H_c' - \wh\Q_c)'\wh\Q_c + (\wh\Q_r-\Q_r\H_r)'\E_t(\wh\Q_c - \Q_c\H_c') \\
    & + (\wh\Q_r - \Q_r\H_r)'\E_t \Q_c\H_c' + \H_r'\Q_r'\E_t(\wh\Q_c - \Q_c\H_c') + \H_r'\Q_r'\E_t\Q_c\H_c' \\
    & + (pq)^{-1} \wh\Q_r'\1_p'\E_t \1_q \1_p\1_q'\wh\Q_c - q^{-1}(\wh\Q_r - \Q_r\H_r')'\E_t\1_q \1_q'\wh\Q_c - q^{-1}\H_r\Q_r'\E_t\1_q \1_q'\wh\Q_c \\
    & - p^{-1}\wh\Q_r'\1_p\1_p'\E_t(\wh\Q_c - \Q_c\H_c') - p^{-1}\wh\Q_r'\1_p\1_p'\E_t\Q_c\H_c' \\
    &=: \cI_{F,1} + \cI_{F,2}+ \cI_{F,3}+ \cI_{F,4}+ \cI_{F,5}+ \cI_{F,6}+ \cI_{F,7}+ \cI_{F,8} \\
    &- \cI_{F,9}- \cI_{F,10}- \cI_{F,11}- \cI_{F,12} ,
    \;\;\; \text{where}
\end{split}
\end{equation}
\begin{align*}
    & \cI_{F,1} := \wh\Q_r'(\Q_r\H_r' - \wh\Q_r)(\H_r^{-1})' \F_{Z,t} \H_c^{-1} (\Q_c\H_c' - \wh\Q_c)'\wh\Q_c , \\
    & \cI_{F,2} := \wh\Q_r'(\Q_r\H_r' - \wh\Q_r)(\H_r^{-1})' \F_{Z,t} \H_c^{-1}, \;\;\;
    \cI_{F,3} := (\H_r^{-1})' \F_{Z,t} \H_c^{-1} (\Q_c\H_c' - \wh\Q_c)'\wh\Q_c, \\
    & \cI_{F,4} := (\wh\Q_r-\Q_r\H_r)'\E_t(\wh\Q_c - \Q_c\H_c'), \;\;\;
    \cI_{F,5} := (\wh\Q_r - \Q_r\H_r)'\E_t \Q_c\H_c', \\
    & \cI_{F,6} := \H_r'\Q_r'\E_t(\wh\Q_c - \Q_c\H_c'), \;\;\;
    \cI_{F,7} := \H_r'\Q_r'\E_t\Q_c\H_c' , \;\;\;
    \cI_{F,8} := (pq)^{-1} \wh\Q_r'\1_p'\E_t \1_q \1_p\1_q'\wh\Q_c , \\
    & \cI_{F,9} := q^{-1}(\wh\Q_r - \Q_r\H_r')'\E_t\1_q \1_q'\wh\Q_c , \;\;\;
    \cI_{F,10} := q^{-1}\H_r\Q_r'\E_t\1_q \1_q'\wh\Q_c , \\
    & \cI_{F,11} := p^{-1}\wh\Q_r'\1_p\1_p'\E_t(\wh\Q_c - \Q_c\H_c'), \;\;\;
    \cI_{F,12} := p^{-1}\wh\Q_r'\1_p\1_p'\E_t\Q_c\H_c' .
\end{align*}

First consider $\F_{Z,t}$, by its definition and Assumption (L1) we have
\[
\|\F_{Z,t}\|_F^2 \leq \|\F_t\|_F^2 \cdot \|\Z_r^{1/2}\|_F^2 \cdot \|\Z_c^{1/2}\|_F^2 = O_P(p^{\delta_{r,1}} q^{\delta_{c,1}}).
\]
Then for $\cI_{F,1}$, we have
\begin{equation*}
\begin{split}
    \|\cI_{F,1} \|_F^2 &= O_P(p^{\delta_{r,1}} q^{\delta_{c,1}}) \cdot \|\Q_r\H_r' - \wh\Q_r\|_F^2 \cdot \|\Q_c\H_c' - \wh\Q_c\|_F^2 \\
    &= O_P\big( T^{-2} p^{2  - 2\delta_{r,k_r}} q^{2 - 2\delta_{c,k_c}}
    + T^{-1} p^{3  - 2\delta_{r,k_r}} q^{2-\delta_{c,1}-3\delta_{c,k_c}} \\
    &\hspace{12pt}
    + T^{-1} q^{3 - 2\delta_{c,k_c}} p^{2-\delta_{r,1}-3\delta_{r,k_r}}
    + p^{3 -\delta_{r,1} -3\delta_{r,k_r}} q^{3-\delta_{c,1}-3\delta_{c,k_c}} \big) ,
\end{split}
\end{equation*}
where we used Lemma \ref{lemma: corollary_asymp_loading} in the last equality. Similarly, we can obtain the following for $\cI_{F,2}$ and $\cI_{F,3}$,
\begin{align*}
    \|\cI_{F,2} \|_F^2 &= O_P\big(T^{-1} p^{1+ 2\delta_{r,1} - 2\delta_{r,k_r}} q^{1-\delta_{c,1}} +  p^{1 +\delta_{r,1} -3\delta_{r,k_r}} q^{2-\delta_{c,1}} \big), \\
    \|\cI_{F,3} \|_F^2 &= O_P\big(T^{-1} q^{1+ 2\delta_{c,1} - 2\delta_{c,k_c}} p^{1-\delta_{r,1}} +  q^{1 +\delta_{c,1} -3\delta_{c,k_c}} p^{2-\delta_{r,1}} \big) .
\end{align*}

For $\cI_{F,4}$, from Assumptions (E1) and (E2) we easily have $\|\E_t\|_F^2 = O(pq)$, so we have
\begin{equation*}
\begin{split}
    \|\cI_{F,4} \|_F^2 &= O_P(pq) \cdot \|\Q_r\H_r' - \wh\Q_r\|_F^2 \cdot \|\Q_c\H_c' - \wh\Q_c\|_F^2 \\
    &= O_P\big( T^{-2} p^{3 -\delta_{r,1} - 2\delta_{r,k_r}} q^{3 -\delta_{c,1} - 2\delta_{c,k_c}}
    + T^{-1} p^{4 -\delta_{r,1} - 2\delta_{r,k_r}} q^{3-2\delta_{c,1}-3\delta_{c,k_c}} \\
    &\hspace{12pt}
    + T^{-1} q^{4 -\delta_{c,1} - 2\delta_{c,k_c}} p^{3-2\delta_{r,1}-3\delta_{r,k_r}}
    + p^{4 -2\delta_{r,1} -3\delta_{r,k_r}} q^{4-2\delta_{c,1}-3\delta_{c,k_c}} \big) ,
\end{split}
\end{equation*}

For $\cI_{F,5}$, consider first
\begin{equation}
\label{eqn: E1_rate}
\begin{split}
    &\hspace{5mm}
    \b{E}\Big\{ \|\E_t\A_c\|_F^2 \Big\}
    \leq p k_c \max_{i\in[p], j\in[k_c]} \b{E}\Big\{ (\E_{t,i\cdot}' \A_{c, \cdot j})^2 \Big\} \\
    &=
    p k_c \max_{i\in[p], j\in[k_c]} \sum_{n=1}^q \sum_{l=1}^q \text{cov}(E_{t, in}, E_{t, il}) \cdot  A_{c,nj} A_{c,lj} = O(pq) ,
\end{split}
\end{equation}
where the last equality used Assumptions (E1) and (E2). Thus,
\begin{equation*}
\begin{split}
    \|\cI_{F,5} \|_F^2 &= O_P(pq) \cdot \|\Q_r\H_r' - \wh\Q_r\|_F^2 \cdot \|\Z_c^{-1/2}\|_F^2 \\
    &= O_P\big(T^{-1} p^{2+ \delta_{r,1} - 2\delta_{r,k_r}} q^{2-2\delta_{c,1}-\delta_{c,k_c}} + p^{2-3\delta_{r,k_r}} q^{3-2\delta_{c,1}-\delta_{c,k_c}} \big) ,
\end{split}
\end{equation*}
where we used Lemma \ref{lemma: corollary_asymp_loading} and Assumption (L1). Similarly for $\cI_{F,6}$,
\begin{equation*}
    \|\cI_{F,6} \|_F^2 = O_P\big(T^{-1} q^{2+ \delta_{c,1} - 2\delta_{c,k_c}} p^{2-2\delta_{r,1}-\delta_{r,k_r}} + q^{2-3\delta_{c,k_c}} p^{3-2\delta_{r,1}-\delta_{r,k_r}} \big) .
\end{equation*}

By Assumptions (E1) and (E2) again, we have
\begin{equation*}
\begin{split}
    &\hspace{5mm}
    \b{E}\Big\{ \|\A_r'\E_t\A_c\|_F^2 \Big\}
    \leq k_r k_c \max_{i\in[k_r], j\in[k_c]} \b{E}\Big\{ (\A_{r,\cdot i}' \E_t \A_{c, \cdot j})^2 \Big\} \\
    &=
    k_r k_c \max_{i\in[k_r], j\in[k_c]} \sum_{m=1}^p \sum_{n=1}^q \sum_{h=1}^p \sum_{l=1}^q \text{cov}(E_{t, mn}, E_{t, hl}) \cdot A_{r, mi} A_{c,nj} A_{r, hi} A_{c,lj} = O(pq) ,
\end{split}
\end{equation*}
hence for $\cI_{F,7}$, it holds that
\begin{equation*}
    \|\cI_{F,7} \|_F^2 = O_P(pq) \cdot \|\Z_r^{-1/2}\|_F^2 \cdot \|\Z_c^{-1/2}\|_F^2 = O_P\big( p^{1-\delta_{r,k_r}} q^{1-\delta_{c,k_c}} \big) .
\end{equation*}

Consider $\cI_{F,8}$, recall that $(\1_p' \E_t \1_q)^2 = O_P(pq)$ from (\ref{eqn: 1E1_bound}) and hence
\begin{equation*}
    \|\cI_{F,8} \|_F^2 = O_P(p^{-2} q^{-2}) \cdot (\1_p' \E_t \1_q)^2 \cdot \|\1_p \1_q'\|_F^2 = O_P(1) .
\end{equation*}

For $\cI_{F,9}$, note that $\b{E}[\|\E_t\1_q\|_F^2]$ has the same rate as (\ref{eqn: E1_rate}) since $k_c$ is fixed, and hence
\begin{equation*}
\begin{split}
    \|\cI_{F,9} \|_F^2 &= O_P(q^{-2}) \cdot \|\Q_r\H_r' - \wh\Q_r\|_F^2 \cdot \|\E_t\1_q\|_F^2 \cdot \|\1_q\|^2 \\
    &= O_P\big( T^{-1} p^{2+ \delta_{r,1} - 2\delta_{r,k_r}} q^{1-2\delta_{c,1}} + p^{2-3\delta_{r,k_r}} q^{2-2\delta_{c,1}} \big).
\end{split}
\end{equation*}
From (\ref{eqn: 1E1_bound}), we also have $\b{E}[\|\A_r'\E_t\1_q\|_F^2] = O(pq)$, so for $\cI_{F,10}$ we have
\begin{equation*}
    \|\cI_{F,10} \|_F^2 = O_P(q^{-2}) \cdot \|\Z_r^{-1/2}\|_F^2 \cdot \|\A_r'\E_t\1_q\|_F^2 \cdot \|\1_q\|^2 = O_P\big( p^{1-\delta_{r,k_r}} \big).
\end{equation*}
Lastly, the rates for $\cI_{F,11}$ and $\cI_{F,12}$ can be obtained similarly as $\cI_{F,9}$ and $\cI_{F,10}$,
\begin{align*}
    & \|\cI_{F,11} \|_F^2 = O_P\big( T^{-1} q^{2+ \delta_{c,1} - 2\delta_{c,k_c}} p^{1-2\delta_{r,1}} + q^{2-3\delta_{c,k_c}} p^{2-2\delta_{r,1}} \big), \\
    & \|\cI_{F,12} \|_F^2 = O_P\big( q^{1-\delta_{c,k_c}} \big).
\end{align*}
Therefore, with all the rates from $\cI_{F,1}$ to $\cI_{F,12}$ in (\ref{eqn: F_Zt_decomp}), by using Assumption (R1) we have
\begin{equation*}
\begin{split}
    \|\wh\F_{Z,t} - (\H_r^{-1})' \F_{Z,t} \H_c^{-1}\|_F^2 &=
    O_P\big( p^{1-\delta_{r,k_r}} q^{1-\delta_{c,k_c}}
    + T^{-1} p^{1+ 2\delta_{r,1} - 2\delta_{r,k_r}} q^{1-\delta_{c,1}} +  p^{1 +\delta_{r,1} -3\delta_{r,k_r}} q^{2-\delta_{c,1}} \\
    & \hspace{36pt}
    + T^{-1} q^{1+ 2\delta_{c,1} - 2\delta_{c,k_c}} p^{1-\delta_{r,1}} +  q^{1 +\delta_{c,1} -3\delta_{c,k_c}} p^{2-\delta_{r,1}}\big) .
\end{split}
\end{equation*}
This shows the first statement of Theorem \ref{thm: consistency_F_C}.

For the remaining proof, consider $\wh{C}_{t,ij} - C_{t,ij}$ for any $t\in[T], i\in[p], j\in[q]$. First, we have
\begin{equation}
\label{eqn: hat_Ct_decomp}
\begin{split}
    &\hspace{5mm}
    \wh{C}_{t,ij} - C_{t,ij} = \wh\Q_{r,i\cdot}' \wh\F_{Z,t} \wh\Q_{c,j\cdot} - \Q_{r,i\cdot}' \F_{Z,t} \Q_{c,j\cdot} \\
    &= (\wh\Q_{r,i\cdot}-\H_r\Q_{r,i\cdot})' (\wh\F_{Z,t} - (\H_r^{-1})' \F_{Z,t} \H_c^{-1}) (\wh\Q_{c,j\cdot} - \H_c\Q_{c,j\cdot}) \\
    & + (\wh\Q_{r,i\cdot}-\H_r\Q_{r,i\cdot})' (\wh\F_{Z,t} - (\H_r^{-1})' \F_{Z,t} \H_c^{-1}) \H_c\Q_{c,j\cdot}
    + (\wh\Q_{r,i\cdot}-\H_r\Q_{r,i\cdot})' (\H_r^{-1})' \F_{Z,t} \H_c^{-1} (\wh\Q_{c,j\cdot}- \H_c\Q_{c,j\cdot}) \\
    & + (\wh\Q_{r,i\cdot}-\H_r\Q_{r,i\cdot})' (\H_r^{-1})' \F_{Z,t}\Q_{c,j\cdot}
    + \Q_{r,i\cdot}'\H_r' (\wh\F_{Z,t} - (\H_r^{-1})' \F_{Z,t} \H_c^{-1}) (\wh\Q_{c,j\cdot} -\H_c \Q_{c,j\cdot}) \\
    & + \Q_{r,i\cdot}'\H_r' (\wh\F_{Z,t} - (\H_r^{-1})' \F_{Z,t} \H_c^{-1}) \H_c \Q_{c,j\cdot}
    + \Q_{r,i\cdot}' \F_{Z,t} \H_c^{-1} (\wh\Q_{c,j\cdot} - \H_c \Q_{c,j\cdot}).
\end{split}
\end{equation}
Notice that from Assumption (L1),
\begin{align*}
    & \|\Q_{r,i\cdot}\|^2 = O(\|\A_{r,i\cdot}\|^2) \cdot \|\Z_{r}^{-1/2}\|_F^2 = O(p^{-\delta_{r,k_r}}) , \\
    & \|\Q_{c,j\cdot}\|^2 = O(\|\A_{c,j\cdot}\|^2) \cdot \|\Z_{c}^{-1/2}\|_F^2 = O(q^{-\delta_{c,k_c}}) .
\end{align*}
Together with Lemma \ref{lemma: corollary_asymp_loading}, the first statement of Theorem \ref{thm: consistency_F_C} and Assumption (R1), we have
\begin{equation*}
\begin{split}
    (\wh{C}_{t,ij} - C_{t,ij})^2
    &= O_P\big(p^{1-2\delta_{r,k_r}} q^{1-2\delta_{c,k_c}}
    + T^{-1} p^{1 + 2\delta_{r,1} -3\delta_{r,k_r}} q^{1 -\delta_{c,1} -\delta_{c,k_c}} + p^{1 +\delta_{r,1} -4\delta_{r,k_r}} q^{2 -\delta_{c,1} -\delta_{c,k_c}} \\
    & \hspace{36pt}
    + T^{-1} q^{1 + 2\delta_{c,1} -3\delta_{c,k_c}} p^{1 -\delta_{r,1} -\delta_{r,k_r}} + q^{1 +\delta_{c,1} -4\delta_{c,k_c}} p^{2 -\delta_{r,1} -\delta_{r,k_r}}
    \big) .
\end{split}
\end{equation*}
This completes the proof of Theorem \ref{thm: consistency_F_C}. $\square$

\textbf{\textit{Proof of Theorem \ref{thm:alpha_normality}.}}
We first consider $\wh\mu_t$, which is given by
\begin{align*}
  \wh\mu_t &= \1_p'\Y_t\1_q/(pq) = \mu_t + \1_p'\E_t\1_q/(pq) = \mu_t + \frac{1}{pq}\sum_{i,j}E_{t,ij}, \text{ so that }\\
  \wh\mu_t - \mu_t &= \frac{1}{pq}\1_p'\A_{e,r}\F_{e,t}\A_{e,c}'\1_q + \frac{1}{pq}\sum_{i,j}\Sigma_{\epsilon,ij}\epsilon_{t,ij} =: I_{\mu,1} + I_{\mu,2}.
\end{align*}
By Assumption (E1), since $\|\A_{e,r}\|_1, \|\A_{e,c}\|_1 = O(1)$, we have $I_{\mu,1} = O_P(1/(pq))$. Also, $I_{\mu,2} = O((pq)^{-1/2})$ since $\b{E}(I_{\mu,2}) = 0$ and $\var(I_{\mu,2}) = (pq)^{-2}\sum_{i,j}\Sigma_{\epsilon,ij}^2 = O(1/(pq))$. Hence $I_{\mu,1}$ is dominated by $I_{\mu,2}$, and
\begin{align*}
  \sqrt{pq}(\wh\mu_t - \mu_t) &=  \sqrt{pq}I_{\mu,2}(1 + o_P(1))
  = \frac{1}{\sqrt{pq}}\sum_{i,j}\Sigma_{\epsilon,ij}\epsilon_{t,ij}(1+o_P(1)) \xrightarrow{\c{D}} \c{N}(0,\gamma_{\mu}^2),
\end{align*}
where we use Theorem 1 in \cite{AyvazyanUlyanov2023} for the convergence in distribution.

For $\wh\balpha_t$, consider the decomposition
\begin{align*}
  \wh\balpha_t - \balpha_t &= (\mu_t - \wh\mu_t)\1_p + q^{-1}\E_t\1_q\\
  &= (\mu_t - \wh\mu_t)\1_p + q^{-1}\A_{e,r}\F_{e,t}\A_{e,c}'\1_q + q^{-1}(\bSigma_{\epsilon}\circ\bepsilon_t)\1_q, \text{ so that }\\
  \wh\alpha_{t,i} - \alpha_{t,i} &= (\mu_t - \wh\mu_t) + q^{-1}(\A_{e,r})_{i\cdot}\F_{e,t}\A_{e,c}'\1_q
  + q^{-1}\sum_{j=1}^q\Sigma_{\epsilon, ij}\epsilon_{t,ij}.
\end{align*}
From what we have proved above we have $|\wh\mu_t - \mu_t| = O_P(1/\sqrt{pq})$. Also, $q^{-1}(\A_{e,r})_{i\cdot}\F_{e,t}\A_{e,c}'\1_q = O_P(1/q)$ since $\|\A_{e,r}\|_1, \|\A_{e,c}\|_1 = O(1)$. Finally, $\b{E}(q^{-1}(\bSigma_{\epsilon}\circ \bepsilon_t)\1_q) = \0$ and $\var(q^{-1}\sum_{j=1}^q\Sigma_{\epsilon,ij}\epsilon_{t,ij}) = q^{-2}\sum_{j=1}^q\Sigma_{\epsilon,ij}^2 = O(q^{-1})$, implying that $q^{-1}(\bSigma_{\epsilon}\circ\bepsilon_t)\1_q = O_P(q^{-1/2})$ elementwise, thus dominating other terms. Hence for $i\in[p]$, using Theorem 1 in \cite{AyvazyanUlyanov2023},
\begin{align*}
  \sqrt{q}(\wh\alpha_{t,i} - \alpha_{t,i}) &= q^{-1/2}\sum_{j=1}^q\Sigma_{\epsilon,ij}\epsilon_{t,ij}(1+o_P(1)) \xrightarrow{\c{D}} \c{N}(0,\gamma_{\alpha,i}^2).
\end{align*}
By Assumption (E1), each element of $q^{-1}(\bSigma_{\epsilon}\circ\bepsilon_t)\1_q$ is independent of each other. Hence with integers $i_1<\cdots<i_m$, $m$ being finite and $\btheta_{\alpha,t} := (\alpha_{t,i_1},\ldots, \alpha_{t,i_m})'$, by Theorem 1 in \cite{AyvazyanUlyanov2023},
\begin{align*}
  \sqrt{q}(\wh\btheta_{\alpha,t} - \btheta_{\alpha,t}) \xrightarrow{\c{D}} \c{N}(\0, \diag(\gamma_{\alpha,i_1}^2,\ldots, \gamma_{\alpha,i_m}^2)).
\end{align*}
We omit the proof of asymptotic normality for $\wh\bbeta_t$ since the arguments used are in parallel to those used for $\wh\balpha_t$, using the independence of the columns in $\bSigma_{\epsilon}\circ\bepsilon_t$ by Assumption (E1).

The rest of the proof is done if we can prove that $\wh\gamma_{\alpha,i}$, $\wh\gamma_{\beta,j}$ and $\wh\gamma_{\mu}$ are consistent estimators for $\gamma_{\alpha,i}$, $\gamma_{\beta,j}$ and $\gamma_{\mu}$ respectively. From (\ref{eqn: Ehat_t}), since we assume $\wh{C}_{t,ij} - C_{t,ij} = o_P(1)$ from Theorem \ref{thm: consistency_F_C}, then elementwise we have
\begin{align*}
  \wh\E_t &= \wh\L_t - \wh\C_t = (\mu_t-\wh\mu_t)\1_p\1_q' + (\balpha_t-\wh\balpha_t)\1_q' + \1_p(\bbeta_t-\wh\bbeta_t)' + (\C_t - \wh\C_t) + \E_t = \E_t(1+o_P(1)).
\end{align*}
Hence we have
\begin{align*}
  q^{-1}(\wh\E_t\wh\E_t')_{ii} &= \{q^{-1}(\A_{e,r})_{i\cdot}\F_{e,t}\A_{e,c}'\A_{e,c}\F_{e,t}'(\A_{e,r})_{i\cdot}' + q^{-1}(\A_{e,r})_{i\cdot}\F_{e,t}\A_{e,c}'(\bSigma_{\epsilon}\circ\bepsilon_t)_{i\cdot}' \\
  &\quad + q^{-1}(\bSigma_{\epsilon}\circ \bepsilon_t)_{i\cdot}(\bSigma_{\epsilon}\circ\bepsilon_t)_{i\cdot}'\}(1+o_P(1))\\
  &= O_P(q^{-1}) + O_P(q^{-1}) + \frac{1}{q}\sum_{j=1}^q\Sigma_{\epsilon,ij}^2(\epsilon_{t,ij}^2-1)(1+o_P(1)) + \frac{1}{q}\sum_{j=1}^q\Sigma_{\epsilon,ij}^2\\
  &\xrightarrow{\c{P}} \gamma_{\alpha,i}^2,
\end{align*}
where we used the Markov inequality to arrive at $q^{-1}\sum_{j=1}^q\Sigma_{\epsilon,ij}^2(\epsilon_{t,ij}^2-1) = O_P(q^{-1/2})$, knowing that each $\Sigma_{\epsilon,ij}$ is bounded away from infinity by Assumption (E1). A parallel argument (omitted) can show that $\wh\gamma_{\beta,j}^2$ is consistent for $\gamma_{\beta,j}^2$. Finally,
\begin{align*}
  \wh\gamma_{\mu}^2 = p^{-1}\sum_{i=1}^p q^{-1}(\wh\E_{t}\wh\E_t')_{ii} \xrightarrow{\c{P}} p^{-1}\sum_{i=1}^p \gamma_{\alpha,i}^2 = \gamma_{\mu}^2.
\end{align*}
This completes the proof of the theorem. $\square$

\textbf{\textit{Proof of Theorem \ref{thm: asymp_loading}.}}
We construct the asymptotic normality for rows of our factor loading estimators. We only prove the result for the row loading estimator, and the proof for the column loading estimator would be similar (and hence omitted). For any $j\in[p]$, consider the following decomposition
\begin{equation}
\label{eqn: Q_asymp_decomp}
\begin{split}
    \wh\Q_{r,j\cdot} - \H_r\Q_{r,j\cdot} &=
    T^{-1}\wh\D_r^{-1} \sum_{i=1}^p
    \wh\Q_{r,i\cdot} \sum_{t=1}^T (\R_{r,t})_{ij}
    \\
    &=
    T^{-1}\wh\D_r^{-1} \sum_{i=1}^p
    (\wh\Q_{r,i\cdot} - \H_r\Q_{r,i\cdot}) \sum_{t=1}^T (\R_{r,t})_{ij}
    + T^{-1}\wh\D_r^{-1} \sum_{i=1}^p
    \H_r\Q_{r,i\cdot} \sum_{t=1}^T (\R_{r,t})_{ij}  \\
    &=
    \cI_1 + \cI_2 + \cI_3 + \cI_4 + \cI_5 - \cI_6 - \cI_7 - \cI_8 - \cI_9 - \cI_{10} - \cI_{11} + \cI_{12} + \cI_{13},
    \;\;\; \text{where}
\end{split}
\end{equation}
\vspace{-24pt}
\begin{align*}
    \cI_1 &:=
    \frac{1}{T}\wh\D_r^{-1} \sum_{i=1}^p
    \H_r \Q_{r,i\cdot} \sum_{t=1}^T (\C_t \E_t')_{ij} , \;\;\;
    \cI_2 :=
    \frac{1}{T}\wh\D_r^{-1} \sum_{i=1}^p
    \H_r \Q_{r,i\cdot} \sum_{t=1}^T (\E_t\C_t')_{ij} , \\
    \cI_3 &:=
    \frac{1}{T}\wh\D_r^{-1} \sum_{i=1}^p
    \H_r \Q_{r,i\cdot} \sum_{t=1}^T (\E_t\E_t')_{ij} , \;\;\;
    \cI_4 :=
    \frac{1}{Tpq}\wh\D_r^{-1} \sum_{i=1}^p
    \H_r \Q_{r,i\cdot} \sum_{t=1}^T (\1_q'\E_t'\1_p\E_t\1_q\1_p')_{ij} , \\
    \cI_5 &:=
    \frac{1}{Tpq}\wh\D_r^{-1} \sum_{i=1}^p
    \H_r \Q_{r,i\cdot} \sum_{t=1}^T (\1_q'\E_t'\1_p\1_p\1_q'\E_t')_{ij} , \;\;\;
    \cI_6 :=
    \frac{1}{Tp}\wh\D_r^{-1} \sum_{i=1}^p
    \H_r \Q_{r,i\cdot} \sum_{t=1}^T (\C_t\E_t'\1_p\1_p')_{ij} , \\
    \cI_7 &:=
    \frac{1}{Tp}\wh\D_r^{-1} \sum_{i=1}^p
    \H_r \Q_{r,i\cdot} \sum_{t=1}^T (\1_p\1_p'\E_t\C_t')_{ij} , \;\;\;
    \cI_8 :=
    \frac{1}{Tp}\wh\D_r^{-1} \sum_{i=1}^p
    \H_r \Q_{r,i\cdot} \sum_{t=1}^T (\E_t\E_t'\1_p\1_p')_{ij} , \\
    \cI_9 &:=
    \frac{1}{Tp}\wh\D_r^{-1} \sum_{i=1}^p
    \H_r \Q_{r,i\cdot} \sum_{t=1}^T (\1_p\1_p'\E_t\E_t')_{ij} , \;\;\;
    \cI_{10} :=
    \frac{1}{Tp}\wh\D_r^{-1} \sum_{i=1}^p
    \H_r \Q_{r,i\cdot} \sum_{t=1}^T (\E_t\1_q\1_q'\E_t')_{ij} , \\
    \cI_{11} &:=
    \frac{1}{T p^2 q}\wh\D_r^{-1} \sum_{i=1}^p
    \H_r \Q_{r,i\cdot} \sum_{t=1}^T ((\1_q'\E_t'\1_p)^2\1_p\1_p')_{ij} , \;\;\;
    \cI_{12} :=
    \frac{1}{T p^2}\wh\D_r^{-1} \sum_{i=1}^p
    \H_r \Q_{r,i\cdot} \sum_{t=1}^T (\1_p'\E_t\E_t'\1_p\1_p\1_p')_{ij} , \\
    \cI_{13} &:=
    \frac{1}{T}\wh\D_r^{-1} \sum_{i=1}^p
    (\wh\Q_{r,i\cdot} - \H_r\Q_{r,i\cdot}) \sum_{t=1}^T (\R_{r,t})_{ij} .
\end{align*}
We shall show that $\cI_1$ is the leading term among the decomposition in (\ref{eqn: Q_asymp_decomp}). To obtain the rate for $\cI_2$, from Assumptions (E1) and (E2) we have for any $i\in[p], h\in[q]$,
\begin{equation*}
    E_{t,ih} = \sum_{w\geq 0}a_{e,w}
    \A_{e,r,i\cdot}'\X_{e,t-w}\A_{e,c,h\cdot}+
    (\bSigma_{\epsilon})_{ih}\sum_{w\geq 0}a_{\epsilon,w}
    (\X_{\epsilon,t-w})_{ih}.
\end{equation*}
Consider first $\sum_{t=1}^T \sum_{h=1}^q(\sum_{w\geq 0}a_{e,w} \A_{e,k,i\cdot}'\X_{e,t-w}\A_{e,\text{-}k,h\cdot}) \A_{c,h\cdot}'\F_t' \A_{r,j\cdot}$. By Assumptions (F1), (E1) and (E2), we have
\begin{equation}
\begin{split}
\label{eqn: I2_step1}
    & \hspace{5mm}
    \b{E}\Big\{\Big(
    \sum_{t=1}^T \sum_{h=1}^q(\sum_{w\geq 0}a_{e,w}\A_{e,k,i\cdot}'\X_{e,t-w}\A_{e,\text{-}k,h\cdot}) \A_{c,h\cdot}'\F_t' \A_{r,j\cdot}
    \Big)^2\Big\} \\
    &=
    \text{cov}\Big(
    \sum_{t=1}^T \sum_{h=1}^q
    \A_{c,h\cdot}'
    (\sum_{w\geq 0}a_{f,w}\X_{f,t-w}' ) \A_{r,j\cdot}
    (\sum_{w\geq 0}a_{e,w}\A_{e,k,i\cdot}'\X_{e,t-w}\A_{e,\text{-}k,h\cdot})
    , \\
    & \hspace{36pt}
    \sum_{t=1}^T \sum_{h=1}^q
    \A_{c,h\cdot}'
    (\sum_{w\geq 0}a_{f,w}\X_{f,t-w}' ) \A_{r,j\cdot}
    (\sum_{w\geq 0}a_{e,w}\A_{e,k,i\cdot}'\X_{e,t-w}\A_{e,\text{-}k,h\cdot}) \Big) \\
    & =
    \sum_{h=1}^q \sum_{l=1}^q \sum_{t=1}^T
    \sum_{w\geq 0} a_{f,w}^2 a_{e,w}^2 \cdot
    \|\A_{r,j\cdot}\|^2 \cdot
    \| \A_{c,h\cdot}\|\cdot
    \| \A_{c,l\cdot}\|\cdot
    \|\A_{e,c,h\cdot}\|\cdot
    \|\A_{e,c,l\cdot}\|\cdot
    \|\A_{e,r,i\cdot}\|^2 \\
    &
    = O(T)\cdot \|\A_{r,j\cdot}\|^2 \cdot
    \|\A_{e,r,i\cdot}\|^2 .
\end{split}
\end{equation}
Consider also $\sum_{i=1}^p \sum_{h=1}^q \sum_{t=1}^T \Q_{r,i\cdot} ((\bSigma_\epsilon)_{ih} \sum_{w\geq 0}a_{\epsilon,w} (\X_{\epsilon,t-w})_{ih}) \A_{c,h\cdot}'\F_t' \A_{r,j\cdot}$. Similarly, by Assumptions (E1), (E2) and (F1), we have
\begin{equation}
\begin{split}
\label{eqn: I2_step2}
    & \hspace{5mm}
    \b{E}\Big\{ \Big\|
    \sum_{i=1}^p \sum_{h=1}^q \sum_{t=1}^T \Q_{r,i\cdot} ((\bSigma_\epsilon)_{ih} \sum_{w\geq 0}a_{\epsilon,w} (\X_{\epsilon,t-w})_{ih}) \A_{c,h\cdot}'\F_t' \A_{r,j\cdot}
    \Big\|^2 \Big\} \\
    &=
    \text{cov}\Big(
    \sum_{i=1}^p \sum_{h=1}^q \sum_{t=1}^T \Q_{r,i\cdot} ((\bSigma_\epsilon)_{ih} \sum_{w\geq 0}a_{\epsilon,w} (\X_{\epsilon,t-w})_{ih}) \A_{c,h\cdot}'\F_t' \A_{r,j\cdot}
    , \\
    & \hspace{36pt}
    \sum_{i=1}^p \sum_{h=1}^q \sum_{t=1}^T \Q_{r,i\cdot} ((\bSigma_\epsilon)_{ih} \sum_{w\geq 0}a_{\epsilon,w} (\X_{\epsilon,t-w})_{ih}) \A_{c,h\cdot}'\F_t' \A_{r,j\cdot} \Big) \\
    &=
    \sum_{i=1}^p \sum_{h=1}^q \sum_{t=1}^T \sum_{w\geq 0} a_{f,w}^2 a_{\epsilon,w}^2
    \cdot \|\A_{r,j\cdot}\|^2 \cdot \|\A_{c,h\cdot}\|^2 \cdot (\bSigma_\epsilon)_{ih}^2
    \cdot \|\Q_{r,i\cdot}\|^2 \\
    &=
    O(T) \cdot \|\A_{r,j\cdot}\|^2 \cdot \|\A_c\|^2 \cdot \|\Q_r\|^2 .
\end{split}
\end{equation}
Hence using Lemma \ref{lemma: norm_hat_D}, it holds that
\begin{equation*}
\begin{split}
    \|\cI_2\|^2 &\leq \frac{1}{T^2} \|\wh\D_r^{-1}\|_F^2 \cdot \|\H_r\|_F^2 \cdot \Big\| \sum_{i=1}^p
    \Q_{r,i\cdot} \sum_{t=1}^T \sum_{h=1}^q E_{t,ih} \A_{c,h\cdot}' \F_t' \A_{r,j\cdot}  \Big\|^2 \\
    &=
    O_P\Big(\frac{1}{T^2} p^{-2\delta_{r,k_r}} q^{-2\delta_{c,1}}\Big)\Big\{
    \Big\| \sum_{i=1}^p \sum_{h=1}^q \sum_{t=1}^T \Q_{r,i\cdot} ((\bSigma_\epsilon)_{ih} \sum_{w\geq 0}a_{\epsilon,w} (\X_{\epsilon,t-w})_{ih}) \A_{c,h\cdot}'\F_t' \A_{j\cdot}
    \Big\|^2 \\
    &+
    \Big(\sum_{i=1}^p \|\Q_{r,i\cdot}\|^2 \Big)
    \sum_{i=1}^p \Big(
    \sum_{t=1}^T \sum_{h=1}^q(\sum_{w\geq 0}a_{e,w}\A_{e,k,i\cdot}'\X_{e,t-w}\A_{e,\text{-}k,h\cdot}) \A_{c,h\cdot}'\F_t' \A_{r,j\cdot} \Big)^2 \Big\} \\
    &=
    O_P\Big(\frac{1}{T^2} p^{-2\delta_{r,k_r}} q^{-2\delta_{c,1}} \cdot Tq^{\delta_{c,1}}\Big)
    = O_P\big(T^{-1} p^{-2\delta_{r,k_r}} q^{-\delta_{c,1}} \big) ,
\end{split}
\end{equation*}
where we used Assumption (L1), (\ref{eqn: I2_step1}) and (\ref{eqn: I2_step2}) in the last equality.

For $\cI_3$, first notice from the noise structure in Assumptions (E1) and (E2), we have
\begin{equation*}
\begin{split}
    &\hspace{5mm}
    \text{Var}\Big(
    \sum_{i=1}^p \sum_{h=1}^q \sum_{t=1}^T
    \Q_{r,i\cdot} E_{t,ih} E_{t,jh}\Big) \\
    & =
    O(1) \cdot
    \sum_{i=1}^p \sum_{u=1}^p \sum_{h=1}^q \sum_{l=1}^q \sum_{t=1}^T \sum_{n=1}^{k_{e,r}} \sum_{m=1}^{k_{e,c}} \sum_{w\geq 0}
    a_{e,w}^4 A_{e,r,in} A_{e,r,un} A_{e,r,jn}^2 A_{e,c,hm}^2 A_{e,c,lm}^2 \\
    &\cdot
    \|\Q_{r,i\cdot}\|
    \cdot \|\Q_{r,u\cdot}\|
    \cdot \text{Var}((\X_{e,t-w})_{nm}^2)\\
    & +
    O(1) \cdot \sum_{i=1}^p \sum_{h=1}^q
    \sum_{t=1}^T \sum_{w\geq 0} a_{\epsilon,w}^4
    (\bSigma_\epsilon)_{ih}^2
    (\bSigma_\epsilon)_{jh}^2
    \cdot \|\Q_{r,i\cdot}\|^2
    \cdot \text{Var}(
    (\X_{\epsilon,t-w})_{ih}
    (\X_{\epsilon,t-w})_{jh}) \\
    &
    = O(T + T q)=O(T q).
\end{split}
\end{equation*}
Moreover, it holds that
\begin{equation*}
    \b{E}\Big( \sum_{i=1}^p \Big| \sum_{h=1}^q \sum_{t=1}^T E_{t,ih} E_{t,jh} \Big|\Big) =
    \sum_{i=1}^p \Big| \sum_{h=1}^q \sum_{t=1}^T \big( \|\A_{e,c,h\cdot}\|^2 \cdot \|\A_{e,r,i\cdot}\| \cdot \|\A_{e,r,j\cdot}\| +
    (\bSigma_\epsilon)_{ih}\b{1}_{\{i=j\}} \big) \Big| =O(Tq),
\end{equation*}
and with $\max_i\|\Q_{r,i\cdot}\|^2 \leq\|\A_{r,j\cdot}\|^2 \cdot \|\Z_r^{-1/2}\|^2 = O_P\big(p^{-\delta_{r,k_r}}\big)$, we thus have
\begin{equation}
\label{eqn: I3_last_step}
\begin{split}
    &\hspace{5mm}
    \|\cI_3\|^2 \leq
    \frac{1}{T^2} \|\wh\D_r^{-1}\|_F^2 \cdot \|\H_r\|_F^2 \cdot \Big\| \sum_{i=1}^p
     \Q_{r,i\cdot} \sum_{t=1}^T \sum_{h=1}^q E_{t,ih} E_{t,jh} \Big\|^2 \\
     &=
     O_P\Big(T^{-2} p^{-2\delta_{r,k_r}} q^{-2\delta_{c,1}} (Tq+ T^2 q^2 p^{-\delta_{r,k_r}}) \Big)
     =
     O_P\big(T^{-1} p^{-2\delta_{r,k_r}} q^{1-2\delta_{c,1}} + p^{-3\delta_{r,k_r}} q^{2-2\delta_{c,1}} \big) .
\end{split}
\end{equation}

Consider now $\cI_4$ and $\cI_5$. From the proof of (\ref{eqn: E11_bound}) in Lemma \ref{lemma: rate_R_rt},
\begin{equation*}
\begin{split}
    &\hspace{5mm}
    \|\cI_4\|^2 \leq
    \frac{1}{T^2 p^2 q^2} \|\wh\D_r^{-1}\|_F^2 \cdot \|\H_r\|_F^2 \cdot \Big\| \sum_{i=1}^p
    \Q_{r,i\cdot} \sum_{t=1}^T (\1_q'\E_t'\1_p\E_t\1_q\1_p')_{ij} \Big\|^2 \\
    &\leq
    \frac{1}{T^2 p^2 q^2} \|\wh\D_r^{-1}\|_F^2 \cdot \|\H_r\|_F^2 \cdot \Big( \sum_{i=1}^p
    \|\Q_{r,i\cdot}\|^2 \Big) \sum_{i=1}^p \Big( \sum_{t=1}^T \1_q'\E_t'\1_p\E_t\1_q\1_p' \Big)_{ij}^2 \\
    &=
    O_P\big(T^{-2} p^{-2-2\delta_{r,k_r}} q^{-2-2\delta_{c,1}} \big) \cdot \Big\| \Big( \sum_{t=1}^T \1_q'\E_t'\1_p\E_t\1_q\1_p' \Big)_{\cdot j} \Big\|^2 \\
    &=
    O_P\big(T^{-2} p^{-2-2\delta_{r,k_r}} q^{-2-2\delta_{c,1}} \big) \cdot
    \sum_{k=1}^p \b{E} \Bigg\{ \Bigg[\sum_{t=1}^T \Big( \sum_{i=1}^p \sum_{u=1}^q E_{t,iu} \Big) \sum_{h=1}^q E_{t,kh} \Bigg]^2\Bigg\} \\
    &=
    O_P\Big(T^{-2} p^{-2-2\delta_{r,k_r}} q^{-2-2\delta_{c,1}} (Tp^2q^2 + T^2pq^2) \Big) \\
    &=
    O_P\big(T^{-1} p^{-2\delta_{r,k_r}} q^{-2\delta_{c,1}} + p^{-1-2\delta_{r,k_r}} q^{-2\delta_{c,1}}\big) ,
\end{split}
\end{equation*}
where we used Lemma \ref{lemma: norm_hat_D} in the third line and (\ref{eqn: E11_bound_last_step}) in the second last equality. Similarly, we have
\begin{equation*}
\begin{split}
    &\hspace{5mm}
    \|\cI_5\|^2 \leq
    \frac{1}{T^2 p^2 q^2} \|\wh\D_r^{-1}\|_F^2 \cdot \|\H_r\|_F^2 \cdot \Big\| \sum_{i=1}^p
    \Q_{r,i\cdot} \sum_{t=1}^T (\1_q'\E_t'\1_p\1_p\1_q'\E_t')_{ij} \Big\|^2 \\
    &=
    O_P\big(T^{-2} p^{-2-2\delta_{r,k_r}} q^{-2-2\delta_{c,1}} \big) \cdot \Big\| \Big( \sum_{t=1}^T \1_q'\E_t'\1_p\E_t\1_q\1_p' \Big)_{j\cdot} \Big\|^2 \\
    &=
    O_P\big(T^{-2} p^{-2-2\delta_{r,k_r}} q^{-2-2\delta_{c,1}} \big) \cdot
    p\cdot \b{E} \Bigg\{ \Bigg[\sum_{t=1}^T \Big( \sum_{i=1}^p \sum_{u=1}^q E_{t,iu} \Big) \sum_{h=1}^q E_{t,jh} \Bigg]^2\Bigg\} \\
    &=
    O_P\Big(T^{-2} p^{-2-2\delta_{r,k_r}} q^{-2-2\delta_{c,1}} (Tp^2q^2 + T^2pq^2) \Big) \\
    &=
    O_P\big(T^{-1} p^{-2\delta_{r,k_r}} q^{-2\delta_{c,1}} + p^{-1-2\delta_{r,k_r}} q^{-2\delta_{c,1}}\big) ,
\end{split}
\end{equation*}
where we used again Lemma \ref{lemma: norm_hat_D} in the third line and (\ref{eqn: E11_bound_last_step}) in the second last equality.

Consider now $\cI_6$, note first from Assumptions (E1) and (E2),
\begin{equation*}
    E_{t,hu} = \sum_{w\geq 0}a_{e,w}
    \A_{e,r,h\cdot}'\X_{e,t-w}\A_{e,c,u\cdot} +
    (\bSigma_{\epsilon})_{hu} \sum_{w\geq 0} a_{\epsilon,w}(\X_{\epsilon,t-w})_{hu}.
\end{equation*}
Thus, we have
\begin{equation*}
\begin{split}
    &\hspace{5mm}
    \text{cov}(E_{t,hu}, E_{s,vl}) =
    \b{E}[\A_{e,r,h\cdot}' \F_{e,t} \A_{e,c,u\cdot} \A_{e,c,l\cdot}' \F_{e,s}'\A_{e,r,v\cdot}]
    + \b{1}_{\{h=v\}} \b{1}_{\{u=l\}} \cdot (\bSigma_{\epsilon})_{hu}^2
    \cdot \b{E}[\bepsilon_{t,\cdot u}'
    \bepsilon_{s,\cdot l}] \\
    &=
    \A_{e,c,l\cdot}' \A_{e,c,u\cdot}
    \A_{e,r,h\cdot}' \A_{e,r,v\cdot}
    \cdot \sum_{w\geq 0}a_{e,w}a_{e,w+|t-s|} +
    \b{1}_{\{h=v\}} \b{1}_{\{u=l\}} \cdot (\bSigma_{\epsilon})_{hu}^2 \sum_{w\geq 0} a_{\epsilon,w} a_{\epsilon,w+|t-s|} .
\end{split}
\end{equation*}
Hence if we fix $t\in[T], h\in[p], u\in[q]$, then for any deterministic vectors $\bf{n}\in \b{R}^{k_r}$ and $\bf{g}, \bf{j}\in \b{R}^{k_c}$, we have
\begin{equation*}
\begin{split}
    &\hspace{5mm}
    \sum_{s=1}^T \sum_{v=1}^p \sum_{l=1}^q
    \b{E} \Big[ E_{t,hu} \bf{n}'\F_t \bf{g}
    \cdot E_{s,vl} \bf{j}' \F_s' \bf{n}\Big]
    =
    \sum_{s=1}^T \sum_{v=1}^p \sum_{l=1}^q
    \text{cov}(E_{t,hu}, E_{s,vl})
    \cdot \b{E}\Big[ \bf{n}' \F_t \bf{g}
    \bf{j}' \F_s' \bf{n} \Big] \\
    &=
    \sum_{v=1}^p \sum_{l=1}^q \Big[
    O(\A_{e,c,l\cdot}' \A_{e,c,u\cdot} \A_{e,r,h\cdot}' \A_{e,r,v\cdot}) \cdot
    \sum_{w\geq 0}\sum_{m\geq 0}\sum_{s=1}^T
    a_{e,w}a_{e,w+|t-s|} a_{f,m}a_{f,m+|t-s|} \\
    & + O(\b{1}_{\{h=v\}} \b{1}_{\{u=l\}} \cdot (\bSigma_{\epsilon})_{hu}^2)\cdot
    \sum_{w\geq 0}\sum_{m\geq 0}\sum_{s=1}^T a_{\epsilon,w} a_{\epsilon,w+|t-s|} a_{f,m}a_{f,m+|t-s|} \Big] \\
    &=
    \sum_{v=1}^p \sum_{l=1}^q
    O(\A_{e,c,l\cdot}' \A_{e,c,u\cdot} \A_{e,r,h\cdot}' \A_{e,r,v\cdot} +
    \b{1}_{\{h=v\}} \b{1}_{\{u=l\}} \cdot (\bSigma_{\epsilon})_{hu}^2)
    = O(1) ,
\end{split}
\end{equation*}
where for the second last equality, we argue for the first term in the second last line only, as the second term could be shown similarly:
\begin{equation*}
\begin{split}
    &\hspace{5mm}
    \sum_{w\geq 0} \sum_{m\geq 0} \sum_{s=1}^T
    a_{e,w} a_{e,w+|t-s|} a_{f,m} a_{f,m+|t-s|} =
    \sum_{w\geq 0} \sum_{m\geq 0}
    a_{e,w} a_{f,m} \sum_{s=1}^T
    a_{e,w+|t-s|}a_{f,m+|t-s|} \\
    &\leq
    \sum_{w\geq 0} \sum_{m\geq 0}
    |a_{e,w}|\cdot|a_{f,m}| \cdot
    \Big(\sum_{s=1}^T a_{e,w+|t-s|}^2 \Big)^{1/2}
    \Big(\sum_{s=1}^T a_{f,m+|t-s|}^2\Big)^{1/2}
    \leq
    \sum_{w\geq 0} \sum_{m\geq 0}
    |a_{e,w}| \cdot|a_{f,m}| \leq c^2,
\end{split}
\end{equation*}
where the constant $c$ is from Assumptions (F1) and (E2). Finally,
\begin{equation}
\label{eqn: I6_step}
\begin{split}
    &\hspace{5mm}
    \b{E}\Big\{\Big( \sum_{h=1}^p \sum_{t=1}^T (\C_t\E_t')_{ih} \Big)^2\Big\} =
    \b{E}\Big\{\Big( \sum_{u=1}^q \sum_{h=1}^p \sum_{t=1}^T  E_{t,hu} \A_{c,u\cdot}'\F_t' \A_{r,i\cdot} \Big)^2\Big\} \\
    &=
    \sum_{t=1}^T \sum_{h=1}^p \sum_{u=1}^q
    \sum_{s=1}^T \sum_{v=1}^p \sum_{l=1}^q
    \b{E} \Big[ E_{t,hu} \A_{r,i\cdot}'\F_t \A_{c,u\cdot} \cdot E_{s,vl} \A_{c,l\cdot}' \F_s' \A_{r,i\cdot} \Big]
    = O(Tpq).
\end{split}
\end{equation}
Thus, we have
\begin{equation*}
\begin{split}
    &\hspace{5mm}
    \|\cI_6\|^2 \leq
    \frac{1}{T^2 p^2} \|\wh\D_r^{-1}\|_F^2 \cdot \|\H_r\|_F^2 \cdot \Big\| \sum_{i=1}^p
    \Q_{r,i\cdot} \sum_{t=1}^T (\C_t\E_t'\1_p\1_p')_{ij} \Big\|^2 \\
    &=
    \frac{1}{T^2 p^2} \|\wh\D_r^{-1}\|_F^2 \cdot \|\H_r\|_F^2 \cdot \Big\| \sum_{h=1}^p \sum_{i=1}^p \Q_{r,i\cdot} \sum_{t=1}^T (\C_t\E_t')_{ih} (\1_p\1_p')_{hj} \Big\|^2 \\
    &\leq
    \frac{1}{T^2 p^2} \|\wh\D_r^{-1}\|_F^2 \cdot \|\H_r\|_F^2 \cdot \|\Q_{r} \|_F^2 \cdot
    \sum_{i=1}^p \Big( \sum_{h=1}^p \sum_{t=1}^T (\C_t\E_t')_{ih} \Big)^2 \\
    &=
    O_P\big( T^{-2} p^{-2-2\delta_{r,k_r}} q^{-2\delta_{c,1}} \cdot Tp^2q \big)
    = O_P\big(T^{-1} p^{-2\delta_{r,k_r}} q^{1-2\delta_{c,1}} \big) ,
\end{split}
\end{equation*}
where we used Cauchy-Schwarz inequality in the third line, both Lemma \ref{lemma: norm_hat_D} and (\ref{eqn: I6_step}) in the second last equality. Similarly, we have the following for $\cI_7$,
\begin{equation*}
\begin{split}
    &\hspace{5mm}
    \|\cI_7\|^2 \leq
    \frac{1}{T^2 p^2} \|\wh\D_r^{-1}\|_F^2 \cdot \|\H_r\|_F^2 \cdot \Big\| \sum_{i=1}^p
    \Q_{r,i\cdot} \sum_{t=1}^T (\1_p\1_p'\E_t\C_t')_{ij} \Big\|^2 \\
    &=
    \frac{1}{T^2 p^2} \|\wh\D_r^{-1}\|_F^2 \cdot \|\H_r\|_F^2 \cdot \Big\| \sum_{h=1}^p \sum_{i=1}^p \Q_{r,i\cdot} \sum_{t=1}^T (\1_p\1_p')_{ih} (\E_t\C_t')_{hj} \Big\|^2 \\
    &\leq
    \frac{1}{T^2 p^2} \|\wh\D_r^{-1}\|_F^2 \cdot \|\H_r\|_F^2 \cdot \|\Q_{r}\|_F^2 \cdot
    p\Big( \sum_{h=1}^p \sum_{t=1}^T (\C_t\E_t')_{jh} \Big)^2 \\
    &=
    O_P\big( T^{-2} p^{-2-2\delta_{r,k_r}} q^{-2\delta_{c,1}} \cdot Tp^2q \big)
    = O_P\big(T^{-1} p^{-2\delta_{r,k_r}} q^{1-2\delta_{c,1}} \big) .
\end{split}
\end{equation*}

For $\cI_8$ and $\cI_9$, their rates can be shown to be the same as that for $\cI_3$ by the following,
\begin{equation*}
\begin{split}
    &\hspace{5mm}
    \|\cI_8\|^2 \leq
    \frac{1}{T^2 p^2} \|\wh\D_r^{-1}\|_F^2 \cdot \|\H_r\|_F^2 \cdot \Big\| \sum_{i=1}^p
    \Q_{r,i\cdot} \sum_{t=1}^T (\E_t\E_t'\1_p\1_p')_{ij} \Big\|^2 \\
    &=
    \frac{1}{T^2 p^2} \|\wh\D_r^{-1}\|_F^2 \cdot \|\H_r\|_F^2 \cdot \Big\| \sum_{h=1}^p \sum_{i=1}^p \sum_{l=1}^q \sum_{t=1}^T \Q_{r,i\cdot} E_{t,il} E_{t,hl} \Big\|^2 \\
    &\leq
    \frac{1}{T^2 p^2} \|\wh\D_r^{-1}\|_F^2 \cdot \|\H_r\|_F^2 \cdot p \sum_{h=1}^p \Big\| \sum_{i=1}^p \sum_{l=1}^q \sum_{t=1}^T \Q_{r,i\cdot} E_{t,il} E_{t,hl} \Big\|^2 \\
    &=
    O_P\big(T^{-2} p^{-2-2\delta_{r,k_r}} q^{-2\delta_{c,1}} \cdot p^2 (Tq + T^2 q^2 p^{-\delta_{r,k_r}}) \big)
    = O_P\big(T^{-1} p^{-2\delta_{r,k_r}} q^{1-2\delta_{c,1}} + p^{-3\delta_{r,k_r}} q^{2-2\delta_{c,1}} \big) ,
\end{split}
\end{equation*}
where we used Lemma \ref{lemma: norm_hat_D} and (\ref{eqn: I3_last_step}) in the second last equality. The proof for $\|\cI_9\|^2$ is similar to the above by using the proof of (\ref{eqn: I3_last_step}) previously and omitted here.

For $\cI_{10}$, first observe from the proof of (\ref{eqn: E11_bound_last_step}), we also have for any $j\in[p]$,
\[
\sum_{i=1}^p \b{E}\Big\{\Big( \sum_{t=1}^T \sum_{l=1}^q E_{t,il} \sum_{h=1}^q E_{t,jh} \Big)^2\Big\} = O\big( Tpq^2 + T^2q^2 \big) ,
\]
then together with Lemma \ref{lemma: norm_hat_D}, it holds that
\begin{equation*}
\begin{split}
    &\hspace{5mm}
    \|\cI_{10}\|^2 \leq
    \frac{1}{T^2 p^2} \|\wh\D_r^{-1}\|_F^2 \cdot \|\H_r\|_F^2 \cdot \Big\| \sum_{i=1}^p
    \Q_{r,i\cdot} \sum_{t=1}^T (\E_t\1_q\1_q'\E_t')_{ij} \Big\|^2 \\
    &=
    \frac{1}{T^2 p^2} \|\wh\D_r^{-1}\|_F^2 \cdot \|\H_r\|_F^2 \cdot \Big\| \sum_{i=1}^p
    \Q_{r,i\cdot} \sum_{t=1}^T
    \sum_{l=1}^q E_{t,il} \sum_{h=1}^q E_{t,jh} \Big\|^2 \\
    &\leq
    \frac{1}{T^2 p^2} \|\wh\D_r^{-1}\|_F^2 \cdot \|\H_r\|_F^2 \cdot \|\Q_r\|_F^2 \cdot \sum_{i=1}^p \Big(\sum_{t=1}^T \sum_{l=1}^q E_{t,il} \sum_{h=1}^q E_{t,jh} \Big)^2 \\
    &=
    O_P\big(T^{-2} p^{-2-2\delta_{r,k_r}} q^{-2\delta_{c,1}} \cdot (Tpq^2 + T^2q^2)\big)
    = O_P\big(
    T^{-1} p^{-1-2\delta_{r,k_r}} q^{2-2\delta_{c,1}} + p^{-2-2\delta_{r,k_r}} q^{2-2\delta_{c,1}} \big) .
\end{split}
\end{equation*}

Consider now $\cI_{11}$, we have
\begin{equation*}
\begin{split}
    &\hspace{5mm}
    \|\cI_{11}\|^2 \leq
    \frac{1}{T^2 p^4 q^2} \|\wh\D_r^{-1}\|_F^2 \cdot \|\H_r\|_F^2 \cdot \Big\| \sum_{i=1}^p
    \Q_{r,i\cdot} \sum_{t=1}^T (\1_q'\E_t'\1_p)^2 \Big\|^2 \\
    &\leq
    \frac{1}{T^2 p^4 q^2} \|\wh\D_r^{-1}\|_F^2 \cdot \|\H_r\|_F^2 \cdot \|\Q_r\|_F^2 \cdot p \Big[\sum_{t=1}^T (\1_q'\E_t'\1_p)^2 \Big]^2 \\
    &=
    O_P\big( T^{-2} p^{-4 -2\delta_{r,k_r}} q^{-2 -2\delta_{c,1}} \cdot T^2 p^3 q^2\big)
    = O_P\big( p^{-1 -2\delta_{r,k_r}} q^{ -2\delta_{c,1}} \big) ,
\end{split}
\end{equation*}
where we used Lemma \ref{lemma: norm_hat_D} and the rate from (\ref{eqn: 1E1_bound}) in the second last equality.

For $\cI_{12}$, we have
\begin{equation*}
\begin{split}
    \|\cI_{12}\|^2 &\leq
    \frac{1}{T^2 p^4} \|\wh\D_r^{-1}\|_F^2 \cdot \|\H_r\|_F^2 \cdot \Big\| \sum_{i=1}^p
    \Q_{r,i\cdot} \sum_{t=1}^T \1_p'\E_t\E_t'\1_p \Big\|^2 \\
    &\leq
    \frac{1}{T^2 p^4} \|\wh\D_r^{-1}\|_F^2 \cdot \|\H_r\|_F^2 \cdot \|\Q_r\|_F^2 \cdot p \Big(\sum_{t=1}^T \1_p'\E_t\E_t'\1_p \Big)^2 \\
    &=
    O_P\big( T^{-2} p^{-4 -2\delta_{r,k_r}} q^{-2\delta_{c,1}} \cdot T^2p^3q^2 \big)
    = O_P\big( p^{-1 -2\delta_{r,k_r}} q^{2-2\delta_{c,1}} \big) ,
\end{split}
\end{equation*}
where the last line used the following result which can be shown similar to (\ref{eqn: cov_EEEE}),
\begin{equation*}
\begin{split}
    &\hspace{5mm}
    \b{E}\Big\{ \Big(\sum_{t=1}^T \1_p'\E_t\E_t'\1_p \Big)^2 \Big\} =
    \b{E}\Big\{ \Big(\sum_{t=1}^T \sum_{h=1}^q \sum_{i=1}^p \sum_{j=1}^p E_{t,ih} E_{t,jh} \Big)^2 \Big\} \\
    &=
    \sum_{t=1}^T \sum_{h=1}^q \sum_{i=1}^p \sum_{j=1}^p \sum_{s=1}^T \sum_{l=1}^q \sum_{m=1}^p \sum_{n=1}^p \text{cov}( E_{t,ih} E_{t,jh}, E_{s,ml} E_{s,nl}) + \Big( \sum_{t=1}^T \sum_{h=1}^q \sum_{i=1}^p \sum_{j=1}^p \b{E}[E_{t,ih} E_{t,jh}]\Big)^2 \\
    &=
    O\big(Tp^2q^2 + T^2p^2q^2 \big) = O\big( T^2p^2q^2 \big).
\end{split}
\end{equation*}

Lastly, $\|\cI_{13}\|^2$ is dominated by the terms from $\cI_1$ to $\cI_{12}$ using Cauchy-Schwarz inequality and Theorem \ref{thm: consistency}. We require the term $\cI_1$ to be truly dominating by using Assumption (AD1) and we equivalently compare the rates without the term $\wh\D_r^{-1}$. Notice the rates for $\|\cI_2\|^2$, $\|\cI_4\|^2$, $\|\cI_5\|^2$, $\|\cI_6\|^2$, $\|\cI_7\|^2$, $\|\cI_{10}\|^2$, $\|\cI_{11}\|^2$ and $\|\cI_{12}\|^2$ are bounded above by the rate for $\|\cI_3\|^2$ which is the same as $\|\cI_8\|^2$ and $\|\cI_9\|^2$. Thus, it suffices to consider the following ratio as $p,q,T \to \infty$,
\begin{align*}
    & \|\cI_3\|^2 / \|\cI_1\|^2 = O_P\big(
    1/ p^{\delta_{r,1}} + Tq/ p^{\delta_{r,1} +\delta_{r,k_r}} \big) = o_P(1),
\end{align*}
by the rate assumption $Tq= o(p^{\delta_{r,1}+ \delta_{r,k_r}})$. Therefore, $\cI_1$ is dominating over other terms in (\ref{eqn: Q_asymp_decomp}) and hence we have
\begin{equation}
\label{eqn: I1_dominate}
\begin{split}
    \wh\Q_{r,j\cdot} - \H_r\Q_{r,j\cdot} &=
    \cI_1 + o_P(1) = \frac{1}{T} \wh\D_r^{-1} \sum_{i=1}^p \H_r \Q_{r,i\cdot} \sum_{t=1}^T (\C_t \E_t')_{ij} + o_P(1) \\
    &\xrightarrow{p}
    \frac{1}{T} \D_r^{-1} \H_r^\ast \sum_{i=1}^p \Q_{r,i\cdot} \sum_{t=1}^T (\C_t \E_t')_{ij} ,
\end{split}
\end{equation}
where the last line used Lemma \ref{lemma: limit_D} and Lemma \ref{lemma: limit_H} in which $\D_r$ and $\H_r^\ast$ are defined, respectively. In the rest of the proof, we show
\begin{equation*}
    \sqrt{T \omega_B} \cdot \frac{1}{T} \D_r^{-1} \H_r^\ast \sum_{i=1}^p \Q_{r,i\cdot} \sum_{t=1}^T (\C_t \E_t')_{ij}
    \xrightarrow{\c{D}} \cN\big(\0, T^{-1} \omega_B \cdot \D_r^{-1} \H_r^\ast \bf{\Xi}_{r,j} (\H_r^\ast)' \D_r^{-1} \big),
\end{equation*}
with $\omega_B:= p^{2\delta_{r,k_r} - \delta_{r,1}} q^{2\delta_{c,1} -1}$ and $\bf{\Xi}_{r,j} := \plim_{p,q,T \to\infty} \text{Var} \big( \sum_{i=1}^p \Q_{r,i\cdot} \sum_{t=1}^T (\C_t \E_t')_{ij} \big)$. We will adapt the central limit theorem for $\alpha$-mixing processes as depicted in Theorem 2.21 in \cite{FanYao2003}. First, define $\B_{j,t} := \sqrt{\omega_B} \cdot \D_r^{-1} \H_r^\ast \sum_{i=1}^p \Q_{r,i\cdot} (\C_t \E_t')_{ij}$, and also let $\bf{b}_{e,t} :=\sum_{w\geq 0} a_{e,w} \X_{e,t-w}$, $\bf{b}_{\epsilon,il,t} :=\sum_{w\geq 0} a_{\epsilon,w}(\X_{\epsilon,t-w})_{il}$ and $\bf{b}_{f,t} :=\sum_{w\geq 0}a_{f,w}\X_{f,t-w}$ which are independent of each other by Assumption (E2).

Since we may write $\B_{j,t} = h(\bf{b}_{e,t}, (\bf{b}_{\epsilon,il,t})_{i\in [p], l\in[q]}, \bf{b}_{f,t})$ for some function $h$, we conclude $\B_{j,t}$ is $\alpha$-mixing using Theorem 5.2 in \cite{Bradley2005}. Observe that $\b{E}[\B_{j,t}]=\0$, and we show in the following that there exists an $m>2$ such that $\b{E}[\|\B_{j,t}\|^m] \leq C$ for some constant $C$,
\begin{equation*}
\begin{split}
    \b{E}[\|\B_{j,t}\|^m] & \leq
    \omega_B^{m/2} \cdot \|\D_r^{-1}\|_F^m \cdot \|\H_r^{\ast}\|_F^m
    \cdot \|\Q_r\|_F^m \cdot
    \b{E}\bigg\{ \Big[\sum_{i=1}^{p}
    \Big(\sum_{l=1}^{q} E_{t,jl} \A_{c,h\cdot}'
    \F_t' \A_{r,i\cdot} \Big)^2 \Big]^{m/2}\bigg\} \\
    &= O\Big( (p^{2 \delta_{r,k_r}} q^{2 \delta_{c,1}})^{m/2}\Big) \cdot\|\D_r^{-1}\|_F^m
    = O(1),
\end{split}
\end{equation*}
where we used Lemma \ref{lemma: correlation_Et_Ft}.2 and the definition of $\omega_B$ in the second last equality, and Lemma \ref{lemma: limit_D} in the last equality. Theorem 2.21 in \cite{FanYao2003} then applies, and hence
\[
\sqrt{T \omega_B} \cdot \frac{1}{T} \D_r^{-1} \H_r^\ast \sum_{i=1}^p \Q_{r,i\cdot} \sum_{t=1}^T (\C_t \E_t')_{ij}
= \frac{1}{\sqrt{T}} \sum_{t=1}^T \B_{j,t}
\xrightarrow{\c{D}} \cN\big(\0, T^{-1} \omega_B \cdot \D_r^{-1} \H_r^\ast \bf{\Xi}_{r,j} (\H_r^\ast)' \D_r^{-1} \big).
\]
Together with (\ref{eqn: I1_dominate}), we arrive at
\begin{equation}
\label{eqn: asymp_loading}
    (T p^{2\delta_{r,k_r} - \delta_{r,1}} q^{2\delta_{c,1} -1})^{1/2} \cdot (\wh\Q_{r,j\cdot} - \H_r\Q_{r,j\cdot})
    \xrightarrow{\c{D}} \cN\big(\0, T^{-1} p^{2 \delta_{r,k_r} -\delta_{r,1}} q^{2\delta_{c,1} -1} \cdot \D_r^{-1} \H_r^\ast \bf{\Xi}_{r,j} (\H_r^\ast)' \D_r^{-1} \big) .
\end{equation}
This completes the proof of Theorem \ref{thm: asymp_loading}. $\square$

\textbf{\textit{Proof of Theorem \ref{thm: HAC}.}}
We only prove the scenario for $Tq= o(p^{\delta_{r,1}+ \delta_{r,k_r}})$, which is on showing the constructed estimator for the row loading estimator is consistent. Note $\wh\D_r$ consistently estimates $\D_r$ by Lemma \ref{lemma: limit_D}, and $\H_r$ consistently estimates $\H_r^\ast$ by Lemma \ref{lemma: limit_H}. Then observe that
\begin{equation*}
\begin{split}
    &\hspace{5mm}
    \H_r \textnormal{Var} \Big\{ \sum_{i=1}^p \Q_{r,i\cdot} \sum_{t=1}^T (\C_t \E_t')_{ij} \Big\} \H_a' = \textnormal{Var} \Big\{ \sum_{i=1}^p \H_r \Q_{r,i\cdot} \sum_{t=1}^T (\C_t \E_t')_{ij} \Big\} \\
    &=
    \textnormal{Var} \Big\{ \sum_{i=1}^p \Big( T^{-1}\wh\D_r^{-1} \wh\Q_r' \Q_r \sum_{s=1}^T \F_{Z,s} \Q_c' \Q_c \F_{Z,s}' \Big) \Q_{r,i\cdot} \sum_{t=1}^T (\C_t \E_t')_{ij} \Big\} \\
    &=
    \textnormal{Var} \Big\{  \sum_{t=1}^T \sum_{i=1}^p \Big( T^{-1} \wh\D_r^{-1} \wh\Q_r' \sum_{s=1}^T \C_s \C_{s,i\cdot} \Big) (\C_t \E_t')_{ij} \Big\} .
\end{split}
\end{equation*}
By Theorem \ref{thm: consistency}, $(\wh\mu_t, \wh\balpha_t, \wh\bbeta_t)$ is consistent for $(\mu_t, \balpha_t, \bbeta_t)$. By Theorem \ref{thm: consistency_F_C} and the rate assumption in the statement of Theorem \ref{thm: HAC}, $\wh\C_t$ is consistent for $\C_t$ and hence $\wh\E_t$ is consistent for $\E_t$. Thus, we conclude that $\wh\bSigma_{r,j}^{HAC}$ is estimating $\H_r \bf{\Xi}_{r,j} \H_r'$ consistently (\cite{NeweyWest1987}) and hence result 1 is implied. Result 2 then follows, and results 3 and 4 can be shown similarly (details omitted). This completes the proof of the Theorem \ref{thm: HAC}. $\square$

\textbf{\textit{Proof of Theorem \ref{thm: testingFMvsMEFM}.}}
Combining (\ref{eqn: Ct_estimator}) and (\ref{eqn: Ehat_t}), we have
\begin{align*}
  \wh\E_t &= \wh\L_t - \wh\C_t = \M_p\Y_t\M_q - \wh\C_t = \M_p\C_t\M_q + \M_p\E_t\M_q -\wh\C_t\\
  &= (\C_t-\wh\C_t) + \M_p\A_{e,r}\F_{e,t}\A_{e,c}'\M_q + \M_p(\bSigma_\epsilon\circ\bepsilon_t)\M_q,
\end{align*}
where the second line comes from (IC1) being satisfied, so that $\M_p\A_r = \A_r$ and $\M_q\A_c = \A_c$. Hence
\begin{align*}
  q^{-1}(\wh\E_t\wh\E_t')_{ii} - q^{-1}\sum_{j=1}^q\Sigma_{\epsilon,ij}^2
  &= \sum_{i=1}^6 I_i, \text{ where }\\
  I_1 &:= q^{-1}\{\M_p(\bSigma_\epsilon\circ\bepsilon_t)\M_q(\bSigma_\epsilon\circ\bepsilon_t)'\M_p\}_{ii}
        -q^{-1}\sum_{j=1}^q\Sigma_{\epsilon,ij}^2,\\
  I_2 &:=  q^{-1}\{(\wh\C_t-\C_t)(\wh\C_t-\C_t)'\}_{ii},\\
  I_3 &:= q^{-1}\{\M_p\A_{e,r}\F_{e,t}\A_{e,c}'\M_q\A_{e,c}\F_{e,t}'\A_{e,r}'\M_p\}_{ii},\\
  I_4 &:= O_P(q^{-1}\{(\C_t-\wh\C_t)\M_q\A_{e,c}\F_{e,t}'\A_{e,r}'\M_p\}_{ii}),\\
  I_5 &:= O_P(q^{-1}\{(\C_t-\wh\C_t)\M_q(\bSigma_\epsilon\circ\bepsilon_t)'\M_p\}_{ii}),\\
  I_6 &:= O_P(q^{-1}\{\M_p\A_{e,r}\F_{e,t}\A_{e,c}'\M_q(\bSigma_\epsilon\circ\bepsilon_t)'\M_p\}_{ii}).
\end{align*}
By an assumption in the statement of the theorem, we have $I_2 = o_P(q^{-1})$. Since $\|\M_p\A_{e,r}\|_1 \leq \|\M_p\|_1\|\A_{e,r}\|_1 = O(1)$, with the finiteness of $k_r$ and $k_c$, we have $I_3, I_6 = O_P(q^{-1})$, and by the Cauchy-Schwarz inequality, $I_4 = O_P(I_2^{1/2}\cdot I_6^{1/2}) = o_P(q^{-1})$.
Writing $\eta_{t,ij} := (\bSigma_\epsilon\circ\bepsilon_t)_{ij}$, we can expand
\begin{align*}
  I_1 &= \frac{1}{q}\sum_{j=1}^q\Sigma_{\epsilon,ij}^2(\epsilon_{t,ij}^2-1) -
  \frac{2}{q}\sum_{j=1}^q\bigg(\frac{1}{p}\sum_{\ell=1}^p\eta_{t,\ell j}\bigg)\eta_{t,ij}\\
  &\quad +  \frac{2}{q}\bigg(\frac{\1_p'\boldeta_t\1_q}{pq}\bigg)\sum_{j=1}^q\eta_{t,ij}
  + \frac{1}{q}\sum_{j=1}^q\bigg(\frac{1}{p}\sum_{i=1}^p\eta_{t,ij}\bigg)^2 - \bigg(\frac{\1_p'\boldeta_t\1_q}{pq}\bigg)^2 - \bigg(\frac{1}{q}\sum_{j=1}^q\eta_{t,ij}\bigg)^2,\\
  &= \frac{1}{q}\sum_{j=1}^q\Sigma_{\epsilon,ij}^2(\epsilon_{t,ij}^2-1) - \frac{2}{pq}\sum_{j=1}^q\eta_{t,ij}^2 - \frac{2}{pq}\sum_{j=1}^q\sum_{\ell\neq i}\eta_{t,\ell j}\eta_{t,ij}\\
  &\quad + O_P(q^{-1}p^{-1/2}) + O_P(p^{-1}) + O_P(p^{-1}q^{-1}) + O_P(q^{-1})\\
  &= \frac{1}{q}\sum_{j=1}^q\Sigma_{\epsilon,ij}^2(\epsilon_{t,ij}^2-1)(1 + o_P(1)),
\end{align*}
where all rates of convergence above are obtained from applying the Markov inequality. Hence
\begin{align*}
  \frac{(\wh\E_t\wh\E_t')_{ii} - \sum_{j=1}^q\Sigma_{\epsilon,ij}^2}{\sqrt{\sum_{j=1}^q\Sigma_{\epsilon,ij}^4\var(\epsilon_{t,ij}^2)}}
  &= \frac{q^{-1}(\wh\E_t\wh\E_t')_{ii} - q^{-1}\sum_{j=1}^q\Sigma_{\epsilon,ij}^2}{\sqrt{q^{-2}\sum_{j=1}^q\var(\epsilon_{t,ij}^2)\Sigma_{\epsilon,ij}^4}}\\
  &= \frac{q^{-1}\sum_{j=1}^q\Sigma_{\epsilon,ij}^2(\epsilon_{t,ij}^2-1)(1+o_P(1))}
  {\sqrt{q^{-2}\sum_{j=1}^q\var(\epsilon_{t,ij}^2)\Sigma_{\epsilon,ij}^4}}\\
  &\xrightarrow{\c{D}} Z_i \xrightarrow{\c{D}} \c{N}(0,1),
\end{align*}
where the last line follows from Theorem 1 in \cite{AyvazyanUlyanov2023}, with
\[Z_i := \frac{q^{-1}\sum_{j=1}^q\Sigma_{\epsilon,ij}^2(\epsilon_{t,ij}^2-1)}
  {\sqrt{q^{-2}\sum_{j=1}^q\var(\epsilon_{t,ij}^2)\Sigma_{\epsilon,ij}^4}}, \;\;\; i\in[p],\]
so that we can easily see that the $Z_i$'s are independent of each other by Assumption (E1). The proof for $(\wh\E_t\wh\E_t')_{ii}$ completes since from the calculations for $I_1$, we see that
\[I_5 = O_P(I_2^{1/2}\cdot 1) = o_P(q^{-1/2}).\]
For $(\check\E_t\check\E_t')_{ii}$ under $H_0$, note that we have
\[\check\E_t = (\C_t - \check\C_t) + \A_{e,r}\F_{e,t}\A_{e,c}' + (\bSigma_{\epsilon}\circ\bepsilon_t),\]
with the rate for $\check{C}_{t,ij} - C_{t,ij}$ the same as that for $\wh{C}_{t,ij} - C_{t,ij}$ since the estimation procedure for FM is essentially the same with the same assumptions on the factor loadings (apart from (IC1) which is not important for FM), the factors and the noise. Hence the proof we employed so far can be replicated with $\M_p$ and $\M_q$ replaced by the corresponding sized identity matrices, and we arrive at the same conclusion with the same $Z_i$'s.

For $(\wh\E_t'\wh\E_t)_{jj}$ under both $H_0$ and $H_1$ and $(\check\E_t'\check\E_t)_{jj}$ under $H_0$, the proofs are parallel to that for $(\wh\E_t\wh\E_t')_{ii}$, and we omit them here. $\square$

\textbf{\textit{Proof of Corollary \ref{cor: test_statistic_FMvsMEFM}.}}
Since both $(\wh\E_t\wh\E_t)_{ii}$ and $(\check\E_t\check\E_t')_{ii}$ are distributed asymptotically the same under $H_0$ by Theorem \ref{thm: testingFMvsMEFM} for each $i\in[p]$ and $t\in[T]$, we have that $x_{\alpha,t}$ and $y_{\alpha,t}$ are asymptotically distributed the same for each $t\in[T]$ under $H_0$. Hence $\wh{q}_{x,\alpha}$ and $\wh{q}_{y,\alpha}$ are asymptotically the same where
$\wh{q}_{y,\alpha}(\theta) := \inf\{c \; | \; \b{F}_{y,\alpha}(c)\geq \theta\}$. This implies that under $H_0$, for each $t\in[T]$,
\begin{align*}
1- \theta \geq \b{P}_{y,\alpha}(y_{\alpha,t} \geq \wh{q}_{y,\alpha}(\theta)) \rightarrow \b{P}_{y,\alpha}(y_{\alpha,t} \geq \wh{q}_{x,\alpha}(\theta)).
\end{align*}
The other inequality is similarly proved, and is omitted here. $\square$

\textbf{\textit{Proof of Theorem \ref{thm: rank}.}}
First consider $\wh{k}_r$, i.e. result 1 in Theorem \ref{thm: rank}. For $j\in[k_r]$, first consider
\begin{equation}
\label{eqn: lambda_j_sample_cov_mat}
\begin{split}
    &\hspace{5mm}
    \lambda_j\Big(\frac{1}{T}\sum_{t=1}^T \C_t\C_t' \Big)
    = \lambda_j \Big(\frac{1}{T}\sum_{t=1}^T \A_r \F_t \A_c' \A_c \F_t' \A_r' \Big)
    = \lambda_j \Big(\A_r'\A_r \cdot \frac{1}{T}\sum_{t=1}^T \F_t \A_c' \A_c \F_t' \Big) \\
    &\asymp_P
    \tr(\A_c' \A_c)\cdot \lambda_j(\A_r'\A_r) = \|\A_c\|_F^2 \cdot \lambda_j (\bSigma_{A,r}^{1/2} \Z_r \bSigma_{A,r}^{1/2})
    \asymp q^{\delta_{c,1}} \cdot \lambda_j(\Z_r)
    = p^{\delta_{r,j}} q^{\delta_{c,1}} ,
\end{split}
\end{equation}
where in the second line, we used Assumption (F1) in the first step and Assumption (L1) in the second. For the second last step, we used Theorem 1 of \cite{Ostrowsk1959} on the eigenvalues of a congruent transformation $\bSigma_{A,r}^{1/2} \Z_r \bSigma_{A,r}^{1/2}$ of $\Z_r$, where we further used Assumption (L1) that all eigenvalues of $\bSigma_{A,r}$ are bounded away from $0$ and infinity.

Since we have $T^{-1} \sum_{t=1}^T \wh\L_t \wh\L_t' = T^{-1} \sum_{t=1}^T \C_t\C_t' + T^{-1} \sum_{t=1}^T \R_{r,t}$ from (\ref{eqn: R_rt}), it holds by Weyl's inequality that for $j\in[k_r]$,
\begin{equation}
\label{eqn: lambda_diff_LL_CC}
\Big|\lambda_j\Big(\frac{1}{T} \sum_{t=1}^T \wh\L_t \wh\L_t'\Big) - \lambda_j\Big(\frac{1}{T} \sum_{t=1}^T \C_t\C_t'\Big) \Big| \leq
\Big\| \frac{1}{T} \sum_{t=1}^T \R_{r,t} \Big\| = o_P(\omega_r),
\end{equation}
where $\omega_r = p^{\delta_{r,k_r}} q^{\delta_{c,1}}$ as defined in Lemma \ref{lemma: norm_hat_D}, and the last equality used $\bgamma'(T^{-1} \sum_{t=1}^T \R_{r,t})\bgamma = o_P(\omega_r)$ for any unit vector $\bgamma\in\b{R}^p$ from the proof of Lemma \ref{lemma: norm_hat_D}. With our choice of $\xi_r$, we also have
\begin{equation}
\label{eqn: xi_over_omega}
\xi_r/ \omega_r \asymp
p^{1-\delta_{r,k_r}} q^{1-\delta_{c,1}}\big[(Tq)^{-1/2} + p^{-1/2}\big] = o(1),
\end{equation}
where we used Assumption (R1) in the last equality. Hence for $k_r>1$, if $j\in[k_r-1]$, using (\ref{eqn: lambda_diff_LL_CC}) and (\ref{eqn: xi_over_omega}) we have
\begin{equation}
\label{eqn: ratio_lambda_j_plus_xi}
\begin{split}
    &\hspace{5mm}
    \frac{\lambda_{j+1}\big(T^{-1} \sum_{t=1}^T \wh\L_t\wh\L_t' \big) + \xi_r}{\lambda_{j}\big(T^{-1} \sum_{t=1}^T \wh\L_t\wh\L_t' \big) + \xi_r} \\
    &\leq
    \frac{\lambda_{j+1}\big(T^{-1} \sum_{t=1}^T \C_t\C_t' \big) + \xi_r + \big|\lambda_{j+1}\big(T^{-1} \sum_{t=1}^T \wh\L_t \wh\L_t'\big) - \lambda_{j+1} \big(T^{-1} \sum_{t=1}^T \C_t\C_t'\big)\big|}{\lambda_{j}\big(T^{-1} \sum_{t=1}^T \C_t\C_t' \big) + \xi_r - \big|\lambda_j\big(T^{-1} \sum_{t=1}^T \wh\L_t \wh\L_t'\big) - \lambda_j\big(T^{-1} \sum_{t=1}^T \C_t\C_t'\big)\big|} \\
    &=
    \frac{\lambda_{j+1}\big(T^{-1} \sum_{t=1}^T \C_t\C_t' \big) + o_P(\omega_r)}{\lambda_{j}\big(T^{-1} \sum_{t=1}^T \C_t\C_t' \big) + o_P(\omega_r)}
    =
    \frac{\lambda_{j+1}\big(T^{-1} \sum_{t=1}^T \C_t\C_t' \big)}{\lambda_{j}\big(T^{-1} \sum_{t=1}^T \C_t\C_t' \big)} \big(1+o_P(1) \big)
    \asymp_P p^{\delta_{r,j+1} - \delta_{r,j}} ,
\end{split}
\end{equation}
where the last line used (\ref{eqn: lambda_j_sample_cov_mat}). Moreover, for any $j\in[k_r-1]$,
\begin{equation}
\label{eqn: ratio_lambda_kr_plus_xi}
\begin{split}
    &\hspace{5mm}
    \frac{\lambda_{k_r +1}\big(T^{-1} \sum_{t=1}^T \wh\L_t\wh\L_t' \big) + \xi_r}{\lambda_{k_r}\big(T^{-1} \sum_{t=1}^T \wh\L_t\wh\L_t' \big) + \xi_r}
    = \frac{\lambda_{k_r +1}\big(T^{-1} \sum_{t=1}^T \wh\L_t\wh\L_t' \big) + \xi_r}{\omega_r (1+ o_P(1))}
    = O_P\Big\{\frac{\lambda_{k_r +1}\big(T^{-1} \sum_{t=1}^T \wh\L_t\wh\L_t' \big)}{\omega_r} + \frac{\xi_r}{\omega_r}\Big\} \\
    &= O_P\Big\{O_P\Big( pq[(Tq)^{-1/2} + p^{-1/2}]/\omega_r \Big) + \xi_r/\omega_r \Big\} =  O_P(\xi_r/\omega_r) = o_P(p^{\delta_{r,j+1} - \delta_{r,j}}) ,
\end{split}
\end{equation}
where we used (\ref{eqn: xi_over_omega}) and the proof of Lemma \ref{lemma: norm_hat_D} in the first equality, our choice of $\xi_r$ in the second last, and the extra rate assumption in the statement of the theorem in the last. In the third equality, we used the following (which will be shown at the end of the this proof),
\begin{equation}
\label{eqn: lambda_after_kr_sample_cov_mat}
\lambda_j\Big(\frac{1}{T}\sum_{t=1}^T \C_t\C_t' \Big) = O_P(T^{-1/2}pq^{1/2} + p^{1/2}q), \;\;\;
j= k_r+1, \dots, p.
\end{equation}

Hence for $j= k_r+1, \dots, \lfloor p/2 \rfloor$ (true also for $k_r=1$),
\begin{equation}
\label{eqn: ratio_lambda_after_kr_plus_xi}
\frac{\lambda_{j+1}\big(T^{-1} \sum_{t=1}^T \wh\L_t\wh\L_t' \big) + \xi_r}{\lambda_{j}\big(T^{-1} \sum_{t=1}^T \wh\L_t\wh\L_t' \big) + \xi_r} \geq
\frac{\xi_r/ \omega_r}{O_P(\xi_r/\omega_r) + \xi_r/ \omega_r} \geq \frac{1}{C}
\end{equation}
in probability for some generic positive constant $C$, where we used again our choice of $\xi_r$ and (\ref{eqn: lambda_after_kr_sample_cov_mat}) in the first inequality. Combining (\ref{eqn: ratio_lambda_j_plus_xi}), (\ref{eqn: ratio_lambda_kr_plus_xi}) and (\ref{eqn: ratio_lambda_after_kr_plus_xi}), we may conclude our proposed $\wh{k}_r$ is consistent for $k_r$.

If $k_r=1$, then from our choice of $\xi_r$, (\ref{eqn: ratio_lambda_kr_plus_xi}) becomes
\[
\frac{\lambda_{k_r +1}\big(T^{-1} \sum_{t=1}^T \wh\L_t\wh\L_t' \big) + \xi_r}{\lambda_{k_r}\big(T^{-1} \sum_{t=1}^T \wh\L_t\wh\L_t' \big) + \xi_r} = O_P(\xi_r/ \omega_r) = o_P(1).
\]
Together with (\ref{eqn: ratio_lambda_after_kr_plus_xi}) which holds true for $k_r=1$, we can also conclude $\wh{k}_r$ is consistent for $k_r$.

It remains to show (\ref{eqn: lambda_after_kr_sample_cov_mat}). To this end, from (\ref{eqn: lambda_j_sample_cov_mat}) and (\ref{eqn: lambda_diff_LL_CC}), the first $k_r$ eigenvalues of $T^{-1}\sum_{t=1}^T \wh\L_t \wh\L_t'$ coincides with those of $T^{-1}\sum_{t=1}^T \C_t\C_t'$ asymptotically, so that the first $k_r$ eigenvectors corresponding to $T^{-1}\sum_{t=1}^T \wh\L_t \wh\L_t'$ coincides with those for $T^{-1}\sum_{t=1}^T \C_t\C_t'$ asymptotically as $T,p \rightarrow \infty$, which are necessarily in $\mathcal{N}^\perp:= \text{Span}(\Q_r)$, the linear span of the columns of $\Q_r$. This means that the $(k_r+1)$-th largest eigenvalue of $T^{-1}\sum_{t=1}^T \wh\L_t \wh\L_t'$ and beyond will asymptotically have eigenvectors in $\mathcal{N}$, the orthogonal complement of $\mathcal{N}^\perp$. Then for any unit vectors $\bgamma \in \mathcal{N}$, we have from (\ref{eqn: CC_decomposition}) and Lemma \ref{lemma: rate_R_rt} that
\begin{align*}
  \bgamma'\Big( T^{-1}\sum_{t=1}^T \wh\L_t \wh\L_t' \Big)\bgamma = \bgamma' \Big(\frac{1}{T}\sum_{t=1}^T \R_{r,t} \Big) \bgamma = O_P\Big(T^{-1} \Big \|\sum_{t=1}^T \R_{r,t} \Big\|_F \Big) = O_P(T^{-1/2}pq^{1/2} + p^{1/2}q),
\end{align*}
which is equivalent to (\ref{eqn: lambda_after_kr_sample_cov_mat}). This completes the proof of Theorem \ref{thm: rank}. $\square$

\subsection{Lemmas and proofs}
As we have the same factor structure as \cite{CenLam2024}, we list the following Lemma \ref{lemma: correlation_Et_Ft} for further use and refer readers to \cite{CenLam2024} for the proof in details.

\begin{lemma}\label{lemma: correlation_Et_Ft}
Let Assumptions (F1), (E1) and (E2) hold. Then
\begin{itemize}
    \item [1.]
    (Weak correlation of noise $\E_t$ across different rows, columns and times) there exists some positive constant $C<\infty$ so that for any $t\in[T], i,j\in[p], h\in[q]$,
    \begin{align*}
    & \sum_{k=1}^{p}\sum_{l=1}^{q}
    \Big|\b{E}[E_{t,ih}
    E_{t,kl}]\Big| \leq C, \\
    & \sum_{l=1}^{q} \sum_{s=1}^T \Big| \textnormal{cov}(E_{t,ih} E_{t,jh}, E_{s,il} E_{s,jl} ) \Big| \leq C.
    \end{align*}
    \item [2.]
    (Weak dependence between factor $\F_t$ and noise $\E_t$) there exists some positive constant $C < \infty$ so that for any $j\in[p], i\in[q]$, and any deterministic vectors $\bf{u}\in\b{R}^{k_r}$ and $\bf{v}\in\b{R}^{k_c}$ with constant magnitudes,
    \begin{equation*}
    \b{E}\Bigg(
        \frac{1}{(qT)^{1/2}}\sum_{h=1}^{q}
        \sum_{t=1}^T
        E_{t,jh}
        \bf{u}' \F_t \bf{v}
    \Bigg)^2 \leq C,
    \;\;\;
    \b{E}\Bigg(
        \frac{1}{(pT)^{1/2}}\sum_{h=1}^{p}
        \sum_{t=1}^T
        E_{t,hi}
        \bf{v}' \F_t' \bf{u}
    \Bigg)^2 \leq C .
    \end{equation*}
    \item [3.]
    (Further results on factor $\F_t$) for any $t\in[T]$, all elements in $\F_t$ are independent of each other, with mean $0$ and unit variance. Moreover,
    \begin{equation*}
    \frac{1}{T}\sum_{t=1}^T \F_t\F_t'
    \xrightarrow{p} \bSigma_r := k_c\I_{k_r}
    , \;\;\;
    \frac{1}{T}\sum_{t=1}^T \F_t'\F_t
    \xrightarrow{p} \bSigma_c := k_r\I_{k_c} ,
    \end{equation*}
    with the number of factors $k_r$ and $k_c$ fixed as $\min\{T,p,q\}\to \infty$.
\end{itemize}
\end{lemma}

\begin{lemma}\label{lemma: rate_R_rt}
(Bounding $\sum_{t=1}^T\R_{r,t}$) Under Assumptions (F1), (L1), (E1) and (E2), it holds that
\begin{align}
    & \Big\| \sum_{t=1}^T \C_t\E_t'\Big\|_F^2 =  O_P(T p^{1+\delta_{r,1}} q),
    \label{eqn: CE_bound} \\
    & \Big\| \sum_{t=1}^T \E_t\E_t' \Big\|_F^2 = O_P(T p^2 q + T^2 p q^2),
    \label{eqn: EE_bound} \\
    & \Big\| \sum_{t=1}^T \1_q' \E_t' \1_p \E_t \1_q \1_p' \Big\|_F^2 =
    O_P(T p^3 q^2 + T^2 p^2 q^2),
    \label{eqn: E11_bound} \\
    & \Big\| \sum_{t=1}^T \C_t \E_t' \1_p \1_p' \Big\|_F^2 = O_P(Tp^{3+\delta_{r,1}} q),
    \label{eqn: CE11_bound} \\
    & \Big\| \sum_{t=1}^T \E_t \E_t' \1_p \1_p' \Big\|_F^2 = O_P(Tp^4q + T^2 p^3 q^2),
    \label{eqn: EE11_bound} \\
    & \Big\| \sum_{t=1}^T \E_t\1_q\1_q' \E_t' \Big\|_F^2 = O_P(T p^2 q^2 + T^2 p q^2),
    \label{eqn: E11E_bound} \\
    & \Big\| \sum_{t=1}^T (\1_q' \E_t' \1_p)^2 \1_p \1_p' \Big\|_F^2 = O_P( T^2 p^4 q^2 ),
    \label{eqn: 1E111_bound} \\
    & \Big\| \sum_{t=1}^T \1_p' \E_t \E_t' \1_p \1_p \1_p' \Big\|_F^2 = O_P(T p^6 q + T^2 p^5 q^2) .
    \label{eqn: 1EE1_bound}
\end{align}
Thus, with $\R_{r,t}$ defined in (\ref{eqn: R_rt}), we have
\[
\Big\| \sum_{t=1}^T \R_{r,t} \Big\|_F^2 = O_P(Tp^2q + T^2pq^2) .
\]
\end{lemma}
\textbf{\textit{Proof of Lemma \ref{lemma: rate_R_rt}.}} Using $\C_t = \A_r \F_t \A_c'$, we have (\ref{eqn: CE_bound}) holds as follows,
\begin{equation*}
\begin{split}
    &\hspace{5mm}
    \Big\| \sum_{t=1}^T \C_t \E_t' \Big\|_F^2 =
    \Big\| \sum_{t=1}^T \A_r \F_t \A_c' \E_t' \Big\|_F^2 =
    \sum_{i=1}^p \sum_{l=1}^p
    \Bigg( \sum_{t=1}^T \A_{r, i\cdot}' \F_t \A_c' \E_{t, l\cdot} \Bigg)^2 \\
    &=
    \sum_{i=1}^p \|\A_{r, i\cdot}\|^2 \cdot \sum_{l=1}^p \Bigg( \sum_{h=1}^q \sum_{t=1}^T E_{t, lh} \frac{1}{\|\A_{r, i\cdot}\|}\A_{r, i\cdot}' \F_t \A_{c, h\cdot} \Bigg)^2
    = O_P\Big(T p^{1+\delta_{r,1}} q\Big) ,
\end{split}
\end{equation*}
where the last equality is from Assumption (L1) and Lemma \ref{lemma: correlation_Et_Ft}.

To show (\ref{eqn: EE_bound}), first notice from Assumption (E1),
\[
E_{t,ij} = \A_{e,r, i\cdot}' \F_{e,t} \A_{e,c, j\cdot} + \Sigma_{\epsilon, ij} \epsilon_{t, ij}.
\]
With Assumption (E2), we have
\[
\text{cov}(E_{t,ij}, E_{t,kj}) = \A_{e,r, i\cdot}'
\A_{e,r, k\cdot} \|\A_{e,c, j\cdot}\|^2 + \Sigma_{\epsilon, ij}^2 \b{1}_{\{i=k\}} ,
\]
and hence using Lemma \ref{lemma: correlation_Et_Ft},
\begin{equation*}
\begin{split}
    &\hspace{5mm}
    \b{E}\Bigg[ \Big\| \sum_{t=1}^T \E_t\E_t' \Big\|_F^2 \Bigg]
    = \sum_{i=1}^p\sum_{k=1}^p \b{E}\Bigg[ \Big(\sum_{t=1}^T \sum_{j=1}^q E_{t,ij} E_{t,kj}\Big)^2 \Bigg] \\
    &= \sum_{i=1}^p\sum_{k=1}^p \Bigg[\sum_{t=1}^T \sum_{j=1}^q \sum_{s=1}^T \sum_{l=1}^q \text{cov}(E_{t,ij} E_{t,kj}, E_{s,il} E_{s,kl}) + \Big( \sum_{t=1}^T \sum_{j=1}^q \b{E}[E_{t,ij} E_{t,kj}] \Big)^2 \Bigg] \\
    &=
    O(Tp^2q) + \sum_{i=1}^p\sum_{k=1}^p  O\Big( T\cdot \A_{e,r, i\cdot}' \A_{e,r, k\cdot} \|\A_{e,c}\|_F^2 + Tq\cdot\b{1}_{\{i=k\}} \Big)^2
    = O(T p^2 q + T^2 p q^2).
\end{split}
\end{equation*}

For (\ref{eqn: E11_bound}), consider first
\begin{equation}
\label{eqn: cov_EEEE}
\begin{split}
    &\hspace{5mm}
    \sum_{t=1}^T\sum_{s=1}^T \text{cov}(E_{t,ij}E_{t,kh}, E_{s,lm}E_{s,kn}) \\
    &=
    \sum_{t=1}^T\sum_{s=1}^T \text{cov}\Big[
    (\A_{e,r, i\cdot}' \F_{e,t} \A_{e,c, j\cdot} + \Sigma_{\epsilon, ij} \epsilon_{t, ij})
    (\A_{e,r, k\cdot}' \F_{e,t} \A_{e,c, h\cdot} + \Sigma_{\epsilon, kh} \epsilon_{t, kh}),\\
    &
    (\A_{e,r, l\cdot}' \F_{e,s} \A_{e,c, m\cdot} + \Sigma_{\epsilon, lm} \epsilon_{s, lm})
    (\A_{e,r, k\cdot}' \F_{e,s} \A_{e,c, n\cdot} + \Sigma_{\epsilon, kn} \epsilon_{s, kn})
    \Big] \\
    &=
    \sum_{t=1}^T\sum_{s=1}^T \text{cov}\Big[
    \A_{e,r, i\cdot}' \F_{e,t} \A_{e,c, j\cdot}
    \A_{e,r, k\cdot}' \F_{e,t} \A_{e,c, h\cdot},
    \A_{e,r, l\cdot}' \F_{e,s} \A_{e,c, m\cdot}
    \A_{e,r, k\cdot}' \F_{e,s} \A_{e,c, n\cdot}
    \Big] \\
    &+
    \sum_{t=1}^T\sum_{s=1}^T \b{E}\Big[\A_{e,r, i\cdot}' \F_{e,t} \A_{e,c, j\cdot} \A_{e,r, l\cdot}' \F_{e,s} \A_{e,c, m\cdot}\Big] \cdot \b{E}\Big[
    \Sigma_{\epsilon, kh} \epsilon_{t, kh}
    \Sigma_{\epsilon, kn} \epsilon_{s, kn}
    \Big] \\
    &+
    \sum_{t=1}^T\sum_{s=1}^T \b{E}\Big[\A_{e,r, i\cdot}' \F_{e,t} \A_{e,c, j\cdot} \A_{e,r, k\cdot}' \F_{e,s} \A_{e,c, n\cdot}\Big] \cdot \b{E}\Big[
    \Sigma_{\epsilon, kh} \epsilon_{t, kh}
    \Sigma_{\epsilon, lm} \epsilon_{s, lm}
    \Big] \\
    &+
    \sum_{t=1}^T\sum_{s=1}^T \b{E}\Big[\A_{e,r, k\cdot}' \F_{e,t} \A_{e,c, h\cdot} \A_{e,r, l\cdot}' \F_{e,s} \A_{e,c, m\cdot}\Big] \cdot \b{E}\Big[
    \Sigma_{\epsilon, ij} \epsilon_{t, ij}
    \Sigma_{\epsilon, kn} \epsilon_{s, kn}
    \Big] \\
    &+
    \sum_{t=1}^T\sum_{s=1}^T \b{E}\Big[\A_{e,r, k\cdot}' \F_{e,t} \A_{e,c, h\cdot} \A_{e,r, k\cdot}' \F_{e,s} \A_{e,c, n\cdot}\Big] \cdot \b{E}\Big[
    \Sigma_{\epsilon, ij} \epsilon_{t, ij}
    \Sigma_{\epsilon, lm} \epsilon_{s, lm}
    \Big] \\
    &+
    \sum_{t=1}^T\sum_{s=1}^T \text{cov}\Big[\Sigma_{\epsilon, ij} \epsilon_{t, ij} \Sigma_{\epsilon, kh} \epsilon_{t, kh}, \Sigma_{\epsilon, lm} \epsilon_{s, lm} \Sigma_{\epsilon, kn} \epsilon_{s, kn} \Big] \\
    &=
    \sum_{t=1}^T\sum_{s=1}^T \text{cov}\Big[
    \A_{e,r, i\cdot}' \Big(\sum_{w\geq 0}a_{e,w}\X_{e,t-w}\Big) \A_{e,c, j\cdot}
    \A_{e,r, k\cdot}' \Big(\sum_{w\geq 0}a_{e,w}\X_{e,t-w}\Big) \A_{e,c, h\cdot}, \\
    &
    \A_{e,r, l\cdot}' \Big(\sum_{w\geq 0}a_{e,w}\X_{e,s-w}\Big) \A_{e,c, m\cdot}
    \A_{e,r, k\cdot}' \Big(\sum_{w\geq 0}a_{e,w}\X_{e,s-w}\Big) \A_{e,c, n\cdot}
    \Big] \\
    &+
    \sum_{t=1}^T\sum_{s=1}^T \b{E}\Big[\A_{e,r, i\cdot}' \Big(\sum_{w\geq 0}a_{e,w}\X_{e,t-w}\Big) \A_{e,c, j\cdot} \A_{e,r, l\cdot}' \Big(\sum_{w\geq 0}a_{e,w}\X_{e,s-w}\Big) \A_{e,c, m\cdot}\Big] \cdot \b{E}\Big[
    \Sigma_{\epsilon, kh} \epsilon_{t, kh}
    \Sigma_{\epsilon, kn} \epsilon_{s, kn}
    \Big] \\
    &+
    \sum_{t=1}^T\sum_{s=1}^T \b{E}\Big[\A_{e,r, i\cdot}' \Big(\sum_{w\geq 0}a_{e,w}\X_{e,t-w}\Big) \A_{e,c, j\cdot} \A_{e,r, k\cdot}' \Big(\sum_{w\geq 0}a_{e,w}\X_{e,s-w}\Big) \A_{e,c, n\cdot}\Big] \cdot \b{E}\Big[
    \Sigma_{\epsilon, kh} \epsilon_{t, kh}
    \Sigma_{\epsilon, lm} \epsilon_{s, lm}
    \Big] \\
    &+
    \sum_{t=1}^T\sum_{s=1}^T \b{E}\Big[\A_{e,r, k\cdot}' \Big(\sum_{w\geq 0}a_{e,w}\X_{e,t-w}\Big) \A_{e,c, h\cdot} \A_{e,r, l\cdot}' \Big(\sum_{w\geq 0}a_{e,w}\X_{e,s-w}\Big) \A_{e,c, m\cdot}\Big] \cdot \b{E}\Big[
    \Sigma_{\epsilon, ij} \epsilon_{t, ij}
    \Sigma_{\epsilon, kn} \epsilon_{s, kn}
    \Big] \\
    &+
    \sum_{t=1}^T\sum_{s=1}^T \b{E}\Big[\A_{e,r, k\cdot}' \Big(\sum_{w\geq 0}a_{e,w}\X_{e,t-w}\Big) \A_{e,c, h\cdot} \A_{e,r, k\cdot}' \Big(\sum_{w\geq 0}a_{e,w}\X_{e,s-w}\Big) \A_{e,c, n\cdot}\Big] \cdot \b{E}\Big[
    \Sigma_{\epsilon, ij} \epsilon_{t, ij}
    \Sigma_{\epsilon, lm} \epsilon_{s, lm}
    \Big] \\
    &+
    \sum_{t=1}^T\sum_{s=1}^T \text{cov}\Big[\Sigma_{\epsilon, ij} \epsilon_{t, ij} \Sigma_{\epsilon, kh} \epsilon_{t, kh}, \Sigma_{\epsilon, lm} \epsilon_{s, lm} \Sigma_{\epsilon, kn} \epsilon_{s, kn} \Big] .
\end{split}
\end{equation}
Consider the six terms in the last equality above, we have the first term as
\begin{equation}
\label{eqn: cov_EEEE_1}
\begin{split}
    &\hspace{5mm}
    \sum_{t=1}^T\sum_{s=1}^T \text{cov}\Big[
    \A_{e,r, i\cdot}' \Big(\sum_{w\geq 0}a_{e,w}\X_{e,t-w}\Big) \A_{e,c, j\cdot}
    \A_{e,r, k\cdot}' \Big(\sum_{w\geq 0}a_{e,w}\X_{e,t-w}\Big) \A_{e,c, h\cdot}, \\
    &
    \A_{e,r, l\cdot}' \Big(\sum_{w\geq 0}a_{e,w}\X_{e,s-w}\Big) \A_{e,c, m\cdot}
    \A_{e,r, k\cdot}' \Big(\sum_{w\geq 0}a_{e,w}\X_{e,s-w}\Big) \A_{e,c, n\cdot}
    \Big] \\
    &=
    \sum_{t=1}^T \sum_{w\geq 0}a_{e,w}^2 \cdot \b{E}\Big[
    \A_{e,r, i\cdot}' \X_{e,t-w} \A_{e,c, j\cdot}
    \A_{e,r, k\cdot}' \X_{e,t-w} \A_{e,c, n\cdot}
    \Big] \cdot \b{E}\Big[
    \A_{e,r, k\cdot}' \X_{e,t-w} \A_{e,c, h\cdot}
    \A_{e,r, l\cdot}' \X_{e,t-w} \A_{e,c, m\cdot}
    \Big] \\
    &+
    \sum_{t=1}^T \sum_{w\geq 0}a_{e,w}^2 \cdot \b{E}\Big[
    \A_{e,r, i\cdot}' \X_{e,t-w} \A_{e,c, j\cdot}
    \A_{e,r, l\cdot}' \X_{e,t-w} \A_{e,c, m\cdot}
    \Big] \cdot \b{E}\Big[
    \A_{e,r, k\cdot}' \X_{e,t-w} \A_{e,c, h\cdot}
    \A_{e,r, k\cdot}' \X_{e,t-w} \A_{e,c, n\cdot}
    \Big] \\
    &+
    \sum_{t=1}^T \sum_{w\geq 0} a_{e,w}^4
    \cdot \b{E}\Big[
    \A_{e,r, i\cdot}' \X_{e,t-w} \A_{e,c, j\cdot}
    \A_{e,r, k\cdot}' \X_{e,t-w} \A_{e,c, h\cdot}
    \A_{e,r, l\cdot}' \X_{e,t-w} \A_{e,c, m\cdot}
    \A_{e,r, k\cdot}' \X_{e,t-w} \A_{e,c, n\cdot}
    \Big] \\
    &-
    \sum_{t=1}^T \sum_{w\geq 0} a_{e,w}^4
    \cdot \b{E}\Big[
    \A_{e,r, i\cdot}' \X_{e,t-w} \A_{e,c, j\cdot}
    \A_{e,r, k\cdot}' \X_{e,t-w} \A_{e,c, h\cdot}
    \Big] \cdot \b{E}\Big[
    \A_{e,r, l\cdot}' \X_{e,t-w} \A_{e,c, m\cdot}
    \A_{e,r, k\cdot}' \X_{e,t-w} \A_{e,c, n\cdot}
    \Big] \\
    &=
    \sum_{t=1}^T
    \A_{e,r, k\cdot}' \A_{e,c, j\cdot}
    \A_{e,r, i\cdot}' \A_{e,c, n\cdot}
    \A_{e,r, l\cdot}' \A_{e,c, h\cdot}
    \A_{e,r, k\cdot}' \A_{e,c, m\cdot} \\
    &+
    \sum_{t=1}^T
    \A_{e,r, l\cdot}' \A_{e,c, j\cdot}
    \A_{e,r, i\cdot}' \A_{e,c, m\cdot}
    \A_{e,r, k\cdot}' \A_{e,c, h\cdot}
    \A_{e,r, k\cdot}' \A_{e,c, n\cdot} \\
    &+
    O(1)\cdot \sum_{t=1}^T \sum_{w\geq 0} a_{e,w}^4
    \cdot \|\A_{e,r, k\cdot}\|^2
    \cdot \|\A_{e,r, i\cdot}\|
    \cdot \|\A_{e,c, j\cdot}\|
    \cdot \|\A_{e,c, h\cdot}\|
    \cdot \|\A_{e,r, l\cdot}\|
    \cdot \|\A_{e,c, m\cdot}\|
    \cdot \|\A_{e,c, n\cdot}\| \\
    &-
    \sum_{t=1}^T \sum_{w\geq 0} a_{e,w}^4 \cdot
    \A_{e,r, k\cdot}' \A_{e,c, j\cdot}
    \A_{e,r, i\cdot}' \A_{e,c, h\cdot}
    \A_{e,r, k\cdot}' \A_{e,c, m\cdot}
    \A_{e,r, l\cdot}' \A_{e,c, n\cdot} ,
\end{split}
\end{equation}
where we used (E2) in the last equality that each entry in $\{\X_{e,t}\}$ is independent with uniformly bounded fourth moment. Similarly, the remaining terms in last equality of (\ref{eqn: cov_EEEE}) are
\begin{align}
    &\hspace{5mm}
    \sum_{t=1}^T\sum_{s=1}^T \b{E}\Big[\A_{e,r, i\cdot}' \Big(\sum_{w\geq 0}a_{e,w}\X_{e,t-w}\Big) \A_{e,c, j\cdot} \A_{e,r, l\cdot}' \Big(\sum_{w\geq 0}a_{e,w}\X_{e,s-w}\Big) \A_{e,c, m\cdot}\Big] \cdot \b{E}\Big[
    \Sigma_{\epsilon, kh} \epsilon_{t, kh}
    \Sigma_{\epsilon, kn} \epsilon_{s, kn}
    \Big]
    \notag \\
    &=
    \sum_{t=1}^T \A_{e,r, l\cdot}' \A_{e,c, j\cdot}
    \A_{e,r, i\cdot}' \A_{e,c, m\cdot}
    \cdot \Sigma_{\epsilon, kh} \Sigma_{\epsilon, kn} \cdot \b{1}_{\{h=n\}} ,
    \label{eqn: cov_EEEE_2} \\
    &\hspace{5mm}
    \sum_{t=1}^T\sum_{s=1}^T \b{E}\Big[\A_{e,r, i\cdot}' \Big(\sum_{w\geq 0}a_{e,w}\X_{e,t-w}\Big) \A_{e,c, j\cdot} \A_{e,r, k\cdot}' \Big(\sum_{w\geq 0}a_{e,w}\X_{e,s-w}\Big) \A_{e,c, n\cdot}\Big] \cdot \b{E}\Big[
    \Sigma_{\epsilon, kh} \epsilon_{t, kh}
    \Sigma_{\epsilon, lm} \epsilon_{s, lm}
    \Big]
    \notag \\
    &=
    \sum_{t=1}^T \A_{e,r, k\cdot}' \A_{e,c, j\cdot}
    \A_{e,r, i\cdot}' \A_{e,c, n\cdot}
    \cdot \Sigma_{\epsilon, kh} \Sigma_{\epsilon, lm} \cdot \b{1}_{\{k=l\}} \b{1}_{\{h=m\}} ,
    \label{eqn: cov_EEEE_3} \\
    &\hspace{5mm}
    \sum_{t=1}^T\sum_{s=1}^T \b{E}\Big[\A_{e,r, k\cdot}' \Big(\sum_{w\geq 0}a_{e,w}\X_{e,t-w}\Big) \A_{e,c, h\cdot} \A_{e,r, l\cdot}' \Big(\sum_{w\geq 0}a_{e,w}\X_{e,s-w}\Big) \A_{e,c, m\cdot}\Big] \cdot \b{E}\Big[
    \Sigma_{\epsilon, ij} \epsilon_{t, ij}
    \Sigma_{\epsilon, kn} \epsilon_{s, kn}
    \Big]
    \notag \\
    &=
    \sum_{t=1}^T \A_{e,r, l\cdot}' \A_{e,c, h\cdot}
    \A_{e,r, k\cdot}' \A_{e,c, m\cdot}
    \cdot \Sigma_{\epsilon, ij} \Sigma_{\epsilon, kn} \cdot \b{1}_{\{i=k\}} \b{1}_{\{j=n\}} ,
    \label{eqn: cov_EEEE_4} \\
    &\hspace{5mm}
    \sum_{t=1}^T\sum_{s=1}^T \b{E}\Big[\A_{e,r, k\cdot}' \Big(\sum_{w\geq 0}a_{e,w}\X_{e,t-w}\Big) \A_{e,c, h\cdot} \A_{e,r, k\cdot}' \Big(\sum_{w\geq 0}a_{e,w}\X_{e,s-w}\Big) \A_{e,c, n\cdot}\Big] \cdot \b{E}\Big[
    \Sigma_{\epsilon, ij} \epsilon_{t, ij}
    \Sigma_{\epsilon, lm} \epsilon_{s, lm}
    \Big]
    \notag \\
    &=
    \sum_{t=1}^T \A_{e,r, k\cdot}' \A_{e,c, h\cdot}
    \A_{e,r, k\cdot}' \A_{e,c, n\cdot}
    \cdot \Sigma_{\epsilon, ij} \Sigma_{\epsilon, lm} \cdot \b{1}_{\{i=l\}} \b{1}_{\{j=m\}} ,
    \label{eqn: cov_EEEE_5} \\
    &\hspace{5mm}
    \sum_{t=1}^T\sum_{s=1}^T \text{cov}\Big[\Sigma_{\epsilon, ij} \epsilon_{t, ij} \Sigma_{\epsilon, kh} \epsilon_{t, kh}, \Sigma_{\epsilon, lm} \epsilon_{s, lm} \Sigma_{\epsilon, kn} \epsilon_{s, kn} \Big]
    \notag \\
    &=
    O(T)\cdot \Big(
    \b{1}_{\{i=l=k\}} \b{1}_{\{j=h=m=n\}} +
    \b{1}_{\{i=l\}} \b{1}_{\{j=m\}} \b{1}_{\{h=n\}} + \b{1}_{\{i=l=k\}} \b{1}_{\{j=n\}} \b{1}_{\{h=m\}}
    \Big) .
    \label{eqn: cov_EEEE_6}
\end{align}
Using (\ref{eqn: cov_EEEE_1}), (\ref{eqn: cov_EEEE_2}), (\ref{eqn: cov_EEEE_3}), (\ref{eqn: cov_EEEE_4}), (\ref{eqn: cov_EEEE_5}) and (\ref{eqn: cov_EEEE_6}), we arrive at an expression for (\ref{eqn: cov_EEEE}). Thus, (\ref{eqn: E11_bound}) can be obtained as
\begin{equation}
\label{eqn: E11_bound_last_step}
\begin{split}
    &\hspace{5mm}
    \b{E}\Bigg[ \Big\| \sum_{t=1}^T \1_q' \E_t' \1_p \E_t \1_q \1_p' \Big\|_F^2 \Bigg] =
    p \sum_{k=1}^p \b{E} \Bigg\{ \Bigg[\sum_{t=1}^T \Big( \sum_{i=1}^p \sum_{j=1}^q E_{t,ij} \Big) \sum_{h=1}^q E_{t,kh} \Bigg]^2\Bigg\} \\
    &= p\sum_{k=1}^p \Bigg[
    \sum_{t=1}^T \sum_{i=1}^p \sum_{j=1}^q \sum_{h=1}^q \sum_{s=1}^T \sum_{l=1}^p \sum_{m=1}^q \sum_{n=1}^q
    \text{cov}(E_{t,ij}E_{t,kh}, E_{s,lm}E_{s,kn})
    +
    \Big(\sum_{t=1}^T \sum_{i=1}^p \sum_{j=1}^q \sum_{h=1}^q \b{E}[E_{t,ij}E_{t,kh}] \Big)^2
    \Bigg] \\
    &=
    O(T p^3 q^2) + p\sum_{k=1}^p
    \Big[\sum_{t=1}^T \sum_{i=1}^p \sum_{j=1}^q \sum_{h=1}^q (\A_{e,r, i\cdot}' \A_{e,r, k\cdot} \A_{e,c, j\cdot}' \A_{e,c, h\cdot} + \Sigma_{\epsilon, ij}^2 \b{1}_{\{i=k\}} \b{1}_{\{j=h\}} ) \Big]^2 \\
    &=
    O(T p^3 q^2 + T^2 p^2 q^2) .
\end{split}
\end{equation}

By (\ref{eqn: CE_bound}) and (\ref{eqn: EE_bound}), we can obtain (\ref{eqn: CE11_bound}) and (\ref{eqn: EE11_bound}), respectively as follows,
\begin{align*}
    & \Big\| \sum_{t=1}^T \C_t \E_t' \1_p \1_p' \Big\|_F^2 \leq
    \Big\| \sum_{t=1}^T \C_t \E_t'\Big\|_F^2
    \cdot \| \1_p \1_p' \|_F^2
    = O_P(T p^{3+\delta_{r,1}} q), \\
    & \Big\| \sum_{t=1}^T \E_t \E_t' \1_p \1_p' \Big\|_F^2 \leq
    \Big\| \sum_{t=1}^T \E_t \E_t'\Big\|_F^2
    \cdot \| \1_p \1_p' \|_F^2
    = O_P(T p^4 q + T^2 p^3 q^2) .
\end{align*}

Similar to the proof of (\ref{eqn: E11_bound}), we can show (\ref{eqn: E11E_bound}) by
\begin{equation*}
\begin{split}
    &\hspace{5mm}
    \b{E}\Bigg[ \Big\| \sum_{t=1}^T \E_t\1_q\1_q' \E_t' \Big\|_F^2 \Bigg] =
    \sum_{i=1}^p \sum_{k=1}^p \b{E}\Bigg\{ \Big[\sum_{t=1}^T \Big(\sum_{j=1}^q E_{t,ij} \Big) \Big(\sum_{h=1}^q E_{t,kh} \Big)\Big]^2 \Bigg\} \\
    &=
    \sum_{i=1}^p \sum_{k=1}^p \Bigg[
    \sum_{t=1}^T \sum_{j=1}^q \sum_{h=1}^q
    \sum_{s=1}^T \sum_{m=1}^q \sum_{n=1}^q
    \text{cov}(E_{t,ij} E_{t,kh}, E_{s,im} E_{s,kn}) + \Big(\sum_{t=1}^T \sum_{j=1}^q \sum_{h=1}^q \b{E}[E_{t,ij} E_{t,kh}]\Big)^2
    \Bigg] \\
    &=
    O(T p^2 q^2) + \sum_{i=1}^p \sum_{k=1}^p
    \Big[\sum_{t=1}^T \sum_{j=1}^q \sum_{h=1}^q (\A_{e,r, i\cdot}' \A_{e,r, k\cdot} \A_{e,c, j\cdot}' \A_{e,c, h\cdot} + \Sigma_{\epsilon, ij}^2 \b{1}_{\{i=k\}} \b{1}_{\{j=h\}} ) \Big]^2 \\
    &= O(T p^2 q^2 + T^2 p q^2) .
\end{split}
\end{equation*}

From (\ref{eqn: 1E1_bound}), we can obtain (\ref{eqn: 1E111_bound}) such that
\[
\Big\| \sum_{t=1}^T (\1_q' \E_t' \1_p)^2 \1_p \1_p' \Big\|_F^2 =
\Big[\sum_{t=1}^T (\1_q' \E_t' \1_p)^2 \Big]^2 \cdot \Big\| \1_p \1_p' \Big\|_F^2 =
O_P( T^2 p^4 q^2 ).
\]
Lastly, from (\ref{eqn: EE_bound}) we have
\[
\Big\| \sum_{t=1}^T \1_p' \E_t \E_t' \1_p \1_p \1_p' \Big\|_F^2 \leq
\|\1_p\|^2 \cdot \|\1_p\|^2 \cdot
\Big\| \sum_{t=1}^T \E_t \E_t' \Big\|_F^2 \cdot
\| \1_p \1_p' \|_F^2 =
O_P(T p^6 q + T^2 p^5 q^2) .
\]
From (\ref{eqn: R_rt}), we have
\begin{equation*}
\begin{split}
    \Big\| \sum_{t=1}^T \R_{r,t} \Big\|_F^2 &=
    O_P\Bigg(
    \Big\| \sum_{t=1}^T \C_t\E_t' \Big\|_F^2 +
    \Big\| \sum_{t=1}^T \E_t\E_t' \Big\|_F^2 +
    (pq)^{-2}\Big\| \sum_{t=1}^T \1_q' \E_t' \1_p \E_t \1_q \1_p' \Big\|_F^2
    \notag \\
     & + p^{-2}\Big\| \sum_{t=1}^T \C_t\E_t'\1_p\1_p' \Big\|_F^2
     + p^{-2}\Big\|\sum_{t=1}^T \E_t\E_t'\1_p\1_p' \Big\|_F^2 + p^{-2}\Big\|\sum_{t=1}^T \E_t\1_q\1_q'\E_t' \Big\|_F^2 \notag \\
     &+ (pq)^{-2}p^{-2} \Big\|\sum_{t=1}^T (\1_q'\E_t'\1_p)^2\1_p\1_p'\Big\|_F^2 + p^{-4} \Big\| \sum_{t=1}^T \1_p' \E_t \E_t' \1_p \1_p \1_p' \Big\|_F^2 \Bigg) \\
     &=
     O_P(Tp^2q + T^2pq^2) ,
\end{split}
\end{equation*}
which completes the proof of Lemma \ref{lemma: rate_R_rt}. $\square$

\newpage

\begin{lemma}\label{lemma: norm_hat_D}
Let Assumptions (M1), (F1), (L1), (E1), (E2) and (R1) hold. Define $\omega_r:= p^{\delta_{r,k_r}} q^{\delta_{c,1}}$ and $\omega_c:= q^{\delta_{c,k_c}} p^{\delta_{r,1}}$. We have
\begin{equation*}
     \|\wh\D_r^{-1} \|_F =
     O_P(\omega_r^{-1}) ,\;\;\;
     \|\wh\D_c^{-1} \|_F =
     O_P(\omega_c^{-1}) .
\end{equation*}
\end{lemma}
\noindent\textbf{\textit{Proof of Lemma \ref{lemma: norm_hat_D}.}}
It suffices to show the bound of $\|\wh\D_r^{-1}\|_F$, since that of $\|\wh\D_c^{-1}\|_F$ would be similar. First, we bound the term $\|\wh\D_r^{-1}\|_F^2$ by finding the lower bound of $\lambda_{k_r}(\wh\D_r)$. To do this, consider the decomposition
\begin{align}
\label{eqn: CC_decomposition}
  \frac{1}{T}\sum_{t=1}^T \wh\L_t\wh\L_t' = \frac{1}{T}\sum_{t=1}^T\C_t\C_t' + \frac{1}{T}\sum_{t=1}^T \R_{r,t},
\end{align}
so that for a unit vector $\bgamma \in\b{R}^p$, we can define
\begin{align}
    S_r(\bgamma) &:=
    \frac{1}{\omega_r}\bm{\gamma}' \Big(\frac{1}{T}\sum_{t=1}^T \wh\L_t\wh\L_t'\Big) \bm{\gamma}
    =: S_r^*(\bgamma) + \wt{S}_r(\bgamma), \; \text{ with} \notag\\
    S_r^*(\bgamma) &:=
    \frac{1}{\omega_r}\bm{\gamma}' \Big(\frac{1}{T}\sum_{t=1}^T \C_t\C_t' \Big) \bm{\gamma}
    , \;\;\;
    \wt{S}_r(\bgamma) :=
    \frac{1}{\omega_r}\bm{\gamma}' \Big(\frac{1}{T}\sum_{t=1}^T \R_{r,t} \Big) \bm{\gamma}
    . \notag
\end{align}
Since $\|\bgamma\|=1$, we have by Lemma \ref{lemma: rate_R_rt},
\[
|\wt{S}_r(\bgamma)|^2 \leq
\frac{1}{\omega_r^2 T^2} \Big\|\sum_{t=1}^T \R_{r,t} \Big\|_F^2 =
O_P\Big( T^{-1}p^{2(1-\delta_{r,k_r})}q^{1-2\delta_{c,1}} + p^{1-2\delta_{r,k_r}}q^{2(1-\delta_{c,1})} \Big) = o_P(1),
\]
where the last equality used Assumption (R1). Next, with Assumption (F1), consider
\begin{align*}
    \lambda_{k_r} \Big(\frac{1}{T}\sum_{t=1}^T \C_t\C_t'\Big) &= \lambda_{k_r}\Big(
    \frac{1}{T}\sum_{t=1}^T \A_r \F_t \A_c' \A_c \F_t' \A_r' \Big)
    \geq \lambda_{k_r}(\A_r'\A_r)\cdot \lambda_{k_r}\Big(\frac{1}{T}\sum_{t=1}^T
    \F_t \A_c'\A_c\F_t' \Big)\\
    &\asymp_P p^{\delta_{r,k_r}} \cdot \lambda_{k_r}(\tr(\A_c'\A_c)\bSigma_r) \asymp_P
    p^{\delta_{r,k_r}} q^{\delta_{c,1}} = \omega_r.
\end{align*}
With this, going back to the decomposition (\ref{eqn: CC_decomposition}),
\begin{align*}
  \omega_r^{-1}\lambda_{k_r}(\wh\D_r) &=
  \omega_r^{-1}\lambda_{k_r} \Big(\frac{1}{T} \sum_{t=1}^T \wh\L_t\wh\L_t'\Big)
  \geq
  \omega_r^{-1}\lambda_{k_r} \Big(\frac{1}{T} \sum_{t=1}^T \C_t\C_t'\Big)
  - \sup_{\|\bgamma\|=1} |\wt{S}_r(\bgamma)|\asymp_P 1,
\end{align*}
and hence finally $\norm{\wh\D_r^{-1}}_F = O_P\big( \lambda_{k_r}^{-1}(\wh\D_r) \big) = O_P(\omega_r^{-1})$, which completes the proof of Lemma \ref{lemma: norm_hat_D}. $\square$

\newpage

\begin{lemma}\label{lemma: limit_D}
(Limit of $\wh\D_r$ and $\wh\D_c$) Let Assumptions (F1), (L1), (E1), (E2) and (R1) hold. With $\wh\D_r, \wh\D_c$ and $\omega_r, \omega_c$ defined in Lemma \ref{lemma: norm_hat_D}, we have
\begin{align*}
    & \omega_r^{-1}\wh\D_r \xrightarrow{p} \omega_r^{-1}\D_r := \omega_r^{-1}\textnormal{tr}(\A_c'\A_c)\cdot \textnormal{diag}\{\lambda_j(\A_r'\A_r) \mid j\in[k_r]\} ,\\
    & \omega_c^{-1}\wh\D_c \xrightarrow{p} \omega_c^{-1}\D_c := \omega_c^{-1}\textnormal{tr}(\A_r'\A_r)\cdot \textnormal{diag}\{\lambda_j(\A_c'\A_c) \mid j\in[k_c]\} .
\end{align*}
\end{lemma}
\textbf{\textit{Proof of Lemma \ref{lemma: limit_D}.}}
It suffices to show the limit of $\wh\D_r$, as the proof for $\wh\D_c$ will be similar. We first show
\begin{equation}
\label{eqn: limit_CC}
\frac{1}{\omega_r T} \sum_{t=1}^T \C_t\C_t'
\xrightarrow{p} \tr(\A_c'\A_c)\cdot \omega_r^{-1} \A_r\A_r' ,
\end{equation}
where $\lambda_{k_r}(\tr(\A_c'\A_c)\cdot \omega_r^{-1} \A_r\A_r') \asymp 1$. By Assumption (F1), we have
\[
\b{E}\Big( \frac{1}{\omega_r T} \sum_{t=1}^T \C_t\C_t' \Big) = \frac{1}{\omega_r T} \sum_{t=1}^T \b{E}[\A_r\F_t\A_c' \A_c\F_t'\A_r'] = \tr(\A_c'\A_c)\cdot \omega_r^{-1} \A_r\A_r',
\]
where we used the independence structure among elements in $\F_t$. Furthermore, for any $i,j\in [p]$,
\begin{equation*}
\begin{split}
    \text{Var}\Big[ \Big(\frac{1}{\omega_r T} \sum_{t=1}^T \C_t\C_t' \Big)_{ij}\Big] &=
    \frac{1}{\omega_r^2 T^2}\sum_{t=1}^T
    \sum_{s=1}^T\text{cov}\bigg(
    \A_{r,i\cdot}'\Big(\sum_{w\geq 0}a_{f,w}\X_{f,t-w}\Big)\A_c'
    \A_c\Big(\sum_{w\geq 0}a_{f,w}\X_{f,t-w}'\Big)
    \A_{r,j\cdot},\\
    &
    \A_{r,i\cdot}'\Big(\sum_{w\geq 0}a_{f,w}\X_{f,s-w}\Big)\A_c'
    \A_c\Big(\sum_{w\geq 0}a_{f,w}\X_{f,s-w}'\Big)
    \A_{r,j\cdot}\bigg) \\
    &=
    \frac{1}{\omega_r^2 T^2}
    \sum_{t=1}^T\sum_{w\geq 0}
    \sum_{l\geq 0}a_{f,w}^2
    a_{f,l}^2\cdot \text{Var}\bigg(
    \A_{r,i\cdot}'\X_{f,t-w}\A_c'
    \A_c\X_{f,t-l}' \A_{r,j\cdot}\bigg)\\
    &=
    \frac{1}{\omega_r^2 T^2}
    \sum_{t=1}^T\sum_{w\geq 0}
    \sum_{l\geq 0}a_{f,w}^2
    a_{f,l}^2\cdot O(\|\A_c\|_F^4)
    =
    O( T^{-1} p^{-2\delta_{r,k_r}}) = o(1),
\end{split}
\end{equation*}
where we used Assumption (F1) in the third last equality, and both (L1) and (F1) in the second last, which concludes (\ref{eqn: limit_CC}). Then it holds that
\begin{equation*}
\begin{split}
    &\hspace{5mm}
    \Big\|\frac{1}{\omega_rT} \sum_{t=1}^T  \wh\L_t\wh\L_t' - \tr(\A_c'\A_c)\cdot \omega_r^{-1} \A_r\A_r' \Big\|_F^2 \\
    &\leq
    2 \cdot \Big\|\frac{1}{\omega_rT} \sum_{t=1}^T  \wh\L_t\wh\L_t' - \frac{1}{\omega_r T} \sum_{t=1}^T \C_t\C_t' \Big\|_F^2
    + 2 \cdot \Big\|\frac{1}{\omega_r T} \sum_{t=1}^T \C_t\C_t' - \tr(\A_c'\A_c)\cdot \omega_r^{-1} \A_r\A_r' \Big\|_F^2 \\
    &=
    2 \cdot \Big\|\frac{1}{\omega_r T} \sum_{t=1}^T \R_{r,t} \Big\|_F^2
    + 2 \cdot \Big\|\frac{1}{\omega_r T} \sum_{t=1}^T \C_t\C_t' - \tr(\A_c'\A_c)\cdot \omega_r^{-1} \A_r\A_r' \Big\|_F^2 \\
    &=
    O_P\Big(
    T^{-1}p^{2(1-\delta_{r,k_r})}q^{1-2\delta_{c,1}} + p^{1-2\delta_{r,k_r}}q^{2(1-\delta_{c,1})} \Big) + o_P(1) = o_P(1) ,
\end{split}
\end{equation*}
where the second last equality used Lemma \ref{lemma: rate_R_rt} and (\ref{eqn: limit_CC}), and the last used Assumption (R1). Using the inequality that for the $i$-th eigenvalue of matrices $\wh\A$ and $\A$, $|\lambda_i(\wh\A) - \lambda_i(\A)|\leq \|\wh\A - \A\|\leq \|\wh\A - \A\|_F$, we have for any $i\in[k_r]$,
\[
|(\omega_r^{-1}\wh\D_r)_{ii} - (\omega_r^{-1}\D_r)_{ii}| = o_P(1).
\]
Thus, $\omega_r^{-1}\wh\D_r \xrightarrow{p} \omega_r^{-1}\D_r$. This completes the proof of Lemma \ref{lemma: limit_D}. $\square$

\begin{lemma}\label{lemma: limit_H}
(Limit of $\H_r$ and $\H_c$) Under Assumptions (F1), (L1), (E1), (E2), (R1) and (L2), we have
\begin{align*}
    & \H_r \xrightarrow{p} \H_r^\ast := (\textnormal{tr}(\A_c'\A_c))^{1/2} \cdot \D_r^{-1/2} (\bGamma_r^\ast)' \Z_r^{1/2} ,\\
    & \H_c \xrightarrow{p} \H_c^\ast := (\textnormal{tr}(\A_r'\A_r))^{1/2} \cdot \D_c^{-1/2} (\bGamma_c^\ast)' \Z_c^{1/2} ,
\end{align*}
where $\D_r, \D_c$ are defined in Lemma \ref{lemma: limit_D}, and $\bGamma_r^\ast, \bGamma_c^\ast$ are the eigenvector matrices of $\textnormal{tr}(\A_c'\A_c)\cdot \omega_r^{-1} \Z_r^{1/2} \bSigma_{A,r} \Z_r^{1/2}$ and $\textnormal{tr}(\A_r'\A_r)\cdot \omega_c^{-1} \Z_c^{1/2} \bSigma_{A,c} \Z_c^{1/2}$, respectively.
\end{lemma}
\textbf{\textit{Proof of Lemma \ref{lemma: limit_H}.}}
The proof of the two limits are similar, and hence we only show the probability limit of $\H_r$ is $\H_r^\ast$. First, left-multiply $\omega_r^{-1}\Z_r^{1/2}\Q_r'$ on (\ref{eqn: row_svd}), we have
\begin{equation*}
\begin{split}
    (\Z_r^{1/2}\Q_r'\wh\Q_r)(\omega_r^{-1}\wh\D_r) &= \omega_r^{-1}\Z_r^{1/2}\Q_r' \Big(T^{-1}\sum_{t=1}^T \wh\L_t\wh\L_t'\Big) \wh\Q_r \\
    &=
    \Bigg[\frac{1}{T}\sum_{t=1}^T \omega_r^{-1}\Z_r^{1/2}\Q_r'\Q_r\F_{Z,t}\Q_c'\Q_c\F_{Z,t}' \Z_r^{-1/2} + \R_{r,\text{res}}
    \Bigg] (\Z_r^{1/2}\Q_r'\wh\Q_r) ,
\end{split}
\end{equation*}
where $\R_{r,\text{res}} := T^{-1}\sum_{t=1}^T \omega_r^{-1}\Z_r^{1/2}\Q_r' \R_{r,t} \wh\Q_r (\Z_r^{1/2} \Q_r' \wh\Q_r)^{-1}$. This implies each column of $(\Z_r^{1/2}\Q_r'\wh\Q_r)$ is an eigenvector of the matrix $\Big[T^{-1} \sum_{t=1}^T \omega_r^{-1} \Z_r^{1/2} \Q_r'\Q_r \F_{Z,t}\Q_c'\Q_c\F_{Z,t}' \Z_r^{-1/2} + \R_{r,\text{res}} \Big]$. Note that
\begin{equation*}
\begin{split}
    &\hspace{5mm}
    \omega_r^{-1}\Big[ (\Z_r^{1/2}\Q_r'\wh\Q_r)' (\Z_r^{1/2}\Q_r'\wh\Q_r) - \tr(\A_c'\A_c)^{-1} \cdot \D_r \Big] \\
    &=
    \omega_r^{-1}\wh\Q'\A_r'\A_r\wh\Q -
    \frac{1}{\tr(\A_c'\A_c)} \wh\Q'\Big( \frac{1}{\omega_r T}\sum_{t=1}^T \C_t\C_t' \Big)\wh\Q +
    \frac{1}{\tr(\A_c'\A_c)} \wh\Q'\Big( \frac{1}{\omega_r T}\sum_{t=1}^T \C_t\C_t' \Big)\wh\Q
    - \tr(\A_c'\A_c)^{-1} \cdot \D_r ,
\end{split}
\end{equation*}
whose Frobenius norm is $o_P(1)$ by (\ref{eqn: limit_CC}) and Lemma \ref{lemma: limit_D}. We arrive at $(\Z_r^{1/2} \Q_r' \wh\Q_r)' (\Z_r^{1/2} \Q_r' \wh\Q_r) \xrightarrow{p} \tr(\A_c'\A_c)^{-1} \cdot \D_r$. Thus, the eigenvalues of $(\Q_r'\wh\Q_r)'(\Q_r'\wh\Q_r)$ are asymptotically bounded away from zero and infinity by Assumption (L1), and hence $\|(\Z_r^{1/2} \Q_r' \wh\Q_r)^{-1}\|_F = O_P(\|\Z_r^{-1/2}\|_F)$. Therefore, we have
\begin{equation}
\label{eqn: rate_R_res}
\begin{split}
    \|\R_{r,\text{res}}\|_F^2 &= O_P(1) \cdot
    \Big\|\frac{1}{\omega_r T} \sum_{t=1}^T \R_{r,t} \Big\|_F^2 \cdot \|\Z_r^{1/2}\|_F^2
    \cdot \|(\Z_r^{1/2} \Q_r' \wh\Q_r)^{-1}\|_F^2\\
    &=
    O_P\Big( T^{-1}p^{2(1-\delta_{r,k_r})}q^{1-2\delta_{c,1}} + p^{1-2\delta_{r,k_r}}q^{2(1-\delta_{c,1})} \Big) = o_P(1),
\end{split}
\end{equation}
where we used Lemma \ref{lemma: rate_R_rt} in the second last equality and Assumption (R1) in the last.

Denote by the following as the normalization of $\Z_r^{1/2} \Q_r' \wh\Q_r$,
\[
\bGamma_r := (\tr(\A_c'\A_c))^{1/2}\cdot (\Z_r^{1/2} \Q_r' \wh\Q_r) \D_r^{-1/2}.
\]
Using the limit in (\ref{eqn: limit_CC}), we have
\begin{equation*}
\begin{split}
    &\hspace{5mm}
    \frac{1}{T}\sum_{t=1}^T \omega_r^{-1}\Z_r^{1/2}\Q_r'\Q_r\F_{Z,t}\Q_c'\Q_c\F_{Z,t}' \Z_r^{-1/2} =
    \frac{1}{\omega_r T}\sum_{t=1}^T \A_r'\A_r\F_{t}\A_c'\A_c\F_{t}'(\A_r'\A_r)(\A_r'\A_r)^{-1} \\
    &\xrightarrow{p}
    \tr(\A_c'\A_c)\cdot \omega_r^{-1} (\A_r'\A_r) (\A_r'\A_r) (\A_r'\A_r)^{-1}
    =
    \tr(\A_c'\A_c)\cdot \omega_r^{-1} \Z_r^{1/2}
    \bSigma_{A,r} \Z_r^{1/2} .
\end{split}
\end{equation*}
By (\ref{eqn: rate_R_res}), Assumption (L2) and eigenvector perturbation theory, there exists a unique eigenvector matrix $\bGamma_r^\ast$ of $\tr(\A_c'\A_c)\cdot \omega_r^{-1} \Z_r^{1/2} \bSigma_{A,r} \Z_r^{1/2}$ such that $\|\bGamma_r^\ast - \bGamma_r\| = o_P(1)$. Thus,
\[
\wh\Q_r' \Q_r = (\tr(\A_c'\A_c))^{-1/2} \D_r^{1/2} \bGamma_r' \Z_r^{-1/2} \xrightarrow{p}
(\tr(\A_c'\A_c))^{-1/2} \D_r^{1/2} (\bGamma_r^\ast)' \Z_r^{-1/2} ,
\]
and hence using (\ref{eqn: limit_CC}) again, we obtain
\begin{equation*}
\begin{split}
    \H_r &= \wh\D_r^{-1} \wh\Q_r' \Big(\frac{1}{T} \sum_{t=1}^T \A_r \F_{t} \A_c' \A_c \F_{t}' \A_r'\Big) \Q_r \bSigma_{A,r}^{-1}
    \xrightarrow{p}
    \D_r^{-1} \wh\Q_r' \Big(\tr(\A_c'\A_c)\cdot \A_r\A_r'\Big) \Q_r \bSigma_{A,r}^{-1} \\
    &=
    \tr(\A_c'\A_c)\cdot \D_r^{-1} \wh\Q_r' \Q_r\Z_r^{1/2}\Z_r^{1/2}\Q_r'\Q_r \bSigma_{A,r}^{-1} \xrightarrow{p}
    (\tr(\A_c'\A_c))^{1/2} \cdot \D_r^{-1/2} (\bGamma_r^\ast)' \Z_r^{1/2} ,
\end{split}
\end{equation*}
which completes the proof of Lemma \ref{lemma: limit_H}. $\square$

\begin{lemma}\label{lemma: corollary_asymp_loading}
Under the assumptions in Theorem \ref{thm: consistency_F_C}, for any $j\in[p], l\in [q]$,
\begin{align}
    & \|\wh\Q_{r,j\cdot}-\H_r \Q_{r,j\cdot}\|^2 =
    O_P\big(T^{-1} p^{\delta_{r,1} - 2\delta_{r,k_r}} q^{1-2\delta_{c,1}} +  p^{-3\delta_{r,k_r}} q^{2-2\delta_{c,1}} \big) ,
    \label{eqn: corollary_asymp_loading1} \\
    & \|\wh\Q_{c,l\cdot}-\H_c \Q_{c,l\cdot}\|^2 =
    O_P\big(T^{-1} q^{\delta_{c,1} - 2\delta_{c,k_c}} p^{1-2\delta_{r,1}} +  q^{-3\delta_{c,k_c}} p^{2-2\delta_{r,1}} \big) .
    \label{eqn: corollary_asymp_loading2}
\end{align}
\end{lemma}

\textbf{\textit{Proof of Lemma \ref{lemma: corollary_asymp_loading}.}}
For $\wh\Q_{r,j\cdot}-\H_r \Q_{r,j\cdot}$, first consider the case when $Tq = o(p^{ \delta_{r,k_r} + \delta_{r,1}} )$. From (\ref{eqn: asymp_loading}) in the proof of Theorem \ref{thm: asymp_loading}, we have
\begin{equation}\label{eqn: corollary_asymp_loading_step1}
\|\wh\Q_{r,j\cdot}-\H_r\Q_{r,j\cdot}\|^2 = O_P\big( (T p^{2\delta_{r,k_r} - \delta_{r,1}} q^{2\delta_{c,1} -1})^{-1} \big) .
\end{equation}

Now suppose $Tq = o(p^{ \delta_{r,k_r} + \delta_{r,1}})$ fails to hold. From the decomposition of $\wh\Q_{r,j\cdot}-\H_r\Q_{r,j\cdot}$ in (\ref{eqn: Q_asymp_decomp}), the leading term among the expressions will be $\cI_3$. It has rate
\[
\|\cI_3\|^2 = O_P\big(T^{-1} p^{-2\delta_{r,k_r}} q^{1-2\delta_{c,1}} + p^{-3\delta_{r,k_r}} q^{2-2\delta_{c,1}} \big) = O_P\big( p^{-3\delta_{r,k_r}} q^{2-2\delta_{c,1}} \big) ,
\]
where the last equality used the fact that $Tq = o\big(p^{ \delta_{r,k_r} + \delta_{r,1}}\big)$ does not hold. Thus we have
\begin{equation}\label{eqn: corollary_asymp_loading_step2}
\|\wh\Q_{r,j\cdot}- \H_r \Q_{r,j\cdot}\|^2 = O_P\big( p^{-3\delta_{r,k_r}} q^{2-2\delta_{c,1}} \big).
\end{equation}
Combining (\ref{eqn: corollary_asymp_loading_step1}) and (\ref{eqn: corollary_asymp_loading_step2}), we arrive at the statement of (\ref{eqn: corollary_asymp_loading1}). The proof for (\ref{eqn: corollary_asymp_loading2}) is similar and omitted here, which ends of the proof of Lemma \ref{lemma: corollary_asymp_loading}. $\square$

\bibliographystyle{apalike}
\bibliography{ref}

\end{document}